\documentclass[12pt,review]{article}
\usepackage{amsfonts}
\usepackage{amsmath}
\usepackage{enumerate}
\usepackage{mathrsfs}
\usepackage{amssymb}
\usepackage{amsthm}   
\usepackage{empheq}   
\usepackage{titlesec} 
\allowdisplaybreaks \textwidth 145mm \textheight 220mm \oddsidemargin 3mm
\usepackage{bm}
\usepackage{bookmark}
\usepackage{arydshln}
\usepackage{lineno,hyperref}
\usepackage{cite}
\modulolinenumbers[5]
\usepackage{cancel}
\usepackage{float}
\usepackage{graphicx}

\newtheorem{thm}{Theorem}[section]
\newtheorem{lem}{Lemma}[section]

\newtheorem{prop}{Proposition}[section]
\usepackage{authblk}
\usepackage{tikz}

\AtEndDocument{%
  \par
  \medskip
  \begin{tabular}{@{}l@{}}%
    \textsc{$^*$Department of Mathematics, National University of Singapore}\\
    \textit{E-mail address}: \texttt{$^1$matax@nus.edu.sg}\\
      \textsc{$^\dag$School of Mathematical Sciences, Fudan University}\\
    \textit{E-mail address}: \texttt{$^2$hy\_chen15@fudan.edu.cn}\\
    \textsc{$^\ddagger$Department of Mathematics, Hangzhou Normal University}\\
    \textit{E-mail address}: \texttt{$^3$yins@hznu.edu.cn}\\
  \end{tabular}}

\title{The title}\title{Low regularity ill-posedness for elastic waves driven by shock formation}

\author{Xinliang An$^*$$^1$}\author{Haoyang Chen$^{\dag}$$^2$}\author{Silu Yin$^{\ddagger}$$^3$}

\affil{}

\date{}

\begin{document}
\maketitle

\begin{abstract}
In this paper, we construct counterexamples to the local existence of low-regularity solutions to elastic wave equations in three spatial dimensions (3D). Inspired by the recent works of Christodoulou, we generalize Lindblad's classic results on the scalar wave equation by showing that the Cauchy problem for 3D elastic waves, a physical system with multiple wave-speeds, is ill-posed in $H^3(\mathbb{R}^3)$. We further prove that the ill-posedness is caused by instantaneous shock formation, which is characterized by the vanishing of the inverse foliation density. The main difficulties of the 3D case come from the multiple wave-speeds and its associated non-strict hyperbolicity. We obtain the desired results by designing and combining a geometric approach and an algebraic approach, equipped with detailed studies and calculations of the structures and coefficients of the corresponding non-strictly hyperbolic system. Moreover, the ill-posedness we depict also applies to 2D elastic waves, which corresponds to a strictly hyperbolic case.
\end{abstract}


\setcounter{section}{0}
\numberwithin{equation}{section}

\section{Introduction}
We study the low regularity ill-posedness problem for elastic waves in three spatial dimensions. For homogeneous isotropic hyperelastic materials, the motion of displacement $U=(U^1,U^2,U^3)$ satisfies a \underline{quasilinear} wave system with \underline{multiple} wave-speeds:
\begin{equation} \label{wave}
    \partial_t^2 U-c_2^2\Delta U-(c_1^2-c_2^2)\nabla(\nabla\cdot U)=G(\nabla U,\nabla^2U),
\end{equation}
where $c_1,c_2$ are two constants satisfying $c_1>c_2>0$. The precise form of $G(\nabla U,\nabla^2U)$ will be discussed in Section \ref{sec1.1}.

The study of system \eqref{wave} was pioneered by John. For the Cauchy problem of the three dimensional elastic wave equations with (smooth) small initial data, he proved that the singularities could arise in the radially symmetric case \cite{john84}. Without symmetry assumptions, the existence of almost-global solutions was later obtained in \cite{john88} by John and \cite{klainerman} by Klainerman-Sideris. In the case that the nonlinearities satisfy the null conditions, the solutions to the Cauchy problem of \eqref{wave} with small initial data exist globally. See Agemi \cite{Agemi} and Sideris \cite{Sideris96,Sideris00}.

In this paper, we focus on the low regularity solutions to elastic wave equations. The main conclusion we obtain is that \emph{The Cauchy problems of the 3D elastic wave equations are ill-posed in $H^3(\mathbb{R}^3)$.} We further show that the mechanism behind ill-posedness is the formation of shock, which is characterized by the vanishing of the inverse foliation density.

Our research is motivated by a series of classic works on the scalar wave equations. With the aid of planar symmetry, Lindblad gave sharp counterexamples to the local existence of low regularity solutions to semilinear and certain quasilinear wave equations in \cite{lindblad93,Lind96,Lind98}. With no symmetry assumption, the first results presenting singularity formation  for solutions to quasilinear wave equations in more than one spatial dimension were due to Alinhac \cite{Alinhac99,Alinhac99II,Alinhac01}.  He employed a Nash-Moser iteration scheme. But this approach is not capable of revealing information beyond the first blow-up point. Hence Alinhac imposed a non-degeneracy condition for initial data to close his arguments. In \cite{christodoulou10}, a breakthrough was made by Christodoulou, where he offered a detailed understanding and a complete description of shock formation in 3D without imposing
 non-degeneracy assumptions. This remarkable work was extended to a large class of equations and data settings, see \cite{speckbk,Speck16,speck20,Speck18,Speck-luk,miao}. Granowski made an important observation in his thesis, for a scalar wave equation, he connected Lindblad's counterexample of local existence in \cite{Ross} to the result of shock formation via the approach of Speck-Holzegel-Luk-Wong in \cite{Speck16}. He found that the low regularity ill-posedness of the quasillinear wave equation is driven by a shock formation. Moreover, he showed that this phenomenon is stable under a perturbation out of planar symmetry.

 Compared with a single wave equation, fewer results are known for quasilinear wave systems. In \cite{lindblad17}, Ettinger-Lindblad studied Einstein's equations and constructed a sharp counterexample for local well-posedness of Einstein vacuum equations in wave coordinates. In \cite{christodoulou}, Christodoulou-Perez studied the propagation of electromagnetic waves in nonlinear crystals. Under planar symmetry, these electromagnetic waves satisfy a first-order genuinely nonlinear and strictly hyperbolic system. By revisiting and further extending John's work \cite{john74}, Christodoulou-Perez gave a more detailed description of the behaviour of solutions at the shock-formation time. Speck \cite{Speck18} proved a stable shock formation result for a class of quasilinear wave systems with multiple wave-speeds under suitable assumptions on the nonlinearities.

In this paper, we study another physical system, the elastic waves. And we prove that the Cauchy problems for 3D elastic waves are ill-posed in $H^3(\mathbb{R}^3)$. We further show that the ill-posedness of elastic waves is also driven by shock formation. Our result extends Lindblad's classic results on scalar wave equations to a physical wave system and generalizes Christodoulou-Perez's work to the non-strictly hyperbolic case.

\subsection{Statements of equations and main theorems}\label{sec1.1}
 The motion of an elastic body in 3D is described by time-dependent orientation-preserving diffeomorphisms, written as $\Psi:\mathbb{R}^3\times\mathbb{R}\to\mathbb{R}^3$, with $\Psi=\Psi(Y,t)$  satisfying $\Psi(Y,0)=Y=(Y^1,Y^2,Y^3)$.
The deformation gradient is defined as $F=\nabla\Psi$  with  $$F^i_{j}:=\frac{\partial \Psi^i}{\partial Y^j}.$$
For a homogeneous isotropic hyperelastic material, the stored energy $W$ depends only on the principal invariants of $FF^T$ (the Cauchy-Green strain tensor). The equations of motion could be derived by employing the least action principle to $\mathcal{L}$:
\begin{align*}
  \mathcal{L}:=\int\int_{\mathbb{R}^2}\Big(\frac12|\partial_t\Psi|^2-W\Big)dYdt.
\end{align*}
Then the Euler-Lagrange equations take the form:
\begin{equation}\label{el}
  \frac{\partial^2\Psi^i}{\partial t^2 }-\frac{\partial}{\partial{Y^l}}\frac{\partial W}{\partial F^i_{l}}=0.
\end{equation}
 Let $U:=\Psi-Y$ be the displacement and $G:=\nabla U=F-I$ be the displacement gradient. Hence $FF^T=I+G+G^T+GG^T$. We thus have the stored energy $W$ being a functional of $C=G+G^T+GG^T$. We further rewrite it as $W=\hat{W}(k_1,k_2,k_3)$, where $k_1,k_2,k_3$ are the principal invariants of $C=G+G^T+GG^T$.
 In detail,
 \begin{equation*}
   \begin{split}
     k_1&=\text{tr} C,\\
     k_2&=\frac12\{(\text{tr} C)^2-\text{tr} C^2\},\\
     k_3&=\frac16\{(\text{tr} C)^3-3(\text{tr} C)(\text{tr} C^2)+2\text{tr} C^3\}.
   \end{split}
 \end{equation*}
  By Taylor expansions up to terms of higher order, we have
\begin{equation}\label{taylor}
  \hat{W}=\gamma_0+\gamma_1k_1+\frac12\gamma_{11}k_1^2+\gamma_2 k_2+\frac16\gamma_{111}k_1^3+\gamma_{12}k_1k_2+\gamma_3k_3+\cdots
\end{equation}
as a functional of $\nabla U$.
Here, $\gamma_0,\gamma_1,\gamma_{11}$, etc., are constant coefficients standing for certain partial derivatives of $\hat{W}$ at $k_i=0$ ($i=1,2,3$). Following \cite{Sideris96}, we impose the stress-free reference state condition, i.e., $\gamma_1=0$.\footnote{Tahvildar-Zadeh in \cite{shadi} studied a case where $\gamma_1\neq 0$ and proved the small data global existence under certain null condition assumptions.} Also note that positive properties of Lam\'e constants imply $\gamma_2<0$ and $2\gamma_{11}+\gamma_2>0$.\footnote{More details are referred to \cite{Agemi,Sideris96}.} Applying \eqref{taylor} into \eqref{el}, we arrive at \eqref{wave}, where
$$c_1^2=4\gamma_{11},\quad c_2^2=-2\gamma_2\quad \text{and}\quad c_1^2>c_2^2>0.$$

Keeping only the quadratic nonlinear terms of \eqref{wave}, we derive a quasilinear wave system with multiple wave-speeds:
\begin{equation}\label{y1}
  \left\{\begin{array}{ll}
    \partial_t^2 U-c_2^2\Delta U-(c_1^2-c_2^2)\nabla(\nabla\cdot U)=N(\nabla U,\nabla^2U),\\
     U(Y,0)=U_0(Y),\ U_t(Y,0)=U_1(Y)
  \end{array}\right.
\end{equation}
with $N(\nabla U,\nabla^2U)$ composed of $\nabla U\nabla^2U$-form\footnote{This means that neither $(\nabla U)^2$ nor $(\nabla^2U)^2$ is included.}.
As Agemi showed in \cite{Agemi}, the quadratic nonlinear terms could be expressed as:
\begin{align*}
  \begin{split}
   &N(\nabla U,\nabla^2U)\\
   =&\sigma_0\nabla(\text{div}\ U)^2+\sigma_1\big(\nabla|\text{curl}\ U|^2-2\text{curl}\ (\text{div}\  U\text{curl}\  U)\big)+Q(U,\nabla U),
  \end{split}
\end{align*}
where $Q$ is a summation of null forms:
\begin{align}\label{17}
  \begin{split}
    Q^i=&\sigma_2\sum_{j,k=1}^3\{2Q_{jk}(\partial_jU^i,U^k)+Q_{jk}(\partial_kU^j,U^i)+Q_{ij}(\partial_jU^k,U^k)\}\\
    &+\sigma_3\sum_{j,k=1}^3\{2Q_{ij}(\partial_kU^k,U^j)+Q_{ji}(\partial_kU^j,U^k)\}\\
    &+\sigma_4\sum_{j,k=1}^3\{Q_{ik}(U^k,\partial_jU^j)-Q_{jk}(U^j,\partial_iU^k)\},
  \end{split}
\end{align}
with $Q_{ij}(f,g)=\partial_if\partial_jg-\partial_jf\partial_ig$.
The constants $\sigma_i\ (i=0,1,2,3,4)$ satisfy
\begin{equation*}
\begin{array}{lll}
  \sigma_0=2(2\gamma_{111}+3\gamma_{11}), & \sigma_1=2(\gamma_{11}-\gamma_{12}), & \sigma_2=2(\gamma_2-\gamma_3), \\
  \sigma_3=2\gamma_3, &\sigma_4=4(\gamma_{11}-2\gamma_{12}). &
\end{array}
\end{equation*}
In this paper, we will study elastic wave equations \eqref{y1} for the general case
\begin{equation}\label{sigma}
  \sigma_0\sigma_1\neq 0.
\end{equation}

\noindent {\bf Remark:} In this paper, with the initial data we construct, the second derivatives of the displacement would tend to be infinite instantaneously, while the first derivatives of the displacement would remain small.\footnote{ This is because the smallness of the first derivatives can be transported along the characteristics.} The cubic and higher-order terms take form of $(\nabla U)^\alpha\nabla U\nabla^2U$ with $\alpha\geq1$, which are negligible compared with terms in $N(\nabla U, \nabla^2 U)$ of the form $\nabla U \nabla^2 U$. So we omit the cubic and higher-order nonlinear terms in \eqref{wave}.\\

We are now ready to state our main result:
\begin{thm} \label{3D}
The Cauchy problems of the 3D elastic wave equations \eqref{y1} are ill-posed in $H^3(\mathbb{R}^3)$ in the following sense: We construct a family of compactly supported  smooth initial data $(U_0^{(\eta)}, U_1^{(\eta)})$ with $${\|U_0^{(\eta)}\|}_{{\dot{H}}^3(\mathbb{R}^3)}+{\|U_1^{(\eta)}\|}_{{\dot{H}}^2(\mathbb{R}^3)} \to 0,\quad \text{as} \quad \eta\to0.$$ Let $T_\eta^*$ be the largest time such that \eqref{y1} $($with a general condition \eqref{sigma}$)$ has a solution $U_\eta\in C^\infty(\mathbb{R}^3\times[0,T_\eta^*) )$. As $\eta \to 0$, we have $T_\eta^* \to 0$.

Moreover, in a spatial region $\Omega_{T_\eta^*}$ the $H^2$ norm of the solution to elastic waves \eqref{y1}  blows up at shock formation time $T^*_{\eta}$:
\begin{equation*}
  \|U_\eta(\cdot,T_\eta^*)\|_{H^2(\Omega_{T_\eta^*})}=+\infty.
\end{equation*}
\end{thm}

We construct the solution of \eqref{y1} by choosing modified
 ``Lindblad-type" initial data. The above result can be extended to elastic waves in 2D as well. There we apply our method to a strictly hyperbolic system. Here is the
   2D version counterpart.
\begin{thm}\label{2D}
The Cauchy problems of 2D elastic wave equations are ill-posed in $H^\frac52(\mathbb{R}^2)$ in the following sense: We construct a family of compactly supported smooth initial data $(U_0^{(\eta)}, U_1^{(\eta)})$ with $${\|U_0^{(\eta)}\|}_{{\dot{H}}^\frac52(\mathbb{R}^2)}+{\|U_1^{(\eta)}\|}_{{\dot{H}}^\frac32(\mathbb{R}^2)} \to 0,\quad \text{as} \quad \eta\to0.$$ Let $T_\eta^*$ be the largest time such that \eqref{y1} has a solution $U_\eta\in C^\infty(\mathbb{R}^2\times[0,T_\eta^*))$. As $\eta \to 0$, we have $T_\eta^* \to 0$.

Moreover, in a spatial region $\Omega_{T_\eta^*}$ the $H^2$ norm of the solution to elastic waves \eqref{y1} blows up at shock formation time $T^*_{\eta}$:
\begin{equation*}
  \|U_\eta(\cdot,T_\eta^*)\|_{H^2(\Omega_{T_\eta^*})}=+\infty.
\end{equation*}
\end{thm}

\subsection{Difficulties and Strategies}

Though our basic strategy to construct ill-posedness is inspired by the aforementioned works, we encounter some new difficulties.

\begin{enumerate}[i)]
  \item {\it Coupling of top order terms in elastic wave equations}

  Let's start with a quasilinear wave equation with a single speed:
  \begin{equation}\label{s}
  \Box \varphi=\partial \varphi\partial^2\varphi.
  \end{equation}
     By taking a rectangular derivative $\partial_l$ $(l=0,1,2,3,)$ of \eqref{s}, we deduce:
      \begin{equation}\label{ss}
       (g^{-1})^{\alpha\beta}(\psi)\partial_\alpha\partial_\beta \psi_l=N(\partial \psi, \partial \psi_l),
      \end{equation}
      where
       $$\psi=(\psi_0,\psi_1,\psi_2,\psi_3):=(\partial_t\varphi,\partial_1\varphi,\partial_2\varphi,\partial_3\varphi).$$
       The geometric approach introduced by Christodoulou \cite{christodoulou10} could be used to prove shock formation for \eqref{ss}. However, this approach is not applicable to the elastic waves with multiple speeds. For a slightly simplified model, Speck \cite{Speck18} studied the following system of two wave equations with multiple speeds:
       \begin{equation}\label{jared1}
         \Box_{g(\Upsilon)}\Upsilon\thicksim Q^g(\partial\Upsilon,\partial\Upsilon)+\sum_{\gamma,\nu=0,1}\Big[N(\partial^\gamma q,\partial^\nu\Upsilon)+N(\partial^\gamma q,\partial^\nu q)\big],
       \end{equation}
       \begin{equation}\label{jared2}
         (h^{-1})^{\alpha\beta}(\Upsilon,\partial q)\partial_\alpha\partial_\beta q\thicksim Q^g(\partial\Upsilon,\partial\Upsilon)+\sum_{\gamma,\nu=0,1}\Big[N(\partial^\gamma q,\partial^\nu\Upsilon)+N(\partial^\gamma q,\partial^\nu q)\Big],
       \end{equation}
     where $\Box_{g(\Upsilon)}$ is the covariant wave operator\footnote{$\Box_{g}f:=\frac1{\sqrt{|\det g|}}\partial_\alpha\big(\sqrt{|\det g|}(g^{-1})^{\alpha\beta}\partial_\beta f\big).$} with respect to $g(\Upsilon)$ and $Q^g$ is the standard null form associated to $g$:
       $$Q^g(\partial\Upsilon,\partial\Upsilon):=(g^{-1})^{\alpha\beta}(\Upsilon)\partial_\alpha\Upsilon\partial_\beta\Upsilon.$$
        Note that no $\partial^2q$ term is included on the right hand side of \eqref{jared1}\footnote{Considered as the fast wave equation.}, and there is no $\partial^2\Upsilon$ term coupled in equation \eqref{jared2}. Employing the geometric method as in \cite{christodoulou10}, Speck proved shock formation of the ``fast wave $\Upsilon$". For the more general quasilinear wave systems in 3D,
        as for
        elastic waves, the $m$ unknowns $\bar{\varphi}=(\bar{\varphi}^1,\cdots,\bar{\varphi}^m)$ satisfy:
       \begin{equation*}
         \Box_i\bar{\varphi}^i=\sum_{j,k}N^i_{jk}(\partial\bar{\varphi}^j,\partial^2\bar{\varphi}^k).
       \end{equation*}
       Take a rectangular derivative, we have
       \begin{equation*}
         \bar{\psi}^i=(\bar{\psi}^i_0,\bar{\psi}^i_1,\bar{\psi}^i_2,\bar{\psi}^i_3):=(\partial_t\bar{\varphi}^i,\partial_1\bar{\varphi}^i,\partial_2\bar{\varphi}^i,\partial_3\bar{\varphi}^i)
       \end{equation*}
       verifying the system of
       \begin{equation}\label{new}
          (g^{-1})^{\alpha\beta}(\bar{\psi})\partial_\alpha\partial_\beta\bar{\psi}^i_l=\sum_{j,k}N^i_{jk}(\partial\bar{\psi}^j_l,\partial\bar{\psi}^k)+\underbrace{\sum_{j,k\atop k\neq i}N^i_{jk}(\bar{\psi}^j,\partial^2\bar{\psi}^k_l)}_{\text{new terms}},
       \end{equation}
       where $1\leq i\leq m$ and $0\leq l\leq 3$. Compared with \eqref{jared1}, there are
        $\partial^2\bar{\psi}^k$ $(k\neq i)$ terms coupled in every $i^{\text{th}}$ equation.
       Because of the new terms in \eqref{new}, there is a loss of derivatives and the geometric approaches in the aforementioned works are not applicable here.
      To study the ill-posedness problem, under the assumption of planar symmetry we use an alternative algebraic approach. The elastic waves $(u^1,u^2,u^3)$ under planar symmetry obey
\begin{align} \label{main}
    \left\{\begin{array}{lll}
    \partial_t^2u^1-c_1^2\partial_x^2u^1=\sigma_0\partial_x(\partial_xu^1)^2+\sigma_1\partial_x(\partial_xu^2)^2+\sigma_1\partial_x(\partial_xu^3)^2,\\
    \partial_t^2u^2-c_2^2\partial_x^2u^2=2\sigma_1\partial_x(\partial_xu^1\partial_xu^2),\\
    \partial_t^2u^3-c_2^2\partial_x^2u^3=2\sigma_1\partial_x(\partial_xu^1\partial_xu^3),
    \end{array}\right.
\end{align}
where $c_1>c_2>0$. Note that we have $\partial_x^2u^2$ and $\partial_x^2u^3$ in the first equation of \eqref{main}, $\partial_x^2u^1$ in the second equation of \eqref{main}, and $\partial_x^2u^1$ in the third equation of \eqref{main}.
       We could rewrite \eqref{main} as a $6\times 6$ first-order hyperbolic system:
      \begin{equation}\label{1order}
        \partial_t\Phi+A(\Phi)\partial_x\Phi=0,
      \end{equation}
      where $\Phi=(\phi_1,\phi_2,\phi_3,\phi_4,\phi_5,\phi_6)^T=(\partial_xu^1,\partial_xu^2,\partial_xu^3,\partial_tu^1,\partial_tu^2,\partial_tu^3)^T$. We then study the ill-posedness problem by exploring the algebraic structures of \eqref{1order} and by combining the geometric method as in \cite{christodoulou}.
  \item {\it Invalidity of Riemann invariants for a $6\times 6$ system}

  For waves with different traveling speeds, we want to prove shock formation for the fastest one. A key ingredient is to understand the interactions of different families of characteristics. Naturally the first and foremost step is to find a proper method to trace these characteristics. For a single wave equation (or a $2\times 2$ first-order hyperbolic system), there exist two Riemann invariants,  which could help. These two Riemann invariants verify two transport equations and can be constructed explicitly with unknown functions. The evolution of other geometric quantities along the characteristics can then be described accordingly. However, for a larger system, such as our $6\times 6$ system, to find proper Riemann invariants with explicit formula is usually impossible. Here we adopt a \underline{different} approach. We appeal to John's classic method \cite{john74}, i.e., decomposition of waves. We compute the eigenvalues of the coefficient matrix $(A(\Phi))_{6\times6}$:
  \begin{align*}
  \begin{split}
    \lambda_6(\Phi)<\lambda_5(\Phi)\leq\lambda_4(\Phi)<\lambda_3(\Phi)\leq\lambda_2(\Phi)<\lambda_1(\Phi).
  \end{split}
\end{align*}
   We give its left eigenvectors $l^i(\Phi)$ and right eigenvectors $r_i(\Phi)$, $i=1,\cdots,6$ and require:
       $$l^i(\Phi)r_j(\Phi)=\delta^i_j.$$
     Then \eqref{1order} can be rewritten as a \underline{diagonalized} system of Riccati-type:
      \begin{equation}\label{w}
          (\partial_t+\lambda_i\partial_x)w^i=-c_{ii}^i\cdot(w^i)^2+\Big(\sum_{m\neq i}(-c_{im}^i+\gamma_{im}^i)w^m\Big)w^i+\sum_{m\neq i,k\neq i\atop m\neq k}\gamma_{km}^iw^kw^m,
      \end{equation}
     where $i=1,\cdots,6$ and
      $$w^i:=l^i(\Phi)\partial_x\Phi.$$
      Here $c_{im}^i$, $\gamma_{km}^i$  are coefficients depending on unknowns.

  \item {\it Non-strict hyperbolicity}

   Equation \eqref{w} is the starting point of our further exploration.  For the aforementioned system \eqref{1order}, for $\Phi$ being a small perturbation around zero, we have that two pairs of characteristic speeds are almost the same
      \begin{equation*}
   \lambda_2(\Phi)\thickapprox\lambda_3(\Phi),\ \lambda_4(\Phi)\thickapprox\lambda_5(\Phi).
      \end{equation*}
This means two pairs of characteristic strips could overlap for a long time. Our strategy to overcome this non-strict hyperbolicity is to consider the wave propagations in four characteristic strips:
 \begin{equation*}
\{\mathcal{R}_1,\mathcal{R}_2\cup\mathcal{R}_3,\mathcal{R}_4\cup\mathcal{R}_5,\mathcal{R}_6\}.
  \end{equation*}
   These four strips will be completely separated when $t>t_0^{(\eta)}$, where $t_0^{(\eta)}$ can be precisely calculated. We then study the inverse foliation density $\rho_1$ corresponding to the fast wave. We show that $\rho_1$ vanishes (shock forms) in the strip $\mathcal{R}_1$ at some time $T_\eta^*$ after $t=t_0^{(\eta)}$. And $\rho_1$ does not vanish in the overlapped characteristic strips. See the following picture:
     
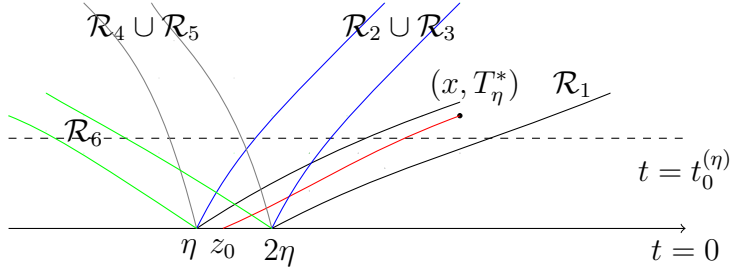
\begin{center}
\begin{tikzpicture}
\draw[->](0,0)--(9,0)node[right,below]{$t=0$};
\draw[dashed](0,1.2)--(9,1.2)node[right,below]{$t=t_0^{(\eta)}$};
\node [below]at(3.6,0){$2\eta$};
\node [below]at(2.4,0){$\eta$};
\filldraw [black] (3.5,0) circle [radius=0.01pt]
(5,0.8) circle [radius=0.01pt]
(6,1) circle [radius=0.01pt]
(8,1.8) circle [radius=0.01pt];
\draw (3.5,0)..controls (5,0.8)and(6,1)..(8,1.8);

\filldraw [black] (2.5,0) circle [radius=0.01pt]
(4,1) circle [radius=0.01pt]
(5.5,1.5) circle [radius=0.01pt]
(6,1.68) circle [radius=0.01pt];
\draw (2.5,0)..controls (4,1) and (5.5,1.5)..(6,1.68);
\node [above] at(7.5,1.6){$\mathcal{R}_1$};

\filldraw [black] (2.85,0) circle [radius=0.01pt]
(4,0.5) circle [radius=0.01pt]
(4.65,1) circle [radius=0.01pt]
(6,1.5) circle [radius=0.8pt];
\draw [color=red](2.85,0)..controls (4,0.5) and (4.65,1)..(6,1.5);
\node [below]at(2.85,0){$z_0$};

\node [above]at(6.2,1.5){$(x,T_\eta^*)$};

\filldraw [green] (3.5,0) circle [radius=0.01pt]
(4,1) circle [radius=0.01pt]
(4.5,1.5) circle [radius=0.01pt]
(6,3) circle [radius=0.01pt];
\draw [color=blue](3.5,0)..controls (4,1) and (4.5,1.5)..(6,3);

\filldraw [green] (2.5,0) circle [radius=0.01pt]
(3,1) circle [radius=0.01pt]
(3.5,1.5) circle [radius=0.01pt]
(5,3) circle [radius=0.01pt];
\draw [color=blue](2.5,0)..controls (3,1) and (3.5,1.5)..(5,3);
\node [below] at(5.2,3){$\mathcal{R}_2\cup\mathcal{R}_3$};

\filldraw [gray] (3.5,0) circle [radius=0.01pt]
(3,2) circle [radius=0.01pt]
(2.2,2.5) circle [radius=0.01pt]
(1.9,3) circle [radius=0.01pt];
\draw [color=gray](3.5,0)..controls (3,2) and (2.2,2.5)..(1.9,3);
\node [below] at(1.8,3){$\mathcal{R}_4\cup\mathcal{R}_5$};

\filldraw [gray] (2.5,0) circle [radius=0.01pt]
(2,2) circle [radius=0.01pt]
(1.5,2.5) circle [radius=0.01pt]
(1,3) circle [radius=0.01pt];
\draw [color=gray] (2.5,0)..controls (2,2) and (1.5,2.5)..(1,3);

\filldraw [green] (3.5,0) circle [radius=0.01pt]
(2,1) circle [radius=0.01pt]
(1,1.5) circle [radius=0.01pt]
(0.5,1.8) circle [radius=0.01pt];
\draw [color=green](3.5,0)..controls (2,1) and (1,1.5)..(0.5,1.8);
\node [below] at(1,1.6){$\mathcal{R}_{6}$};

\filldraw [green] (2.5,0) circle [radius=0.01pt]
(1,1) circle [radius=0.01pt]
(0.5,1.3) circle [radius=0.01pt]
(0,1.5) circle [radius=0.01pt];
\draw [color=green] (2.5,0)..controls (1,1) and (0.5,1.3)..(0,1.5);
\end{tikzpicture}
\end{center}
     Subtle structures of \eqref{w} are the key in our analysis, which makes it possible to trace $\rho_i$ in the above four characteristic strips up to time $T^*_{\eta}$. We list some crucial structures here:
      \begin{equation}\label{wsketch}
 \left\{ \begin{split}
    (\partial_t+\lambda_1\partial_x)w^1\thicksim&(w^1)^2+\cancel{w^2w^3}+\cancel{w^4w^5}+\cdots,\\
    (\partial_t+\lambda_2\partial_x)w^2\thicksim&\cancel{(w^2)^2}+(\lambda_2-\lambda_3)w^2w^3+(\lambda_4-\lambda_5)w^4w^5+\cdots,\\
    (\partial_t+\lambda_3\partial_x)w^3\thicksim&(w^3)^2+\cancel{w^2w^3}+\cancel{w^4w^5}+\cdots,\\
    (\partial_t+\lambda_4\partial_x)w^4\thicksim&(w^4)^2+\cancel{w^2w^3}+\cancel{w^4w^5}+\cdots,\\
    (\partial_t+\lambda_5\partial_x)w^5\thicksim&\cancel{(w^5)^2}+(\lambda_2-\lambda_3)w^2w^3+(\lambda_4-\lambda_5)w^4w^5+\cdots,\\
    (\partial_t+\lambda_6\partial_x)w^6\thicksim&(w^6)^2+\cancel{w^2w^3}+\cancel{w^4w^5}+\cdots.\\
  \end{split}\right.
\end{equation}
 The deleted terms with cancelation symbols in \eqref{wsketch} mean their coefficients are \underline{zero}. That is to say there is no interaction of the almost-repeated characteristic waves ($w^2$ and $w^3$, $w^4$ and $w^5$) appearing in the equations of $\{w^i\}_{i=1,3,4,6}$.

Moreover, in the equations for $w^2$ and $w^5$ we also have $\lambda_2-\lambda_3=\lambda_4-\lambda_5$ being small, and it shows that the related interactions are weak. These small coefficients allow us to use the bi-characteristic transformation to get desired bounds for non-strictly hyperbolic systems. For instance, with the equation $(\partial_t+\lambda_2\partial_x)w^2$, to estimate $w^2$ we will need to bound the integration of $(\lambda_2-\lambda_3)w^2 w^3$ within $\mathcal{R}_2\cup\mathcal{R}_3$. See the exact estimates in \eqref{y18} and \eqref{553}. There we will use bi-characteristic transformation as in \cite{christodoulou}
\begin{equation*}
    dt=\frac{\rho_2}{\lambda_3-\lambda_2}dy_2+\frac{\rho_3}{\lambda_2-\lambda_3}dy_3.
  \end{equation*}
This transformation could be singular owing to the \underline{non-strict} hyperbolicity, since $\lambda_2(0)-\lambda_3(0)=0$. Nevertheless, the additional small coefficient $\lambda_2-\lambda_3$ in front of $w^2 w^3$ save us!

\item {\it Meaning of blow-ups}

 Notice that \eqref{wsketch} are a set of Riccati-type equations. The key structures mentioned above ensure that singularities (blow-ups) will form in finite time. A new ingredient of this paper is to study the meaning of the blow-ups and to study the ill-posedness mechanism of the first-order non-strictly hyperbolic system. Based on the decomposition of waves in an algebraic manner, we trace the evolution of inverse foliation density $\rho_i$ ($i=1,\cdots,6$) coming from geometry. Here $\rho_i$ depicts the density of nearby characteristics in $\mathcal{R}_i$ and it verifies
 \begin{equation*}
   (\partial_t+\lambda_i\partial_x)\rho_i=c_{ii}^iv^i+\Big(\sum_{m\neq i}c_{im}^iw^m\Big)\rho_i,
 \end{equation*}
  where $v^i=\rho_iw^i$ for fixed $i$.
  By constructing suitable initial data, we get a positive lower bound for $\{\rho_i\}_{i=2,\cdots,6}$:
   \begin{equation*}
     \min_{i\in\{2,\cdots,6\}}\inf_{0\leq t\leq T_\eta^*} \rho_i\geq \frac{(1-\varepsilon)^2}{2},
   \end{equation*}
   for some sufficiently small $\varepsilon$. Furthermore, $\rho_1$ is proved to obey
   \begin{equation*}
   \begin{split}
     (1-\varepsilon)\Big(1-(1+\varepsilon)^3|c_{11}^1(0)|tW_0^{(\eta)}\Big)\leq\rho_1(X_1&(z_0,t),t)\\
     \leq& (1+\varepsilon)\Big(1-(1-\varepsilon)^4|c_{11}^1(0)|tW_0^{(\eta)}\Big),
     \end{split}
   \end{equation*}
   where $W_0^{(\eta)}=\max_i\sup_z|w^i_{(\eta)}(z,0)|=w^1_{(\eta)}(z_0,0)$.
   In view of the above inequalities, we conclude that $\rho_1(z_0,t)\to 0$ as $t\to T_\eta^*$ with $T_\eta^*\sim1/W_0^{(\eta)}$. A shock forms at time $T_\eta^*$ and we can show that $T_\eta^*$ is the first time when blow-up appears.

We then check the $H^2(\mathbb{R}^3)$ norm of solutions to \eqref{y1} at time $T_\eta^*$. In a suitable constructed spatial region $\Omega_{T^*_\eta}$, we have
 \begin{equation*}
 \begin{split}
 \|U(\cdot,T^*_\eta)\|_{H^2(\Omega_{T^*_\eta})}^2\gtrsim & C_\eta\int_{z_0}^{z_0^*}\frac{1}{\rho_1(z,T_\eta^*)}dz\\
 \gtrsim & C_\eta\int_{z_0}^{z_0^*}\frac{1}{(\sup_{z_1\in(z_0,z_0^*]}|\partial_{z_1}\rho_1|)(z-z_0)}dz\\
 =&+\infty.
 \end{split}
\end{equation*}
In the above inequalities, we utilize
\begin{equation*}
  \rho_1(z,T_\eta^*)=\rho_1(z,T_\eta^*)-\rho_1(z_0,T_\eta^*)=\partial_{z_1}\rho_1(z',T_\eta^*)(z-z_0)
\end{equation*}
for some $z'\in(z_0,z_0^*]$ and employ the uniform bound of $\partial_{z_1}\rho_1$ obtained in Section \ref{supbound}. This shows that the ill-posedness is driven by the shock formation: $\rho_1(z_0,t)\to 0$ as $t\to T_\eta^*.$
 Moreover, as $\eta \to 0$, we have $T_\eta^*\to 0$.
\end{enumerate}

\subsection{New ingredients}
The method we develop in this paper might be useful for studying other problems. We list some of the key points:

\begin{enumerate}
\item For elastic waves (\ref{y1}), both our lower regularity ill-posedness result and our exploration of the ill-posedness mechanism in this paper are \underline{new}. We extend the previous result on a single quasilinear wave equation to a physical wave system.

\item For elastic waves with multiple speeds, under plane symmetry we give a \underline{complete} description of the wave dynamics up to the time $T^*_{\eta}$, when the first (shock) singularity happens. For $t\leq T^*_{\eta}$, there is \underline{no} other singular point in the spacetime region. And the solution is smooth before time $T^*_{\eta}$. We summarize it in Proposition \ref{wave dynamics}.

\item Our algebraic approach of rewriting (\ref{main}) into a $6\times 6$ system and using $r_j$ in (\ref{right}) and $l^i$ in (\ref{left}) as right and left eigenvectors for wave decomposition is \underline{new}. By algebraic calculations, we find the system (\ref{main}) is \underline{non-strictly} hyperbolic, with eigenvalues of the coefficient matrix $(A(\Phi))_{6\times6}$ satisfying
  \begin{align*}
  \begin{split}
    \lambda_6(\Phi)<\lambda_5(\Phi)\leq\lambda_4(\Phi)<\lambda_3(\Phi)\leq\lambda_2(\Phi)<\lambda_1(\Phi).
  \end{split}
\end{align*}
For $i,j=1,\cdots,6$, for the left eigenvectors $l^i(\Phi)$ and right eigenvectors $r_i(\Phi)$ we construct, it holds
$$l^i(\Phi)r_j(\Phi)=\delta^i_j,$$
but there is \underline{not} requirement of $l^i(\Phi)$ and $r_i(\Phi)$ being unit. This \underline{flexibility} enables us to use bi-characteristic coordinates for our non-strictly hyperbolic system.

\item It is the \underline{first} time the subtle structures of (\ref{wsketch}) are explored.  See also the structures in \eqref{infmw} and \eqref{infmrho}.  As aforementioned, these subtle structures are really the \underline{key} to our proof. And utilizing them may lead to other applications.

\item Tracing dynamics of the $6\times 6$ non-strictly hyperbolic system with four characteristic strips $\{\mathcal{R}_1, \mathcal{R}_2\cup\mathcal{R}_3, \mathcal{R}_4\cup\mathcal{R}_5, \mathcal{R}_6\}$ is \underline{new}. Here $\mathcal{R}_2$ and $\mathcal{R}_3$ (or $\mathcal{R}_4$ and $\mathcal{R}_5$) could overlap with each other for a long time. This treatment is applied together with using the subtle structures of (\ref{main}) and (\ref{wsketch}).

For instance, $\mathcal{R}_2$ and $\mathcal{R}_3$ could overlap for a long time. But for the equation $(\partial_t+\lambda_2\partial_x)w^2$ in (\ref{wsketch}), there is a small coefficient $\lambda_2-\lambda_3$ in front of $w^2w^3$. This means that the worrisome nonlinear interaction of $w^2$ and $w^3$ inside $\mathcal{R}_2\cup\mathcal{R}_3$ is small. In addition, the corresponding coefficient of $w^2w^3$ in the equation of $(\partial_t+\lambda_3\partial_x)w^3$ is \underline{zero}. We incorporate the use of these subtle structures in our proof.

\item Because of the non-strict hyperbolicity, we need to include $\underline{S}$, i.e., estimates of the lower bound of $\{\rho_i\}_{i=2,\cdots, 6}$ as part of the bootstrap argument. This is an extension of Christodoulou-Perez \cite{christodoulou}. To get the desired results, we also employ the \underline{modified} Lindblad-type initial data by introducing a small parameter $\eta$. This allows us to prove that our constructed initial data are in $H^3(\mathbb{R}^3)$ for 3D and are in $H^{\frac52}(\mathbb{R}^2)$ for 2D and for $T^*_{\eta}$ being the first blow-up time, it holds $T^*_{\eta}\rightarrow 0$ as $\eta\rightarrow 0$.

\item For elastic waves (\ref{y1}), not only we prove $H^3$ ill-posedness  in 3D and $H^{\frac52}$ ill-posedness in 2D, but also we demonstrate the ill-posedness mechanism, and it is driven by shock formation. In addition, at the shock formation time $T^*_{\eta}$, the $H^2$ norms of the solutions in 3D and 2D are \underline{infinity}.

For elastic waves (\ref{y1}) in 3D, its critical norm is $H^{\frac52}(\mathbb{R}^3)$. And for the 2D case, its critical norm is  $H^{2}(\mathbb{R}^2)$. Our ill-posedness results are with Sobolev norms $\frac12$-derivative higher than the critical norms. Let's take a look at another physical quasilinear wave system \textit{Einstein's equations}. In $3+1$ dimensions, its critical norm is $H^{\frac32}(\mathbb{R}^3)$, but the sharp local well-posedness and ill-posedness results are all at the level of $H^2(\mathbb{R}^3)$, i.e., a $\frac12$-derivative higher than the critical norm. As an analogue, our $H^3(\mathbb{R}^3)$ ill-posedess of elastic waves (\ref{y1}) in 3D is a \underline{desired} result.

\end{enumerate}

\subsection{A single wave model}
In this section, we demonstrate the basic geometric and algebraic approaches with a single wave model\footnote{For a simpler scenario, if we set the waves $u^2$ and $u^3$ to vanish, the first equation of \eqref{main} takes the form of \eqref{single}. The constants $c_1$ and $\sigma_0$ have no essential influence on these two approaches.}:
\begin{equation} \label{single}
  \partial_t^2\varphi-\partial_x^2 \varphi=\partial_x(\partial_x \varphi)^2.
\end{equation}
This is a simplified case of the 3D quasilinear scalar wave equation:
\begin{equation*}
  \partial_t^2\varphi-\Delta \varphi=\partial(\partial \varphi)^2,\quad \partial\ \text{is a rectangular spatial derivative.}
\end{equation*}
To study the ill-posedness of a scalar wave equation, Lindblad considered in \cite{Lind98} the following quaslinear wave model
\begin{equation*}
  \partial_t^2\varphi-\Delta \varphi=D(D \varphi)^2,
\end{equation*}
where $D=\partial_x-\partial_t$. In this case, the explicit formula of solutions can be obtained by solving the initial data problem along the characteristics.
However, this is not applicable for \eqref{single} to establish explicit formula along characteristics. In order to gain some insights about solutions to \eqref{single}, we take a glance at two different methods. One is the geometric method. The shock formation can be described in a more explicit manner by studying the corresponding geometric wave equations of \eqref{single}. The other is an algebraic approach, with which we rewrite \eqref{single} as a $2\times 2$ system and we algebraically diagonalize this system by introducing and finding Riemann invariants. This leads to two Riccati-type  equations, which may have singularities formed in finite time.
\subsubsection{Geometric approach}
We first rewrite the model equation \eqref{single} into the following form,
\begin{equation*}
  (g^{-1}(\partial_x \varphi))^{\alpha\beta} \partial_\alpha\partial_\beta\varphi=0.
\end{equation*}
The ill-posedness and shock formation of the above equation have been studied in Speck\cite{Speck18}\cite{Speck16} and Granowski\cite{Ross}. For demonstration purpose, we review the ideas and methods therein.

Take $\psi=\partial_x \varphi$. We have $\psi$ satisfies a geometrically covariant wave equation:
\begin{equation} \label{singlegeo}
  \Box_{g(\psi)} \psi=\mathscr{Q}(\partial \psi, \partial \psi).
\end{equation}
According to \eqref{single}, the nonlinear terms are given by $$\mathscr{Q}(\partial \psi, \partial \psi)=-{(1+2\psi)}^{-1}\mathscr{Q}_0(\partial \psi, \partial \psi),$$
 where $\mathscr{Q}_0(\partial \psi, \partial \psi)$ is the standard null form with respect to the geometric metric $g(\psi)$.

To understand the causal structure with respect to $g(\psi)$, one can use the geometric coordinates $(t,u)$, with $u$ an eikonal function satisfying the eikonal equation
\begin{equation*}
g^{\alpha\beta}\partial_\alpha u\partial_\beta u=0.
 \end{equation*}
 The intersection of characteristics is then described by the vanishing of a geometric quantity called \emph{inverse foliation density} $\mu$,
\begin{equation*}
\mu:=-\frac{1}{ {(g^{-1})}^{\alpha\beta}(\psi)\partial_{\alpha}t\partial_{\beta}u}=\frac1{\partial_tu}.
\end{equation*}
We then define geometric frame $(L,\check{\underline{L}})$, where $\check{\underline{L}}$ and $L$ represent incoming and outgoing null directions respectively. In particular, we have $L=\partial_t$ in $(t,u)$ coordinates. The equation \eqref{singlegeo} for $\psi$ can be reformulated in this frame. And it also holds that $\mu$ satisfies a transport equation along the outgoing characteristic. In detail, equation \eqref{single} can be written as a system with the following form:
\begin{equation*}
\begin{cases}
L\check{\underline{L}}\psi=-{(1+2\psi)}^{-1} L\psi \cdot \check{\underline{L}}\psi,\\
\check{\underline{L}}L\psi=-\frac{\mu}{2}{(1+2\psi)}^{-1} {(L\psi)}^2-\frac{1}{2}{(1+2\psi)}^{-1} L\psi \cdot \check{\underline{L}}\psi,\\
L\mu=-\frac{\mu}{2}{(1+2\psi)}^{-1} L\psi-\frac{1}{2}{(1+2\psi)}^{-1} \check{\underline{L}}\psi.
\end{cases}
\end{equation*}

By properly choosing initial data, and by setting and proving a bootstrap argument for $\psi$, $L\psi$, $\check{\underline{L}}\psi$ and $\mu$, one can arrive at a description of the inverse foliation density
$$\mu \sim 1-(\sup\limits_{0\leq u\leq 1}{[f(u)]}_{-})t,$$
where $f$ is the initial data given at $t=0$ along $\Sigma_0^1:=\{(x,0)|0\leq u(x,0)\leq 1\}$. This leads to shock formation  at time $T_{(shock)} \sim 1/{(\sup\limits_{0\leq u\leq 1}{[f(u)]}_{-})}$. We refer to \cite{Speck18,Speck16} for more details.

Granowski \cite{Ross} studied ill-posedness theory of the above geometric wave equation \eqref{singlegeo}. By prescribing Lindblad-type initial data, he showed the instantaneous blow-up of the $H^2(\mathbb{R}^3)$ norm. In his proof, based on the setup of Speck\cite{Speck18}, the inverse foliation density also satisfies a transport equation:
\begin{equation*}
L\mu=\frac{1}{2}G_{LL} \check{X}\Phi,\footnote{The explicit forms of $G_{LL}$ and $\check{X}\Phi$ are given in \cite{Ross}.}
\end{equation*}
in geometric coordinates $(t,u)$. We hence obtain
\begin{equation*}
\mu(t,u)=\mu(0,u)+\frac{t}{2}G_{LL}(u) \check{X}\Phi(u).
\end{equation*}
From the calculation of $H^2(\mathbb{R}^3)$ norm, we can read its instantaneous blow-up is driven by the vanishing of $\mu$, which demonstrates the shock formation. See \cite{Ross} for more details.

\subsubsection{Algebraic approach}

We also outline another approach to study \eqref{single}.
Taking $V^{(1)}=\varphi_t$ and $V^{(2)}=\varphi_x$, equation \eqref{single} can be transformed into a first-order genuinely nonlinear\footnote{See the definition of genuinely nonlinear in sense of Lax given in Lemma \ref{GN}.} strictly hyperbolic system\footnote{Recall the definition of strictly hyperbolic in Section \ref{nstricth}.}:
\begin{equation}\label{riem}
  \left\{\begin{array}{lll}
    \partial_tV^{(1)}-(1+2V^{(2)})\partial_xV^{(2)}=0,\\
    \partial_tV^{(2)}-\partial_xV^{(1)}=0.
  \end{array}\right.
\end{equation}
By direct calculation, the eigenvalues of its coefficient matrix are obtained:
\begin{equation*}
\lambda_1(V)=-\sqrt{1+2V^{(2)}},\quad \lambda_2(V)=\sqrt{1+2V^{(2)}}.
\end{equation*}
Now we can introduce $W=(W^{(1)},W^{(2)})$ and they are defined through
\begin{equation*}
W^{(1)}=V^{(1)}-\frac{1}{3}{(1+2V^{(2)})}^{\frac{3}{2}}+\frac{1}{3},
\end{equation*}
\begin{equation*}
W^{(2)}=V^{(1)}+\frac{1}{3}{(1+2V^{(2)})}^{\frac{3}{2}}-\frac{1}{3}.
\end{equation*}
Then the eigenvalues $\lambda_1$, $\lambda_2$ can be viewed as functions depending on $W$. It is a straight forward check that $W^{(1)}$ and $W^{(2)}$ diagonalize \eqref{riem}, i.e., they satisfy the following transport equations:
\begin{equation} \label{Riem}
\begin{cases}
\partial_t W^{(1)}+\lambda_2(W) \partial_x W^{(1)}=0, \\
\partial_t W^{(2)}+\lambda_1(W) \partial_x W^{(2)}=0.
\end{cases}
\end{equation}
Then one can define two families of characteristics, which are the solutions to the following ODEs:
\begin{equation*}
\begin{cases}
\frac{dx^{(1)}(t)}{dt}=\lambda_1(W(x^{(1)}(t),t)),\\
x^{(1)}(0)=y_1,
\end{cases}
\quad
\begin{cases}
\frac{dx^{(2)}(t)}{dt}=\lambda_2(W(x^{(2)}(t),t)),\\
x^{(2)}(0)=y_2.
\end{cases}
\end{equation*}
By \eqref{Riem}, $W^{(1)}$ and $W^{(2)}$ are invariant along the characteristics $( x^{(2)}(t),t)$ and $( x^{(1)}(t),t)$, respectively. We hence call $W^{(1)}$ and $W^{(2)}$ Riemann invariants for equation \eqref{riem}. Note that for the Cauchy problem of \eqref{Riem}, $W( x,t)$ remains bounded for any given bounded initial data $W(x,0)$. But $\partial_xW^{(1)}$ or $\partial_xW^{(2)}$ could blow up. It can be deduced that $\partial_xW^{(1)}$ and $\partial_xW^{(2)}$ also satisfy Riccati-type equations. Along the second characteristic $(x^{(2)}(t),t)$, we have
\begin{equation} \label{riccati2}
\begin{split}
\frac{dW^{(1)}_x( x^{(2)}(t),t)}{dt}+\frac{\partial \lambda_2}{\partial W^{(2)}} \frac{1}{\lambda_2-\lambda_1}\frac{dW^{(2)}_x(x^{(2)}(t),t)}{dt}W^{(1)}_x( x^{(2)}(t),t)\\
=-\frac{\partial \lambda_2}{\partial W^{(1)}}[W^{(1)}_x( x^{(2)}(t),t)]^2.
\end{split}\end{equation}
And along the first characteristic $( x^{(1)}(t),t)$, $\partial_xW^{(1)}$ verities:
\begin{equation} \label{riccati1}
\begin{split}
\frac{dW^{(2)}_x(x^{(1)}(t),t)}{dt}+\frac{\partial \lambda_1}{\partial W^{(1)}} \frac{1}{\lambda_1-\lambda_2}\frac{dW^{(1)}_x(x^{(1)}(t),t)}{dt}{W^{(2)}_x( x^{(1)}(t),t)}\\
=-\frac{\partial \lambda_1}{\partial W^{(2)}}[W^{(2)}_x( x^{(1)}(t),t)]^2.
\end{split}\end{equation}
The solution of \eqref{riccati2} can be expressed as:
\begin{equation}\label{129}
\begin{split}
  \partial_xW^{(1)}(x^{(2)}(t),t)=\frac{e^{-I_t}\partial_xW^{(1)}_0(y_2)}{1+\partial_xW^{(1)}_0(y_2)e^{-I_t}\int_0^te^{-I_\tau}\frac{\partial\lambda_2}{\partial W^{(1)}}(W(x^{(2)}(\tau),\tau))d\tau},
\end{split}
\end{equation}
where $\partial_xW^{(1)}_0(y_2)=\partial_xW^{(1)}( x^{(2)}(0),0)=\partial_xW^{(1)}( y_2,0)$ and $$I_t=\int_{W^{(2)}( x^{(2)}(0),0)}^{W^{(2)}(x^{(2)}(t),t)}\frac1{\lambda_2-\lambda_1}\frac{\partial\lambda_2(W^{(1)}_0,W^{(2)})}{\partial W^{(2)}}dW^{(2)}.$$
Moreover, $e^{I_t}$ has positive lower bound and upper bound. If further assuming that there exist a uniform positive constant $C$ and a point $x_0$ such that
\begin{equation*}
 \frac{\partial\lambda_2}{\partial W^{(1)}}>C,\quad \partial_xW^{(1)}_0(x_0)<0,
\end{equation*}
then from \eqref{129}, we have $ \partial_xW^{(1)}$ blows up in finite time.

Picking Lindblad-type initial data as in \cite{Lind98}, by further analysis we can show that at the blow-up point, the corresponding inverse foliation density $\Big(\frac{\partial x^{(2)}(y_2,t)}{\partial y_2}\Big)$ also vanishes, and it renders $H^2(\mathbb{R}^3)$ norm to blow up, which is driven by shock formation.

\noindent {\bf Remark:} The classic algebraic approach on constructing explicit Riemann invariants and using them for a genuinely nonlinear strictly hyperbolic system, relies heavily on its being a $2\times2$ first-order system. However, for more complicated cases, such as our setting of a $6\times6$ system for elastic waves, it is almost impossible to construct proper Riemann invariants with explicit forms. John \cite{john74} extended this idea of constructing Riemann invariants to calculating and analyzing the decomposition of waves, which could be useful for larger hyperbolic systems.\footnote{In \cite{Zhou}, with John's formula, Zhou studied the local well-posedness of low regularity solutions for linearly degenerate hyperbolic systems in one spatial dimension.} In this paper, we will also adopt John's approach.

\subsection{Main steps in the proof}
Here we outline our proof of Theorem \ref{3D}. We take an algebraic approach first and we decompose the waves as \cite{john74,christodoulou}. The components of decomposed waves satisfy a quasilinear hyperbolic system. It is worthwhile to mention that our system is \underline{not} strictly hyperbolic: the coefficient matrix has almost-the-same eigenvalues for small perturbation of the unknowns around zero.
In the end, we will show that \eqref{y1} is ill-posed in $H^3(\mathbb{R}^3)$ and it is caused by shock formation of the fast characteristic wave.

\textbf{Step 1: Reduction to a first-order quasilinear hyperbolic system}

We first transform the equations of elastic plane waves into a first-order {\color{black}quasilinear} hyperbolic system 
\begin{align}\label{0y2}
  \partial_t\Phi+A(\Phi)\partial_x\Phi=0,
\end{align}
for $\Phi=(\phi_1,\phi_2,\phi_3,\phi_4,\phi_5,\phi_6)^T=(\partial_xu^1,\partial_xu^2,\partial_xu^3,\partial_tu^1,\partial_tu^2,\partial_tu^3)^T$.
This system is \textit{not} uniformly strictly hyperbolic in a small ball $B_{2\delta}^6(0)$ where the amplitude of $\Phi$ is bounded by a small parameter $\delta$. In particular, the eigenvalues of the coefficients matrix $(A(\Phi))_{6\times 6}$ satisfy
\begin{align*}
  \begin{split}
    \lambda_6(\Phi)<\lambda_5(\Phi)\leq\lambda_4(\Phi)<\lambda_3(\Phi)\leq\lambda_2(\Phi)<\lambda_1(\Phi),
  \end{split}
\end{align*}
where $\lambda_2(\Phi)$ and $\lambda_3(\Phi)$ or $\lambda_4(\Phi)$ and $\lambda_5(\Phi)$ could be the same for some $\Phi\in B_{2\delta}^6(0)$. Based on the eigenvalues, the left eigenvectors $\{l^i(\Phi)\}_{i=1,\cdots,6}$ and the right eigenvectors $\{r_i(\Phi)\}_{i=1,\cdots,6}$ can also be calculated, and we require their verifying the property
\begin{equation*}
  l^i(\Phi)r_j(\Phi)=\delta^i_j.
\end{equation*}

\textbf{Step 2: Decomposition of waves}

The next is to decompose the waves and to derive the equations for their components. As in \cite{christodoulou}, we use characteristic coordinates and bi-characteristic coordinates. For any $(x,t)\in \mathbb{R}\times [0,T]$, there is a unique $(z_i,s_i)\in \mathbb{R}\times [0,T]$, called characteristic coordinates, such that $$(x,t)=\big(X_i(z_i,s_i),s_i\big),$$
where the flow map $X_i(z_i,s_i)$ is defined by
\begin{align*}
  \left\{\begin{array}{ll}
  \frac{\partial}{\partial t}X_i(z,t)=\lambda_i\big(\Phi(X_i(z_i,t),t)\big),\quad t\in[0,T],\\
  X_i(z,0)=z_i.
  \end{array}\right.
\end{align*}
Given $(x,t)\in \mathbb{R}\times [0,T]$, along different characteristics, there is a unique $(y_i,y_j)\in\mathbb{R}^2$, called bi-characteristic coordinates, such that $t=t'(y_i, y_j) \mbox{ and }$
\begin{equation*}
  (x,t)=\big(X_i(y_i,t'(y_i,y_j)),t'(y_i,y_j)\big)=\big(X_j(y_j,t'(y_i,y_j)),t'(y_i,y_j)\big).
\end{equation*}
Let $\mathcal{C}_i(z_i)$ be the $i^{\text{th}}$ characteristic with propagation speed $\lambda_i$ starting at $z_i$. Denote the corresponding $i^{\text{th}}$ characteristic strip to be $\mathcal{R}_i:=\cup_{z_i\in I_0}\mathcal{C}_i(z_i)$, where $I_0$ is the support of initial data. Define $$\rho_i:=\partial_{z_i}X_i$$ and use it to describe the inverse foliation density of the $i^{\text{th}}$ characteristics. For fixed $i$, let
\begin{equation*}
  w^i:=l^i\partial_x\Phi, \quad \text{and}\quad
  v^i:=\rho_iw^i.
\end{equation*}
It can be checked that these geometric quantities satisfy
\begin{align}
  \partial_{s_i}\rho_i=&c_{ii}^iv^i+\Big(\sum_{m\neq i}c_{im}^iw^m\Big)\rho_i,\\
  \partial_{s_i}w^i=&-c_{ii}^i(w^i)^2+\Big(\sum_{m\neq i}(-c_{im}^i+\gamma_{im}^i)w^m\Big)w^i+\sum_{m\neq i,k\neq i\atop m\neq k}\gamma_{km}^iw^kw^m,\label{123}\\
  \partial_{s_i}v^i=&\Big(\sum_{m\neq i}\gamma_{im}^iw^m\Big)v^i+\sum_{m\neq i,k\neq i\atop m\neq k}\gamma_{km}^iw^kw^m\rho_i,
\end{align}
where the characteristic vectorfields $\{\partial_{s_i}\}$ are given by
$$\partial_{s_i}=\partial_t+\lambda_i\partial_x.$$
We note that the non-zero of $c_{ii}^i$ means the genuine non-linearity in \eqref{123}, while the vanishing of $c_{ii}^i$ denotes degeneration. Since $\lambda_1$ is the largest eigenvalue, we have $w^1$ being the wave of the fastest speed. Being genuinely nonlinear means $c_{11}^1\neq0$ and without loss of generality we assume $c_{11}^1<0$, and clearly \eqref{123} gives a Riccati-type equation. In \cite{christodoulou}, Christodoulou-Perez studied the shock formation of the uniformly strictly hyperbolic systems. While for elasticity a lot of attention should be paid to the structures of elastic wave equations related to the $2^{\text{nd}}$, $3^{\text{rd}}$, $4^{\text{th}}$, and $5^{\text{th}}$ characteristics, which make the system not uniformly strictly hyperbolic. We explore the structures and calculate the coefficients $\{c_{im}^i\}$ and $\gamma_{km}^i$ very carefully. We observe that some crucial coefficients are zero and some key terms have almost-zero factor in front. These structures are vital in the proof of the main theorem. See Section \ref{yf} for the details.

\textbf{Step 3: Construction of initial data}

We employ a family of Lindblad-type smooth initial data $ \{w^i_{(\eta)}(z_i,0)\}$ and modify them a bit such that they are supported in $[\eta,2\eta]$ for a given small $\eta$ and satisfy
 \begin{equation*}
 \left\{\begin{split}
  &W_0^{(\eta)}:=\max_{i=1,\cdots,6}\sup_{z_i}|w^i_{(\eta)}(z_i,0)|=w^1_{(\eta)}(z_0,0),\\
   &\max_{i=3,4}\sup_{z_i}|w^i_{(\eta)}(z_i,0)|\leq \min\Big\{\frac{(1-\varepsilon)^4|c_{11}^1(0)|}{(1+\varepsilon)^3}W_0^{(\eta)},W_0^{(\eta)}\Big\},\\
   &\sup_{z_6}|w^6(z_6,0)|\leq \frac{(1-\varepsilon)^4}{2(1+\varepsilon)^3}W_0^{(\eta)},
  \end{split}\right.
\end{equation*}
with
\begin{equation*}
   w^1_{(\eta)}(z,0):=\theta\int_\mathbb{R}\zeta_{\frac\eta{10}}(y)| \ln (z-y)|^\alpha \chi(z-y)dy, \quad 0<\alpha<\frac12,
\end{equation*}
where $\chi$ is a characteristic function, $\zeta_{\frac\eta{10}}$ is a test function and $\theta$ is a small constant to be chosen. Note that we can choose $\hat{w}(x,Y^2,Y^3)=w^1_{(\eta)}(x,0)$ such that $\hat{w}\in H^1(\mathbb{R}^3)$ by Lemma \ref{alem1} in Appendix. This implies that for \eqref{y1} the corresponding initial data $(U_0,U_1)\in  H^3(\mathbb{R}^3)\times H^2(\mathbb{R}^3)$.

As mentioned before, two pairs of characteristic strips ($\mathcal{R}_2$ and $\mathcal{R}_3$, $\mathcal{R}_4$ and $\mathcal{R}_5$) may overlap for a long time. Our strategy is to consider the following four separated strips:
 \begin{equation*}
\{\mathcal{R}_1,\mathcal{R}_2\cup\mathcal{R}_3,\mathcal{R}_4\cup\mathcal{R}_5,\mathcal{R}_6\},
  \end{equation*}
  and to derive desired estimates in $\mathcal{R}_1$.
  By calculation, we find these four strips are separated when $t>t_0^{(\eta)}:=\frac\eta\sigma$. Here, $\sigma$ is a constant that describes the amplitude of the eigenvalues, a uniform nonzero constant. The main quantities to be estimated are as follows:
\begin{align}
  S_i(t):=&\sup_{(z'_i,s'_i)\atop z'_i\in[\eta,2\eta],\ 0\leq s'_i\leq t}\rho_i(z'_i,s'_i),&S(t):=&\max_{i\in\{1,2,3,4,5,6 \}}S_i(t),\label{ingr}\\
  J_i(t):=&\sup_{(z'_i,s'_i)\atop z'_i\in[\eta,2\eta]\ 0\leq s'_i\leq t}|v^i(z'_i,s'_i)|,&J(t):=&\max_{i\in\{1,2,3,4,5,6 \}}J_i(t),\\
  \bar{U}(t):=&\sup_{(x',t')\atop 0\leq t'\leq t}|\Phi(x',t')|,&W(t):=&\max_{i\in\{1,2,3,4,5,6 \}}\sup_{(x',t')\atop 0\leq t'\leq t}|w^i(x',t')|,\\
 V_1(t):= & \sup_{(x',t')\notin\mathcal{R}_1,\atop 0\leq t'\leq t}|w^1(x',t')|,&V_{\bar{2}}(t):=&\sup_{(x',t')\notin\mathcal{R}_2\cup\mathcal{R}_3,\atop 0\leq t'\leq t}\{|w^2(x',t')|,|w^3(x',t')|\},\\
V_6(t):=&\sup_{(x',t')\notin\mathcal{R}_6,\atop 0\leq t'\leq t}|w^6(x',t')|,
  &V_{\bar{5}}(t):=&\sup_{(x',t')\notin\mathcal{R}_4\cup\mathcal{R}_5,\atop 0\leq t'\leq t}\{|w^4(x',t')|,|w^5(x',t')|\},\\
V(t):=&\max_{i=1,\bar{2},\bar{5},6}V_i(t). &&\label{norm}
\end{align}
See Section \ref{datanorm} for the precise expressions and the details.

\textbf{Step 4: A priori estimates}

We give estimates of the quantities in \eqref{ingr}-\eqref{norm}. Assume, for some $T>0$, the Cauchy problem of system \eqref{0y2} has a solution $\Phi\in C^\infty(\mathbb{R}\times [0,T]; B_{2\delta}^6(0))$. We then derive estimates of these quantities in the non-separated characteristic region $t\in[0,t_0^{(\eta)}]$ and in the separated characteristic regions $t\in[t_0^{(\eta)},T]$. We will show that the formation of shock happens in the separated regions.

\begin{itemize}
\item{\bf Estimates for $\mathbf{t\in[0,t_0^{(\eta)}]}$}

In the non-separated region $t\in[0,t_0^{(\eta)}]$, all the characteristic strips $\{\mathcal{R}_i\}$ overlap. All the quantities defined in \eqref{ingr}-\eqref{norm} are bounded in such a way
\begin{align*}
  W(t)=&O(W_0^{(\eta)}),\ V(t)=O(\eta [W_0^{(\eta)}]^2),\ S(t)=O(1),\\
   J(t)=&O(W_0^{(\eta)}),\ \bar{U}(t)=O(\eta W_0^{(\eta)}).
\end{align*}
In particular, for the inverse foliation densities $\{\rho_i\}$, we have
$$ \rho_i(z_i,t)\geq 1-\varepsilon,\quad i=1,\cdots,6,\quad \forall\ t\in[0,t_0^{(\eta)}],$$
for some given small parameter $\varepsilon$.
See Section \ref{3d} for the details.

\item{\bf Estimates for $\mathbf{t\in[t_0^{(\eta)},T]}$}

For $t\in[t_0^{(\eta)},T]$, though $\mathcal{R}_2$ and $\mathcal{R}_3$ or $\mathcal{R}_4$ and $\mathcal{R}_5$ are not completely separated due to the non-strict hyperbolicity, the aforementioned four characteristic strips $\mathcal{R}_1,\mathcal{R}_2\cup\mathcal{R}_3,\mathcal{R}_4\cup\mathcal{R}_5,\mathcal{R}_6$ are well separated. We give estimates in different strips, and show that the first singularity (a shock) happens in $\mathcal{R}_1$ along $\mathcal{C}_1$. Subtle cancellations are explored. When $ t\in[t_0^{(\eta)},T]$, a priori estimates are obtained in these strip regions:
\begin{align*}
    S=&O(1+tVS+tJ+\eta SJ),\\
    J=&O\Big(W_0^{(\eta)}+t VJ+tV^2 S+\eta J^2+\frac{tV SJ}{\underline{S}}\Big),\\
    V=&O\Big(\eta [W_0^{(\eta)}]^2+tV^2+\eta VJ+\eta\frac{\varepsilon}{\underline{S}}J^2\Big), \\
    \bar{U}=&O(\eta J+\eta V+\eta tV).
\end{align*}
Here $\underline{S}(t)$ is the infimum of inverse foliation densities except for $\rho_1$:
\begin{equation*}
  \underline{S}(t):=\min_{i\in\{2,\cdots,6\}}\inf_{z'\in[\eta,2\eta]\atop 0\leq t'\leq t} \rho_i(z',t').
\end{equation*}

\end{itemize}

\textbf{Step 5: Bounds of the norms and a positive lower bound for $\underline{S}$}

Based on the above estimates, we obtain a desired positive lower bound of $\underline{S}$ via a bootstrap argument. We prove:\footnote{See Section \ref{bootargu} and Section \ref{lows} for the details.}
\begin{equation*}
J=O(W_0^{(\eta)}),\, S=O(1),\, tV=O(\eta W_0^{(\eta)}+\eta\varepsilon\theta^{-\frac13} W_0^{(\eta)}),
\end{equation*}
\begin{equation*}
V=O(\eta  [W_0^{(\eta)}]^2+\eta\varepsilon\theta^{-\frac13} [W_0^{(\eta)}]^2)
\end{equation*}
and
\begin{equation*}
    \underline{S}(s)\geq \frac{(1-\varepsilon)^2}{2}.
\end{equation*}

\textbf{Step 6: Bound for $\partial_{z_1}\rho_1$}

With characteristic coordinates and bi-characteristic coordinates, we prove that $\tau_{1}^{(6)}:=\partial_{y_1}\rho_1$ and $\pi_{1}^{(6)}:=\partial_{y_1}v^1$ satisfy a linear system:
  \begin{equation*}
  \left\{ \begin{split}
  \partial_{s_1}\tau_{1}^{(6)}=&B_{11}^\eta\tau_1^{(6)}+B_{12}^\eta\pi_1^{(6)}+B_{13}^\eta\\
  \partial_{s_1}\pi_{1}^{(6)}=&B_{21}^\eta\tau_1^{(6)}+B_{22}^\eta\pi_1^{(6)}+B_{23}^\eta,
  \end{split}\right.
\end{equation*}
where $\{B_{ij}^\eta\}_{i=1,2;j=1,2,3}$ are uniformly bounded constants depending on $\eta$. Hence $\tau_{1}^{(6)}:=\partial_{y_1}\rho_1$ is bounded. With expression
\begin{equation*}
 \partial_{z_1}\rho_1=\partial_{y_1}\rho_1+\frac{\rho_1}{2\lambda_1}\partial_{s_1}\rho_1=\partial_{y_1}\rho_1+\frac{\rho_1}{2\lambda_1}\big(c_{11}^1v^1+\sum_{k\neq 1}c_{1k}^1w^k\rho_1\big),
\end{equation*}
and estimate for $v^1$ and estimates for $w^k$ where $k\neq1$ in Step 4, we then have a bound for $\partial_{z_1}\rho_1$ by a uniform constant depending on $\eta$.

\textbf{Step 7: Ill-posedness mechanism}

The formation of shock happens in the first characteristic strip $\mathcal{R}_1$. With the obtained estimates, we control $\rho_1$ by a sharp description:
\begin{equation*}
\rho_1 \sim 1-t W_0^{(\eta)}.
\end{equation*}
This implies that the shock forms at a time $T_{\eta}^* $ with
\begin{equation*}
  T_{\eta}^*  \sim \frac{1}{W_0^{(\eta)}}.
\end{equation*}
With our constructed initial data, which are finite in $H^3(\mathbb{R}^3)$, we calculate the Sobolev norms at later moments. We observe that as $t \to  T_\eta^* $, the $H^2(\mathbb{R}^3)$ norm of solutions to \eqref{y1} at time $t$ approaches to infinite, and it is driven by shock formation. That is
\begin{equation*}
  \|U_\eta(\cdot,T_\eta^*)\|_{H^2(\Omega_{T_\eta^*})}\gtrsim\|w^1(\cdot,T_\eta^*)\|_{L^2(\Omega_{T_\eta^*})}=+\infty.
\end{equation*}
 Moreover, as $\eta \to 0$, we have $T_\eta^*\to 0$. This is the desired ill-posedness result for 3D elastic waves in $H^3(\mathbb{R}^3)$. See Section \ref{ill} for the details.

\subsection{Other related works}
In this section, we refer to some related works on the local well-posedness of low regularity solutions.

For quasilinear wave equations, Bahouri-Chemin \cite{bahouri-chemin1,bahouri-chemin2} and Tataru \cite{tataru1}  showed the local well-posedness in $H^s(\mathbb{R}^n)$ with $s>\frac{n}{2}+\frac{3}{4}$(later enhanced to $s>\frac{n}{2}+\frac{1}{2}+\frac{1}{6}$ in \cite{tataru3}) with Strichartz estimates. In \cite{klainerman-igor} Klainerman-Rodnianski improved their results to $H^{s_*}(\mathbb{R}^3)$ with $s_*>s_0=2+(2-\sqrt{3})/2$ by combing a geometric method with a paradifferential calculation. Tataru-Smith obtained in \cite{tataru} sharper results where they used a parametrix construction to improve the estimates. They showed that for $n$ dimensional quasilinear wave equations, the Cauchy problems are locally well-posed in $H^s(\mathbb{R}^n)$ with $s>n/2+3/4$ for $n=2$ and $s>(n+1)/2$ for $n=3,4,5$. The corresponding low regularity result for Einstein vacuum equations was obtained by Klainerman-Rodnianski \cite{klainerman-igor05} via a vectorfield method. In \cite{Wang}, Wang adopted the geometric vectorfield approach and show the $H^s(\mathbb{R}^3)$ local well-posedness of three dimensional quasilinear wave equations for any $s>2$. This result is generalized to the 3D compressible Euler equations in \cite{Disconzi} by Disconzi-Luo-Mazzone-Speck, via a decomposition of Euler flow into the wave part and the transport part. Recently, Wang \cite{Wang19} revisited this topic and enhanced the regularity of the transport part in \cite{Disconzi} by $1/2$-derivative.

 For Einstein vacuum equations, it was shown in \cite{klainerman-igor05} that the local existence holds for initial data $g_0$ being in $H^s(\mathbb{R}^3)$ for $s>2$. And it can be seen from the counterexamples Ettinger-Lindblad constructed in \cite{lindblad17} that the above result is sharp in wave coordinates. Through a series of remarkable works by Klainerman, Rodnianski, Szeftel \cite{S1,S2,S3,S4,klainerman-igor-szeftel}, the celebrated $L^2$ bounded curvature theorem with  Yang-Mills frames was achieved. They proved that Einstein vacuum equations are well-posedness in $H^2(\mathbb{R}^3)$.

For elastic wave equations, local well-posedness in $H^s(\mathbb{R}^3)$ with $s>3/2+2$ was obtained by Hughes-Kato-Marsden \cite{kato}.  For the radially symmetric case, Hidano-Zha \cite{zha} proved almost global existence for small weighted initial data in $H^3_{\text{rad}}(\mathbb{R}^3)$.

\section{Reduction of Equations}\label{nstricth}
For $Y=(Y^1,Y^2,Y^3)\in\mathbb{R}^3$, $U(Y,t)=(U^1(Y,t),U^2(Y,t),U^3(Y,t))$ being a solution to \eqref{y1}. Under planar symmetry(with respect to $Y^2$ and $Y^3$), we have
\begin{equation*}
\begin{split}
U^1(Y,t)=u^1(Y^1,t),\ U^2(Y,t)=u^2(Y^1,t),\ U^3(Y,t)=u^3(Y^1,t).
\end{split}
\end{equation*}
Denote $Y^1$ to be $x$. Then $u(x,t)=(u^1(x,t),u^2(x,t),u^3(x,t))$ is a solution to \eqref{y1}. In the following proposition, we will show that $u(x,t)$ satisfies a quasilinear wave system with multiple speeds:
\begin{prop}\label{prop}
  Assume that $u(x,t)$ is a solution of \eqref{y1} under a planar symmetry. Then $u=(u^1,u^2,u^3)$ verifies the following quasilinear wave system:
  \begin{align}\label{y2.2}
    \left\{\begin{array}{lll}
    \partial_t^2u^1-c_1^2\partial_x^2u^1=\sigma_0\partial_x(\partial_xu^1)^2+\sigma_1\partial_x(\partial_xu^2)^2+\sigma_1\partial_x(\partial_xu^3)^2,\\
    \partial_t^2u^2-c_2^2\partial_x^2u^2=2\sigma_1\partial_x(\partial_xu^1\partial_xu^2),\\
    \partial_t^2u^3-c_2^2\partial_x^2u^3=2\sigma_1\partial_x(\partial_xu^1\partial_xu^3),
    \end{array}\right.
\end{align}
with initial data $u(x,0)=U_0(x)$, $u_t(x,0)=U_1(x)$.
\end{prop}
{\it Proof.}
For $U=u(x,t)=(u^1(x,t),u^2(x,t),u^3(x,t))$, we revisit \eqref{y1}. With definition of $Q_{jk}(f,g)=\partial_jf\partial_kg-\partial_kf\partial_jg$ ($j\neq k$) for $f,g\in C^1(\mathbb{R})$, it is a straight forward check that $Q^1,Q^2,Q^3$ given in \eqref{17} are all zero. Hence
 \begin{align*}
Q(U,\nabla U)=0.
\end{align*}
We then calculate the remaining terms, since
\begin{equation*}
  \text{curl}\ U=(0,-\partial_xu^3,\partial_xu^2),
\end{equation*}
it holds that
\begin{equation*}
  \nabla|\text{curl}\  U|^2=\big(\partial_x(\partial_xu^2)^2+\partial_x(\partial_xu^3)^2,0,0\big)
  \end{equation*}
  and
  \begin{equation*}
                            \text{curl}\ (\text{div}\  U\text{curl}\  U)=\big( 0, -\partial_x(\partial_xu^1\partial_xu^2), -\partial_x(\partial_xu^1\partial_xu^3)\big).
\end{equation*}
Back to \eqref{y1}-\eqref{17}, Proposition \ref{prop} is proved. \hfill$\Box$

We then change \eqref{y2.2} into a first-order hyperbolic system.
Let us recall some important definitions. Assume unknowns $\Phi\in\mathbb{R}^n$ and a coefficient matrix $A(\Phi)\in\mathbb{R}^n\times\mathbb{R}^n$ satisfy a first-order system:
\begin{equation}\label{y2.1}
  \partial_t\Phi+A(\Phi)\partial_x\Phi=0.
\end{equation}
The system \eqref{y2.1} is called {\it hyperbolic} if $A(\Phi)$ has $n$ real eigenvalues noted as: $$\lambda_1(\Phi),\cdots,\lambda_n(\Phi)$$ and $A(\Phi)$ is diagonalizable.
If furthermore all of the eigenvalues are distinct, then \eqref{y2.1} is called {\it strictly hyperbolic}. Otherwise, system \eqref{y2.1} is called {\it non-strictly hyperbolic}.

The following Lemma shows that the elastic plane waves in 3D verify a non-strictly hyperbolic system. We will use a notation $B_{2\delta}^6(0)$ to mean an open ball of radius $2\delta$ around $0\in\mathbb{R}^6$.
Let
\begin{align*}
  \phi_1:=\partial_xu^1,\quad \phi_2:=\partial_xu^2,\quad \phi_3:=\partial_xu^3,\quad \phi_4:=\partial_tu^1,\quad \phi_5:=\partial_tu^2,\quad \phi_6:=\partial_tu^3.
\end{align*}
\begin{lem}
  Let $\Phi=(\phi_1,\phi_2,\phi_3,\phi_4,\phi_5,\phi_6)^T$. Then system \eqref{y2.2} is equivalent to
\begin{align}\label{y2}
  \partial_t\Phi+A(\Phi)\partial_x\Phi=0,
\end{align}
where
\begin{align*}
  \begin{split}
    A(\Phi)=\left(
              \begin{array}{cccccc}
                0 & 0 & 0 &-1 & 0  & 0 \\
                0 & 0 & 0 &0  & -1 & 0 \\
                0 & 0 & 0 &0  & 0  & -1 \\
                -(c_1^2+2\sigma_0\phi_1) & -2\sigma_1\phi_2 & -2\sigma_1\phi_3 &0  & 0  & 0 \\
                -2\sigma_1\phi_2 & -(c_2^2+2\sigma_1\phi_1) & 0 &0  & 0  & 0 \\
                -2\sigma_1\phi_3 & 0 & -(c_2^2+2\sigma_1\phi_1) &0  & 0  & 0 \\
              \end{array}
            \right).
  \end{split}
\end{align*}
Moreover, system \eqref{y2} is \textit{not} uniformly strictly hyperbolic for $\Phi\in B_{2\delta}^6(0)$.
\end{lem}
{\it Proof.}
For notational simplicity, we set
\begin{equation}\label{y2.6}
  a=c_1^2+2\sigma_0\phi_1,\quad b=c_2^2+2\sigma_1\phi_1,\quad c=2\sigma_1\phi_2,\quad d=2\sigma_1\phi_3.
\end{equation}
We then calculate the eigenvalues of $A(\Phi)$. By direct computing, we have
\begin{align*}
  \begin{split}
      \det{(\lambda I-A)}&=\lambda^6-(a+2b)\lambda^4+(2ab+b^2-c^2-d^2)\lambda^2-ab^2+b(c^2+d^2)\\
      &=(\lambda^2-b)[(\lambda^2-a)(\lambda^2-b)-(c^2+d^2)].
  \end{split}
\end{align*}
And its six roots are:
\begin{align}
    \lambda_{1}=&\sqrt{\frac12(a+b)+\frac12\sqrt{(a-b)^2+4(c^2+d^2)}},\label{lambda1}\\
    \lambda_2=&\sqrt{b},\label{211}\\
\lambda_3=&\sqrt{\frac12(a+b)-\frac12\sqrt{(a-b)^2+4(c^2+d^2)}},\label{212}\\
\lambda_4=&-\sqrt{\frac12(a+b)-\frac12\sqrt{(a-b)^2+4(c^2+d^2)}},\\
\lambda_5=&-\sqrt{b},\\
\lambda_6=&-\sqrt{\frac12(a+b)+\frac12\sqrt{(a-b)^2+4(c^2+d^2)}}.\label{lambda6}
\end{align}
Note that for $|\Phi|<2\delta$ being small, we have
\begin{align}\label{y3}
  \begin{split}
    \lambda_6(\Phi)<\lambda_5(\Phi)\leq\lambda_4(\Phi)<\lambda_3(\Phi)\leq\lambda_2(\Phi)<\lambda_1(\Phi).
  \end{split}
\end{align}
  Equalities in \eqref{y3} hold when $\phi_2=\phi_3=0$, i.e., $c=d=0$. In this case
$a-b=c_1^2+2(\sigma_0-\sigma_1)\phi_1-c_2^2>0$ and
\begin{equation*}
  \lambda_3=\sqrt{\frac12(a+b)-\frac12\sqrt{(a-b)^2}}=\sqrt{b}=\lambda_2,
\end{equation*}
\begin{equation*}
  \lambda_4=-\sqrt{\frac12(a+b)-\frac12\sqrt{(a-b)^2}}=-\sqrt{b}=\lambda_5.
\end{equation*}
For $\Phi=0\in\mathbb{R}^6$, we further have
\begin{align}\label{219}
  \begin{split}
  \lambda_6=-c_1,\quad  \lambda_5(0)=\lambda_4(0)=-c_2,\quad\lambda_3(0)=\lambda_2(0)=c_2,\quad\lambda_1=c_1.
  \end{split}
\end{align}
Though $A(0)$ has repeated eigenvalues, we still can find six linearly independent right eigenvectors. By calculation, we have
\begin{align*}
  \begin{array}{llll}
    &e_{01}=\left(
             \begin{array}{c}
               1 \\
               0 \\
               0 \\
               -c_1 \\
               0 \\
               0 \\
             \end{array}
           \right),\ e_{02}=\left(
                              \begin{array}{c}
                                0 \\
                                1 \\
                                0 \\
                                0 \\
                                -c_2 \\
                                0 \\
                              \end{array}
                            \right),\ e_{03}=\left(
                                               \begin{array}{c}
                                                 0 \\
                                                 0 \\
                                                 1 \\
                                                 0 \\
                                                 0 \\
                                                 -c_2 \\
                                               \end{array}
                                             \right),\\
                            & e_{04}=\left(
                                     \begin{array}{c}
                                                                  0 \\
                                                                  0 \\
                                                                  1 \\
                                                                  0 \\
                                                                  0 \\
                                                                  c_2 \\
                                                                \end{array}
                                                              \right),\ e_{05}=\left(
                                                                                 \begin{array}{c}
                                                                                   0 \\
                                                                                   1 \\
                                                                                   0 \\
                                                                                   0 \\
                                                                                   c_2 \\
                                                                                   0 \\
                                                                                 \end{array}
                                                                               \right),\ e_{06}=\left(
                                                                                                  \begin{array}{c}
                                                                                                    1 \\
                                                                                                    0 \\
                                                                                                    0 \\
                                                                                                    c_1 \\
                                                                                                    0 \\
                                                                                                    0 \\
                                                                                                  \end{array}
                                                                                                \right).
  \end{array}
\end{align*}
Similarly, for eigenvalues $\lambda_1$ to $\lambda_6$ as in (\ref{lambda1})-(\ref{lambda6}) we compute their corresponding right eigenvectors with $A(\Phi)$:
\begin{equation}\label{right}
\begin{split}
  r_1=\left(
                \begin{array}{c}
                  \frac{\lambda_1^2-b}{2\sigma_1} \\
                  \phi_2 \\
                  \phi_3 \\
                  -\frac{\lambda_1(\lambda_1^2-b)}{2\sigma_1} \\
                  -\lambda_1\phi_2 \\
                  -\lambda_1\phi_3 \\
                \end{array}
              \right),\
                r_2=\left(
                \begin{array}{c}
                  0 \\
                 \phi_3 \\
                  -\phi_2 \\
                  0 \\
                  -\lambda_2\phi_3 \\
                  \lambda_2\phi_2 \\
                \end{array}
              \right),\
               r_3=\left(
                \begin{array}{c}
                  \frac{\lambda_3^2-b}{2\sigma_1}\\
                  \phi_2 \\
                  \phi_3 \\
                  -\frac{\lambda_3(\lambda_3^2-b)}{2\sigma_1} \\
                  -\lambda_3\phi_2 \\
                  -\lambda_3\phi_3 \\
                \end{array}
              \right),\\
 r_4=\left(
                                     \begin{array}{c}
                                       \frac{\lambda_3^2-b}{2\sigma_1} \\
                                       \phi_2 \\
                                       \phi_3\\
                                        \frac{\lambda_3(\lambda_3^2-b)}{2\sigma_1} \\
                                       \lambda_3\phi_2 \\
                                       \lambda_3\phi_3 \\
                                     \end{array}
                                   \right),\
                            r_5=\left(
                                      \begin{array}{c}
                                        0 \\
                                        \phi_3 \\
                                        -\phi_2 \\
                                        0 \\
                                        \lambda_2\phi_3 \\
                                        -\lambda_2\phi_2 \\
                                      \end{array}
                                    \right),\
                                    \ r_6=\left(
                                     \begin{array}{c}
                                       \frac{\lambda_1^2-b}{2\sigma_1} \\
                                       \phi_2 \\
                                       \phi_3\\
                                       \frac{\lambda_1(\lambda_1^2-b)}{2\sigma_1} \\
                                       \lambda_1\phi_2\\
                                       \lambda_1\phi_3 \\
                                     \end{array}
                                   \right).
 \end{split}
\end{equation}
 \hfill$\Box$

In order to use John's approach in \cite{john74}, we also compute the left eigenvectors $\{l^i(\Phi)\}_{i=1,\cdots,6}$ and require
 \begin{equation}\label{lfproduct}
  l^i(\Phi)\cdot r_j(\Phi)=\delta^i_j,\qquad\text{for}\quad i,j\in\{1,\cdots,6\}.
 \end{equation}
Requirement \eqref{lfproduct} is important for deriving the decomposition of waves in the next Section. Via calculations, we have
\begin{equation}\label{left}
\begin{split}
l^1=&\frac1{K}\Big(\frac{\lambda_1^2-b}{2\sigma_1},\phi_2,\phi_3,-\frac{\lambda_1^2-b}{2\sigma_1\lambda_{1}},-\frac{\phi_2}{\lambda_{1}},-\frac{\phi_3}{\lambda_{1}}\Big),\\ l^2=&\frac1{M}\Big(0,\phi_3,-\phi_2,0,-\frac{\phi_3}{\lambda_{2}},\frac{\phi_2}{\lambda_{2}}\Big),\\
l^3=&\frac1{N}\Big(\frac{\lambda_3^2-b}{2\sigma_1},\phi_2,\phi_3,-\frac{\lambda_3^2-b}{2\sigma_1\lambda_{3}},-\frac{\phi_2}{\lambda_{3}},-\frac{\phi_3}{\lambda_{3}}\Big),\\  l^4=&\frac1{N}\Big(\frac{\lambda_3^2-b}{2\sigma_1},\phi_2,\phi_3,\frac{\lambda_3^2-b}{2\sigma_1\lambda_{3}},\frac{\phi_2}{\lambda_{3}},\frac{\phi_3}{\lambda_{3}}\Big),\\
l^5=&\frac1{M}\Big(0,\phi_3,-\phi_2,0,\frac{\phi_3}{\lambda_{2}},-\frac{\phi_2}{\lambda_{2}}\Big),\\
l^6=&\frac1{K}\Big(\frac{\lambda_1^2-b}{2\sigma_1},\phi_2,\phi_3,\frac{\lambda_1^2-b}{2\sigma_1\lambda_{1}},\frac{\phi_2}{\lambda_{1}},\frac{\phi_3}{\lambda_{1}}\Big),
\end{split}
\end{equation}
where
 \begin{align*}
   K=&\frac{(a-b)^2+4(c^2+d^2)}{4\sigma_1^2}+\frac{(a-b)\sqrt{(a-b)^2+4(c^2+d^2)}}{4\sigma_1^2},\\
   M=&2(\phi_2^2+\phi_3^2),\\
   N=&\frac{(a-b)^2+4(c^2+d^2)}{4\sigma_1^2}-\frac{(a-b)\sqrt{(a-b)^2+4(c^2+d^2)}}{4\sigma_1^2}.
 \end{align*}

\section{Decomposition of Waves}\label{yf}
In this section, we introduce useful characteristic coordinates and bi-characteristic coordinates. We assume that $\Phi\in C^2(\mathbb{R}\times[0,T],B_{2\delta}^6(0))$ is a solution to \eqref{y2} for some $T>0$. We define the $i^{\text{th}}$ characteristic $\mathcal{C}_i(z)$ to be the image $(X_i(z,t),t)$ of solutions of the following ODE:
\begin{align}\label{flow}
  \left\{\begin{array}{ll}
  \frac{\partial}{\partial t}X_i(z,t)=\lambda_i\big(\Phi(X_i(z,t),t)\big),\quad t\in[0,T],\\
  X_i(z,0)=z.
  \end{array}\right.
\end{align}
Define $\mathcal{R}_i$ to be the $i^{\text{th}}$ characteristic strip evolving from initial data supported in $I_0$. That is,
\begin{equation*}
  \mathcal{R}_i=\bigcup_{z\in I_0}\mathcal{C}_i(z).
\end{equation*}
For any given $(x,t)\in\mathbb{R}\times[0,T]$, there is a unique $(z_i,s_i)\in\mathbb{R}\times[0,T]$ such that along the $i^{\text{th}}$ characteristic $\mathcal{C}_i(z)$ with $X_i(z_i,0)=z_i$ we have
\begin{equation}\label{cc}
  (x,t)=\big(X_i(z_i,s_i),s_i\big).
\end{equation}
We then define inverse foliation density
\begin{equation}\label{dense}
  \rho_i:=\frac{\partial}{\partial z_i}X_i(z_i,s_i),
\end{equation}
and it implies
\begin{equation}\label{35}
  \partial_{z_i}=\rho_i\partial_x,\quad \partial_{s_i}=\lambda_i\partial_x+\partial_t
\end{equation}
and
\begin{equation}
  dx=\rho_idz_i+\lambda_ids_i,\quad dt=ds_i.
\end{equation}

The bi-characteristic coordinates are introduced to study the intersection of the $i^{\text{th}}$ and  $j^{\text{th}}$ characteristics when $i\neq j$. Set $t=t'(y_i,y_j)$ to be a function of $y_i,y_j$ and require
\begin{equation}\label{37}
  (x,t)=\big(X_i(y_i,t'(y_i,y_j)),t'(y_i,y_j)\big)=\big(X_j(y_j,t'(y_i,y_j)),t'(y_i,y_j)\big)
\end{equation}
with
\begin{equation}\label{38}
  \partial_{y_i}t'=\frac{\rho_i}{\lambda_j-\lambda_i},\quad\partial_{y_j}t'=\frac{\rho_j}{\lambda_i-\lambda_j}.
\end{equation}
Note that \eqref{38} implies $[\partial_{y_i},\partial_{y_j}]=0$ for $1\leq i\neq j\leq6$. We then have
\begin{equation}\label{y17}
  dx=\frac{\rho_i\lambda_j}{\lambda_j-\lambda_i}dy_i+\frac{\rho_j\lambda_i}{\lambda_i-\lambda_j}dy_j,\quad dt=\frac{\rho_i}{\lambda_j-\lambda_i}dy_i+\frac{\rho_j}{\lambda_i-\lambda_j}dy_j
\end{equation}
and
\begin{equation}
  dz_i=dy_i,\quad dz_j=dy_j.
\end{equation}

Now we are ready to study the decomposition of waves. For fixed $i\in\{1,\cdots,6\}$, let
\begin{equation}\label{y4}
  w^i:=l^i\partial_x\Phi,\quad\text{and}\quad
  v^i:=l^i\partial_{z_i}\Phi=\rho_iw^i.
\end{equation}
The transformation $w^1=l^1\partial_x\Phi$ is well-defined for $|\Phi|<2\delta$, since
  \begin{equation*}
    l^1(0)=\Big(\frac{\sigma_1}{c_1^2-c_2^2},0,0,-\frac{\sigma_1}{c_1(c_1^2-c_2^2)},0,0\Big)
  \end{equation*}
being a non-zero vector.
Using \eqref{lfproduct} and \eqref{y4}, we have
\begin{equation}\label{uxx}
  \partial_x\Phi=\sum_k w^kr_k.
\end{equation}
By John's formula \cite{john74}, we diagonalize the system \eqref{y2} as:
\begin{align}\label{y12}
  \partial_{s_i}w^i=&-c_{ii}^i(w^i)^2+\Big(\sum_{m\neq i}(-c_{im}^i+\gamma_{im}^i)w^m\Big)w^i+\sum_{m\neq i,k\neq i\atop m\neq k}\gamma_{km}^iw^kw^m,
\end{align}
where
\begin{equation}\label{y8}
  c_{im}^i=\nabla_\Phi\lambda_i\cdot r_m,
\end{equation}
and
\begin{align}\label{gamma}
\begin{array}{lll}
  \gamma_{im}^i=&-(\lambda_i-\lambda_m)l^i \cdot(\nabla_\Phi r_i \cdot r_m-\nabla_\Phi r_m \cdot r_i),\quad &\text{for}\quad m\neq i,\\
  \gamma_{km}^i=&-(\lambda_k-\lambda_m)l^i \cdot (\nabla_\Phi r_k \cdot r_m),\quad& \text{for}\quad k\neq i,\ \text{and}\ m\neq i.
  \end{array}
\end{align}
Moreover, computed as in Christoudoulou-Perez \cite{christodoulou}, the inverse foliation density of characteristics $\rho_i$ and quantity $v^i$ satisfy:
\begin{align}
  \partial_{s_i}\rho_i=&c_{ii}^iv^i+\Big(\sum_{m\neq i}c_{im}^iw^m\Big)\rho_i,\label{y14}\\
  \partial_{s_i}v^i=&\Big(\sum_{m\neq i}\gamma_{im}^iw^m\Big)v^i+\sum_{m\neq i,k\neq i\atop m\neq k}\gamma_{km}^iw^kw^m\rho_i.\label{y13}
\end{align}

We next analyze the detailed algebraic structures of \eqref{y12}, \eqref{y14} and \eqref{y13} when applying to elastic waves. We have
\begin{lem}\label{GN}
The $1^{\text{st}}$ and $6^{\text{th}}$ characteristics of system \eqref{y2} are genuinely non-linear in the sense of Lax:
 \begin{equation*}
   \nabla_\Phi\lambda_i(\Phi) \cdot r_i(\Phi)\neq 0,\quad i=1,6,\quad\forall\ \Phi\in B_{2\delta}^6(0).
 \end{equation*}
 The $2^{\text{nd}}$ and $5^{\text{th}}$ characteristics of system \eqref{y2} are linearly degenerate in the sense of Lax:
\begin{equation*}
  \nabla_\Phi\lambda_j \cdot r_j=0,\quad j=2,5,\quad\forall\ \Phi\in B_{2\delta}^6(0).
\end{equation*}
\end{lem}
{\it Proof.}
For notational simplicity, we define
\begin{equation*}
  \Delta:=(a-b)^2+4(c^2+d^2).
\end{equation*}
By direct calculation, we get
\begin{align}\label{gradlam}
\begin{array}{llll}
  &\nabla_\Phi\lambda_1=-\nabla_\Phi\lambda_6=\Big(\frac{(\sigma_0+\sigma_1)\sqrt{\Delta}+(a-b)(\sigma_0-\sigma_1)}{2\lambda_1\sqrt{\Delta}},\frac{4\sigma_1^2\phi_2}{\lambda_1\sqrt{\Delta}},\frac{4\sigma_1^2\phi_3}{\lambda_1\sqrt{\Delta}},0,0,0\Big),\\
  &\nabla_\Phi\lambda_2=-\nabla_\Phi\lambda_5=\Big(\frac{\sigma_1}{\lambda_2},0,0,0,0,0\Big),\\
  &\nabla_\Phi\lambda_3=-\nabla_\Phi\lambda_4=\Big(\frac{(\sigma_0+\sigma_1)\sqrt{\Delta}-(a-b)(\sigma_0-\sigma_1)}{2\lambda_3\sqrt{\Delta}},-\frac{4\sigma_1^2\phi_2}{\lambda_3\sqrt{\Delta}},-\frac{4\sigma_1^2\phi_3}{\lambda_3\sqrt{\Delta}},0,0,0\Big).
\end{array}\end{align}
With \eqref{right} and \eqref{gradlam}, in the expression of \eqref{y8} we have
\begin{equation*}
  c_{22}^2(\Phi)=c_{55}^5(\Phi)=0.
\end{equation*}
And it also holds
\begin{equation*}
  c_{11}^1(\Phi)=-c_{66}^6(\Phi)=\frac{2\sigma_0(a-b)(\lambda_1^2-b)+(2\sigma_0+6\sigma_1)(c^2+d^2)}{4\sigma_1\lambda_1\sqrt{\Delta}}.
\end{equation*}
By a direct check of
\begin{equation*}
c_{11}^1(0)=-c_{66}^6(0)=\frac{\sigma_0(c_1^2-c_2^2)}{2\sigma_1c_1}\neq 0,
\end{equation*}
we have
\begin{equation*}
 c_{11}^1(\Phi)=-c_{66}^6(\Phi)\neq 0
\end{equation*}
for any $\Phi\in B_{2\delta}^6(0)$ with sufficiently small $\delta$.
\hfill$\Box$

Without loss of generality, we assume that
\begin{equation*}
 c_{11}^1(0)<0\quad\text{and}\quad c_{11}^1(\Phi)<0,\quad \forall\ \Phi\in B_{2\delta}^6(0).
\end{equation*}
Correspondingly, we require the initial data of $w^1(z,0)$ to be positive at some $z_0\in I_0$. This condition ensures the positivity of $v^1$ in a bootstrap argument.

In \cite{christodoulou}, Christodoulou-Perez studied mechanism for shock formation of strictly hyperbolic systems. However, here elastic waves are reduced to a non-strictly hyperbolic system below. We explore some subtle structures in \eqref{y12}. These structures are critical for the proof of the main theorem.
\begin{prop}
For the formula \eqref{y12} of the non-strictly hyperbolic system \eqref{y2}, with $\Phi\in B_{2\delta}^6(0)$ we have the following properties:
\begin{itemize}
  \item For $i=1,3,4,6$, terms of $w^2w^3$ and $w^4w^5$ disappear;
  \item For $i=2,5$, the coefficients of $w^2w^3$ and $w^4w^5$ always have a small factor $\lambda_2(\Phi)-\lambda_3(\Phi)$ in front.
\end{itemize}
\end{prop}
{\it Proof.} We may worry about the interaction of waves with almost-the-same speeds. The corresponding terms are $w^2w^3$ and $w^4w^5$.
To get the coefficient of $w^2w^3$, with \eqref{right}, \eqref{left}, \eqref{y12} and \eqref{gamma}, we first calculate
\begin{equation*}
  \nabla_\Phi r_2=\left(
                    \begin{array}{cccccc}
                      0 & 0 & 0 & 0 & 0 & 0 \\
                      0 & 0 & 1 & 0 & 0 & 0 \\
                      0 & -1 & 0 & 0 & 0 & 0 \\
                      0 & 0 & 0 & 0 & 0 & 0 \\
                      -\phi_3\partial_{\phi_1}\lambda_2 & 0 & -\lambda_2 & 0 & 0 & 0 \\
                      \phi_2\partial_{\phi_1}\lambda_2 & \lambda_2 & 0 & 0 & 0 & 0 \\
                    \end{array}
                  \right)
\end{equation*}
and
\begin{equation*}
\begin{split}
  &\nabla_\Phi r_3=\\
  &\left(
                    \begin{array}{cccccc}
                      \frac{\partial_{\phi_1}(\lambda_3^2-b)}{2\sigma_1} & \frac{\partial_{\phi_2}(\lambda_3^2-b)}{2\sigma_1} & \frac{\partial_{\phi_3}(\lambda_3^2-b)}{2\sigma_1} & 0 & 0 & 0 \\
                      0 & 1 & 0 & 0 & 0 & 0 \\
                      0 & 0 & 1 & 0 & 0 & 0 \\
                      -\frac{\lambda_3\partial_{\phi_1}(\lambda_3^2-b)+(\lambda_3^2-b)\partial_{\phi_1}\lambda_3}{2\sigma_1} & -\frac{\lambda_3\partial_{\phi_2}(\lambda_3^2-b)+(\lambda_3^2-b)\partial_{\phi_2}\lambda_3}{2\sigma_1} & -\frac{\lambda_3\partial_{\phi_3}(\lambda_3^2-b)+(\lambda_3^2-b)\partial_{\phi_3}\lambda_3}{2\sigma_1} & 0 & 0 & 0 \\
                      -\phi_2\partial_{\phi_1}\lambda_3 & -\lambda_3-\phi_2\partial_{\phi_2}\lambda_3 & -\phi_2\partial_{\phi_3}\lambda_3 & 0 & 0 & 0 \\
                      -\phi_3\partial_{\phi_1}\lambda_3 &-\phi_3\partial_{\phi_2}\lambda_3 & -\lambda_3-\phi_3\partial_{\phi_3}\lambda_3& 0 & 0 & 0 \\
                    \end{array}
                  \right).
  \end{split}
\end{equation*}
These imply
\begin{equation}\label{r23}
   \nabla_\Phi r_2 \cdot r_3=\left(
                \begin{array}{c}
                  0 \\
                 \phi_3 \\
                  -\phi_2 \\
                  0 \\
                  -\lambda_2\phi_3-\frac{\phi_3(\lambda_3^2-b)}{2\lambda_2} \\
                  \lambda_2\phi_2+\frac{\phi_2(\lambda_3^2-b)}{2\lambda_2} \\
                \end{array}
              \right),\quad   \nabla_\Phi r_3 \cdot r_2=\left(
                \begin{array}{c}
                  0 \\
                 \phi_3 \\
                  -\phi_2 \\
                  0 \\
                  -\lambda_3\phi_3 \\
                  \lambda_3\phi_2 \\
                \end{array}
              \right).
\end{equation}
For $w^4w^5$, by similar calculations we get
\begin{equation*}
\begin{split}
  &\nabla_\Phi r_4=\\
  &\left(
                    \begin{array}{cccccc}
                      \frac{\partial_{\phi_1}(\lambda_3^2-b)}{2\sigma_1} & \frac{\partial_{\phi_2}(\lambda_3^2-b)}{2\sigma_1} & \frac{\partial_{\phi_3}(\lambda_3^2-b)}{2\sigma_1} & 0 & 0 & 0 \\
                      0 & 1 & 0 & 0 & 0 & 0 \\
                      0 & 0 & 1 & 0 & 0 & 0 \\
                      \frac{\lambda_3\partial_{\phi_1}(\lambda_3^2-b)+(\lambda_3^2-b)\partial_{\phi_1}\lambda_3}{2\sigma_1} & \frac{\lambda_3\partial_{\phi_2}(\lambda_3^2-b)+(\lambda_3^2-b)\partial_{\phi_2}\lambda_3}{2\sigma_1} & \frac{\lambda_3\partial_{\phi_3}(\lambda_3^2-b)+(\lambda_3^2-b)\partial_{\phi_3}\lambda_3}{2\sigma_1} & 0 & 0 & 0 \\
                      \phi_2\partial_{\phi_1}\lambda_3 & \lambda_3+\phi_2\partial_{\phi_2}\lambda_3 & \phi_2\partial_{\phi_3}\lambda_3 & 0 & 0 & 0 \\
                      \phi_3\partial_{\phi_1}\lambda_3 &\phi_3\partial_{\phi_2}\lambda_3 & \lambda_3+\phi_3\partial_{\phi_3}\lambda_3& 0 & 0 & 0 \\
                    \end{array}
                  \right),
\end{split}
\end{equation*}
\begin{equation*}
  \nabla_\Phi r_5=\left(
                    \begin{array}{cccccc}
                      0 & 0 & 0 & 0 & 0 & 0 \\
                      0 & 0 & 1 & 0 & 0 & 0 \\
                      0 & -1 & 0 & 0 & 0 & 0 \\
                      0 & 0 & 0 & 0 & 0 & 0 \\
                      \phi_3\partial_{\phi_1}\lambda_2 & 0 & \lambda_2 & 0 & 0 & 0 \\
                      -\phi_2\partial_{\phi_1}\lambda_2 & -\lambda_2 & 0 & 0 & 0 & 0 \\
                    \end{array}
                  \right),
\end{equation*}
and then
\begin{equation}\label{r45}
  \nabla_\Phi r_4 \cdot r_5=\left(
                \begin{array}{c}
                  0 \\
                 \phi_3 \\
                  -\phi_2 \\
                  0 \\
                  \lambda_3\phi_3 \\
                  -\lambda_3\phi_2 \\
                \end{array}
              \right),\quad
               \nabla_\Phi r_5 \cdot r_4=\left(
                \begin{array}{c}
                  0 \\
                 \phi_3 \\
                  -\phi_2 \\
                  0 \\
                  \lambda_2\phi_3+\frac{\phi_3(\lambda_3^2-b)}{2\lambda_2} \\
                  -\lambda_2\phi_2-\frac{\phi_2(\lambda_3^2-b)}{2\lambda_2} \\
                \end{array}
              \right).
\end{equation}

Now we derive the coefficients of  $w^2w^3$ and $w^4w^5$ in the formulas of $\partial_{s_i}w^i (i=1,\cdots,6)$.

Take $\partial_{s_1}w^1$ and  $\partial_{s_2}w^2$  for examples. For $\partial_{s_1}w^1$, the coefficient of $w^2w^3$ is $\gamma_{23}^1+\gamma_{32}^1$, with $\gamma^{i}_{km}$ given in \eqref{gamma}. Since
\begin{equation*}
  l^1 \cdot (\nabla_\Phi r_2 \cdot r_3)=\frac1K\big[\phi_2\phi_3-\phi_2\phi_3-(\lambda_2+\frac{\lambda_3^2-b}{2\lambda_2})(\phi_2\phi_3-\phi_2\phi_3)\big]=0,
\end{equation*}
we have
\begin{equation}\label{gam1}
  \gamma_{23}^1=-(\lambda_2-\lambda_3)l^1 \cdot (\nabla_\Phi r_2 \cdot r_3)=0.
\end{equation}

We proceed to calculate $\gamma^1_{32}$. With the expression for $\nabla_{\Phi}\lambda_3$ in \eqref{gradlam}, it holds that
\begin{equation*}
 (\phi_3\partial_{\phi_2}-\phi_2\partial_{\phi_3})\lambda_3^2=-\frac{8\sigma_1^2}{\sqrt{\Delta}}\phi_2\phi_3+\frac{8\sigma_1^2}{\sqrt{\Delta}}\phi_2\phi_3=0.
\end{equation*}
Hence we have
\begin{align}\label{132}
\begin{array}{ll}
  &\gamma_{32}^1=-(\lambda_3-\lambda_2)l^1 \cdot(\nabla_\Phi r_3 \cdot r_2)\\
  =&\frac{\lambda_2-\lambda_3}{K}
  \left(
    \begin{array}{c}
      \diamondsuit\\
      \phi_2+\frac{\lambda_3}{\lambda_1}\phi_2+\big[\frac{\lambda_1^2-b}{4\sigma_1^2}+\frac{\lambda_3(\lambda_1^2-b)}{4\sigma_1^2\lambda_1}\big]\partial_{\phi_2}(\lambda_3^2-b)+\big[\frac{(\lambda_1^2-b)(\lambda_3^2-b)}{4\sigma_1^2\lambda_1}+\frac{\phi_2^2+\phi_3^2}{\lambda_1}\big]\partial_{\phi_2}\lambda_3 \\
      \phi_3+\frac{\lambda_3}{\lambda_1}\phi_3+\big[\frac{\lambda_1^2-b}{4\sigma_1^2}+\frac{\lambda_3(\lambda_1^2-b)}{4\sigma_1^2\lambda_1}\big]\partial_{\phi_3}(\lambda_3^2-b)+\big[\frac{(\lambda_1^2-b)(\lambda_3^2-b)}{4\sigma_1^2\lambda_1}+\frac{\phi_2^2+\phi_3^2}{\lambda_1}\big]\partial_{\phi_3}\lambda_3 \\
      0 \\
      0 \\
      0 \\
    \end{array}
  \right)^\top\\
  &\cdot\left(
                \begin{array}{c}
                  0 \\
                 \phi_3 \\
                  -\phi_2 \\
                  0 \\
                  -\lambda_2\phi_3 \\
                  \lambda_2\phi_2 \\
                \end{array}
              \right)\\
  =&\frac{\lambda_2-\lambda_3}{K}\Big[\frac{\lambda_1^2-b}{4\sigma_1^2}+\frac{\lambda_3(\lambda_1^2-b)}{4\sigma_1^2\lambda_1}+\frac{(\lambda_1^2-b)(\lambda_3^2-b)}{8\sigma_1^2\lambda_1\lambda_3}+\frac{\phi_2^2+\phi_3^2}{2\lambda_1\lambda_3}\Big](\phi_3\partial_{\phi_2}-\phi_2\partial_{\phi_3})\lambda_3^2\\
  =&0.
\end{array}
\end{align}
Notation ``$\diamondsuit$" in \eqref{132} represents certain functions depending on the unknown $\Phi$. The detailed expression is not used in the above calculations, so we omit it.
From \eqref{gam1} and \eqref{132}, we conclude the coefficient of $w^2w^3$ in $\partial_{s_1}w^1$ is
\begin{equation*}
 \gamma_{23}^1(\Phi)+\gamma_{32}^1(\Phi)=0.
\end{equation*}
Similarly, with \eqref{r45}, we  have the coefficient of $w^4w^5$ in $\partial_{s_1}w^1$ is
\begin{equation*}
  \gamma_{45}^1(\Phi)+\gamma_{54}^1(\Phi)=0.
\end{equation*}
Hence, there is no $w^2w^3$ and $w^4w^5$ terms in $\partial_{s_1}w^1$, i.e.,
\begin{equation*}
  \partial_{s_1}w^1\thicksim(w^1)^2+\cancel{w^2w^3}+\cancel{w^4w^5}+\cdots.
\end{equation*}

For the formula of $\partial_{s_2}w^2$ in \eqref{y12}, by definition of $c_{23}^2$, $\gamma_{23}^2$, $\gamma_{45}^2$, $\gamma_{54}^2$ in \eqref{gamma}, we get the coefficient of $w^2w^3$ is
\begin{equation*}
  \begin{split}
    -c_{23}^2(\Phi)+\gamma_{23}^2(\Phi)&=-\frac{\lambda_3^2-\lambda_2^2}{2\lambda_2}-(\lambda_2-\lambda_3)\big[l^2 \cdot(\nabla_\Phi r_2 \cdot r_3)-l^2 \cdot(\nabla_\Phi r_3 \cdot r_2)\big]\\
    &=(\lambda_2-\lambda_3)\Big(1+\frac{\lambda_3^2-\lambda_2^2}{4\lambda_2^2}\Big),
  \end{split}
\end{equation*}
and the coefficient of $w^4w^5$ is
\begin{equation*}
  \gamma_{45}^2(\Phi)+\gamma_{54}^2(\Phi)=-(\lambda_4-\lambda_5)\big[l^2 \cdot(\nabla_\Phi r_4 \cdot r_5)-l^2 \cdot(\nabla_\Phi r_5 \cdot r_4)\big]=\frac{(\lambda_3-\lambda_2)^3}{4\lambda_2^2}.
\end{equation*}
Since $\lambda_2(0)=\lambda_3(0)=c_2$, it holds that $\lambda_2(\Phi)-\lambda_3(\Phi)$ is small for $\Phi\in B_{2\delta}^6(0)$ with small $\delta$. Combing with Lemma \ref{GN}, we have
\begin{equation*}
  \partial_{s_2}w^2\thicksim\cancel{(w^2)^2}+(\lambda_2-\lambda_3)w^2w^3+(\lambda_2-\lambda_3)w^4w^5+\cdots.
\end{equation*}

Similarly, with \eqref{right}, \eqref{left}, \eqref{gradlam}, \eqref{r23} and \eqref{r45}, we obtain
\begin{equation*}
\begin{split}
  -c_{32}^3(\Phi)&=\gamma_{32}^3(\Phi)=0, \quad \gamma_{45}^3(\Phi)=\gamma_{54}^3(\Phi)=0,\\
  -c_{45}^4(\Phi)&=\gamma_{45}^4(\Phi)=0, \quad \gamma_{23}^4(\Phi)=\gamma_{32}^4(\Phi)=0,\\
\end{split}\end{equation*}
\begin{equation*}
  \begin{split}
    -c_{54}^5(\Phi)+\gamma_{54}^5(\Phi)&=(\lambda_2-\lambda_3)(-\frac{\lambda_3}{\lambda_2}+\frac{\lambda_3^2-\lambda_2^2}{4\lambda_2^2}),\\
    \gamma_{23}^5(\Phi)+\gamma_{32}^5(\Phi)&=\frac{(\lambda_2-\lambda_3)^3}{4\lambda_2^2},
  \end{split}
\end{equation*}
\begin{equation*}
\gamma_{23}^6(\Phi)=\gamma_{32}^6(\Phi)=0,\quad \gamma_{45}^6(\Phi)=\gamma_{54}^6(\Phi)=0.
\end{equation*}
This means that expressions of $\partial_{s_i}w^i$ admit the following property: when $i=1,3,4,6$, terms of $w^2w^3$ and $w^4w^5$ in $\partial_{s_i}w^i$ disappear; when $i=2,5$, the coefficients in front of $w^2w^3$ and $w^4w^5$ always have a factor $\lambda_2(\Phi)-\lambda_3(\Phi)$.
By \eqref{219}, for small $\varepsilon\in(0,\frac1{100}]$, we can choose sufficiently small $\delta$ such that
\begin{equation}\label{347}
  \begin{split}
    |\lambda_2(\Phi)-\lambda_3(\Phi)|\leq\varepsilon,\quad\text{for}\quad \Phi\in B_{2\delta}^6(0).
  \end{split}
\end{equation}
 \hfill$\Box$

In summary, we list the important structures\footnote{Here we only list the terms, whose nonlinear interactions might be large, because they are from the same or almost-the-same eigenvalues.} of $\partial_{s_i}w^i$ here:
\begin{equation}\label{infmw}
 \left\{ \begin{split}
    \partial_{s_1}w^1\thicksim&(w^1)^2+\cancel{w^2w^3}+\cancel{w^4w^5}+\cdots,\\
    \partial_{s_2}w^2\thicksim&\cancel{(w^2)^2}+(\lambda_2-\lambda_3)w^2w^3+(\lambda_2-\lambda_3)w^4w^5+\cdots,\\
    \partial_{s_3}w^3\thicksim&(w^3)^2+\cancel{w^2w^3}+\cancel{w^4w^5}+\cdots,\\
    \partial_{s_4}w^4\thicksim&(w^4)^2+\cancel{w^2w^3}+\cancel{w^4w^5}+\cdots,\\
    \partial_{s_5}w^5\thicksim&\cancel{(w^5)^2}+(\lambda_2-\lambda_3)w^2w^3+(\lambda_2-\lambda_3)w^4w^5+\cdots,\\
    \partial_{s_6}w^6\thicksim&(w^6)^2+\cancel{w^2w^3}+\cancel{w^4w^5}+\cdots.\\
  \end{split}\right.
\end{equation}
The deleted terms mean their coefficients are \underline{zero}. That is to say that no interactions of the almost-repeated characteristic waves ($w^2$ and $w^3$, $w^4$ and $w^5$) appear in the equations of $\{\partial_{s_i}w^i\}_{i=1,3,4,6}$. Moreover, we have $\lambda_2-\lambda_3=\lambda_4-\lambda_5$ being small, and it shows that the interactions $w^2w^3$, $w^4w^5$ in $\partial_{s_2}w^2$ and  $\partial_{s_5}w^5$ are also weak. This allows us to use the bi-characteristic transformation to get desired bounds.

Back to \eqref{y14}-\eqref{y13} with the above coefficients, we also derive the equations for $\{\rho_i\}$ and $\{v^i\}$:
\begin{equation}\label{infmrho}
 \left\{ \begin{split}
    \partial_{s_1}\rho_1\thicksim&v^1+\rho_1(w^i)_{i\neq 1},\\
    \partial_{s_2}\rho_2\thicksim&\cancel{v^2}+(\lambda_2-\lambda_3)\rho_2w^3+\rho_2(w^i)_{i\neq 2,3},\\
    \partial_{s_3}\rho_3\thicksim& v^3+\cancel{\rho_3 w^2}+\rho_3(w^i)_{i\neq 2,3},\\
    \partial_{s_4}\rho_4\thicksim& v^4+\cancel{\rho_4 w^5}+\rho_4(w^i)_{i\neq 4,5},\\
    \partial_{s_5}\rho_5\thicksim&\cancel{v^5}+(\lambda_2-\lambda_3)\rho_5w^4+\rho_5(w^i)_{i\neq 4,5},\\
    \partial_{s_6}\rho_6\thicksim&v^6+\rho_6(w^i)_{i\neq 6},\\
  \end{split}\right.
\end{equation}
and
\begin{equation}\label{infmv}
 \left\{ \begin{split}
    \partial_{s_1}v^1\thicksim&v^1(w^i)_{i\neq 1}+\rho_1(\cancel{w^2w^3}+\cancel{w^4w^5}+\cdots),\\
    \partial_{s_2}v^2\thicksim&v^2[(\lambda_2-\lambda_3)w^3+\cdots]+\rho_2[(\lambda_2-\lambda_3)w^4w^5+\cdots],\\
    \partial_{s_3}v^3\thicksim&v^3(\cancel{w^2}+\cdots)+\rho_3(\cancel{w^4w^5}+\cdots),\\
    \partial_{s_4}v^4\thicksim&v^4(\cancel{w^5}+\cdots)+\rho_4(\cancel{w^2w^3}+\cdots),\\
    \partial_{s_5}v^5\thicksim&v^5[(\lambda_2-\lambda_3)w^4+\cdots]+\rho_5[(\lambda_2-\lambda_3)w^2w^3+\cdots],\\
    \partial_{s_6}v^6\thicksim&v^6(w^i)_{i\neq 6}+\rho_6(\cancel{w^2w^3}+\cancel{w^4w^5}+\cdots).\\
  \end{split}\right.
\end{equation}
The deleted terms mean their coefficients are zero. And we only list the terms, whose nonlinear interactions might be large.

\section{Construction of Initial Data}\label{datanorm}
Let $\varepsilon\in(0,\frac1{100}]$ be a small parameter. Given a small fixed parameter $\eta$ $(0<\eta\ll 1$), we choose initial data $w^i_{(\eta)}(z_i,0)$ such that
$\text{supp}\ w^i_{(\eta)}(z_i,0)\in[\eta,2\eta]$,
and
\begin{equation}\label{W0}
  W_0^{(\eta)}:=\max_i\sup_{z_i}|w^i_{(\eta)}(z_i,0)|=w^1_{(\eta)}(z_0,0)>0.
\end{equation}
Furthermore, we require
\begin{align}\label{273}
  \sup_{z_6}|w^6(z_6,0)|\leq \frac{(1-\varepsilon)^4}{2(1+\varepsilon)^3}W_0^{(\eta)}
  \end{align}
and
  \begin{align}\label{initial}
   \max_{i=3,4}\sup_{z_i}|w^i_{(\eta)}(z_i,0)|\leq \min\Big\{\frac{(1-\varepsilon)^4|c_{11}^1(0)|}{(1+\varepsilon)^3}W_0^{(\eta)},W_0^{(\eta)}\Big\}.
\end{align}
For $w^1_{(\eta)}(z,0)$, we choose
\begin{equation} \label{data}
   w^1_{(\eta)}(z,0)=\theta\int_\mathbb{R}\zeta_{\frac\eta{10}}(y)| \ln (z-y)|^\alpha \chi(z-y)dy, \quad 0<\alpha<\frac12,
\end{equation}
where $\theta$ is a small parameter to be determined later. Here $\chi$ is the characteristic function
\begin{equation*}
  \chi(z)=\left\{\begin{array}{ll}
  1,\quad z\in [\frac65\eta,\frac{9}5\eta],\\
  0,\quad z\notin(\frac65\eta,\frac{9}5\eta),
  \end{array}\right.
\end{equation*}
and $\zeta_{\frac\eta{10}}(z)$ is a test function in $C_0^\infty(\mathbb{R})$ and satisfying
\begin{equation*}
  \text{supp}\ \zeta_{\frac\eta{10}}(z)=\{z:|z|\leq\frac\eta{10}\},\quad\text{and}\quad\int_{\mathbb{R}}\zeta_{\frac\eta{10}}(z)dz=1.
\end{equation*}

We then define the following norms for $i =1,2,3,4,5,6,$
\begin{align}
  S_i(t):=&\sup_{(z'_i,s'_i)\atop z'_i\in[\eta,2\eta],\ 0\leq s'_i\leq t}\rho_i(z'_i,s'_i),&S(t):=&\max_{i}S_i(t),\label{46}\\
  J_i(t):=&\sup_{(z'_i,s'_i)\atop z'_i\in[\eta,2\eta]\ 0\leq s'_i\leq t}|v^i(z'_i,s'_i)|,&J(t):=&\max_{i}J_i(t),\label{47}\\
  W(t):=&\max_i\sup_{(x',t')\atop 0\leq t'\leq t}|w^i(x',t')|,&\bar{U}(t):=&\sup_{(x',t')\atop 0\leq t'\leq t}|\Phi(x',t')|.\label{48}
\end{align}
We also denote
\begin{align*}
  & V_1(t):=\sup_{(x',t')\notin\mathcal{R}_1,\atop 0\leq t'\leq t}|w^1(x',t')|,\\
  &V_{\bar{2}}(t):=\max\sup_{(x',t')\notin\mathcal{R}_2\cup\mathcal{R}_3,\atop 0\leq t'\leq t}\{|w^2(x',t')|,|w^3(x',t')|\},\\
 &V_{\bar{5}}(t):=\max\sup_{(x',t')\notin\mathcal{R}_4\cup\mathcal{R}_5,\atop 0\leq t'\leq t}\{|w^4(x',t')|,|w^5(x',t')|\},\\
  & V_6(t):=\sup_{(x',t')\notin\mathcal{R}_6,\atop 0\leq t'\leq t}|w^6(x',t')|
\end{align*}
and
\begin{equation}\label{410}
V(t):=\max_{i}V_i(t),\quad \text{for} \ i=1,\bar{2},\bar{5},6.
\end{equation}
Finally, we set
\begin{equation*}
  \underline{S}(t):=\min_{i\in\{2,\cdots,6\}}\inf_{(z'_i,s'_i)\atop z'\in[\eta,2\eta],\ 0\leq s'_i\leq t} \rho_i(z'_i,s'_i).
\end{equation*}

\section{Estimates of Norms}\label{3d}
Recall that for the aforementioned system \eqref{y2.1}, two pairs of characteristic wave speeds are almost the same,
      \begin{equation*}
   \lambda_2(\Phi)\thickapprox\lambda_3(\Phi),\ \lambda_4(\Phi)\thickapprox\lambda_5(\Phi),
      \end{equation*}
   for $\Phi\in B_{2\delta}^6(0)$. This means two pairs of characteristic strips $\mathcal{R}_2$ and $\mathcal{R}_3$ or $\mathcal{R}_4$ and $\mathcal{R}_5$ could overlap for a long time. Our strategy is to consider the wave propagations in four characteristic strips:
 \begin{equation}\label{4strip2}
\{\mathcal{R}_1,\mathcal{R}_2\cup\mathcal{R}_3,\mathcal{R}_4\cup\mathcal{R}_5,\mathcal{R}_6\}.
  \end{equation}
   These four strips will be completely separated when $t>t_0^{(\eta)}$, where $t_0^{(\eta)}$ can be precisely calculated as in \eqref{t0} below.

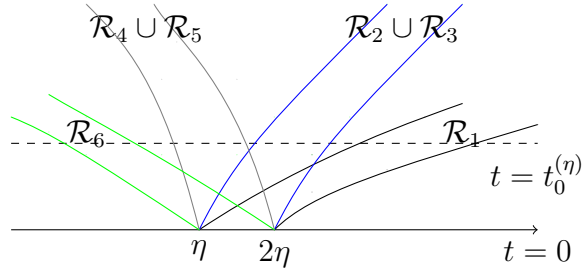
\begin{center}
\begin{tikzpicture}
\draw[->](0,0)--(7,0)node[left,below]{$t=0$};
\draw[dashed](0,1.15)--(7,1.15)node[right,below]{$t=t_0^{(\eta)}$};
\filldraw [black] (3.5,0) circle [radius=0.01pt]
(4,0.5) circle [radius=0.01pt]
(5,0.8) circle [radius=0.01pt]
(7,1.4) circle [radius=0.01pt];
\draw (3.5,0)..controls (4,0.5) and (5,0.8)..(7,1.4);
\node [below]at(3.5,0){$2\eta$};

\filldraw [black] (2.5,0) circle [radius=0.01pt]
(4,1) circle [radius=0.01pt]
(5.5,1.5) circle [radius=0.01pt]
(6,1.68) circle [radius=0.01pt];
\draw (2.5,0)..controls (4,1) and (5.5,1.5)..(6,1.68);
\node [below]at(2.5,0){$\eta$};
\node [below] at(6,1.6){$\mathcal{R}_1$};

\filldraw [green] (3.5,0) circle [radius=0.01pt]
(4,1) circle [radius=0.01pt]
(4.5,1.5) circle [radius=0.01pt]
(6,3) circle [radius=0.01pt];
\draw [color=blue](3.5,0)..controls (4,1) and (4.5,1.5)..(6,3);

\filldraw [green] (2.5,0) circle [radius=0.01pt]
(3,1) circle [radius=0.01pt]
(3.5,1.5) circle [radius=0.01pt]
(5,3) circle [radius=0.01pt];
\draw [color=blue](2.5,0)..controls (3,1) and (3.5,1.5)..(5,3);
\node [below] at(5.2,3){$\mathcal{R}_2\cup\mathcal{R}_3$};

\filldraw [gray] (3.5,0) circle [radius=0.01pt]
(3,2) circle [radius=0.01pt]
(2.2,2.5) circle [radius=0.01pt]
(1.9,3) circle [radius=0.01pt];
\draw [color=gray](3.5,0)..controls (3,2) and (2.2,2.5)..(1.9,3);
\node [below] at(1.8,3){$\mathcal{R}_4\cup\mathcal{R}_5$};

\filldraw [gray] (2.5,0) circle [radius=0.01pt]
(2,2) circle [radius=0.01pt]
(1.5,2.5) circle [radius=0.01pt]
(1,3) circle [radius=0.01pt];
\draw [color=gray] (2.5,0)..controls (2,2) and (1.5,2.5)..(1,3);

\filldraw [green] (3.5,0) circle [radius=0.01pt]
(2,1) circle [radius=0.01pt]
(1,1.5) circle [radius=0.01pt]
(0.5,1.8) circle [radius=0.01pt];
\draw [color=green](3.5,0)..controls (2,1) and (1,1.5)..(0.5,1.8);
\node [below] at(1,1.6){$\mathcal{R}_{6}$};

\filldraw [green] (2.5,0) circle [radius=0.01pt]
(1,1) circle [radius=0.01pt]
(0.5,1.3) circle [radius=0.01pt]
(0,1.5) circle [radius=0.01pt];
\draw [color=green] (2.5,0)..controls (1,1) and (0.5,1.3)..(0,1.5);
\end{tikzpicture}
\end{center}

For notational simplicity, we denote
\begin{equation*}
  \mathcal{R}_{\bar{2}}:=\mathcal{R}_2\bigcup\mathcal{R}_3,\quad \mathcal{R}_{\bar{5}}:=\mathcal{R}_4\bigcup\mathcal{R}_5.
\end{equation*}
Let
\begin{equation*}
  \bar{\lambda}_i:=\sup_{\Phi\in B_{2\delta}^6(0)}\lambda_i(\Phi),\quad \underline{\lambda}_i:=\inf_{\Phi\in B_{2\delta}^6(0)}\lambda_i(\Phi),\quad \text{for}\quad i=1,6,
\end{equation*}
\begin{equation*}
  \bar{\lambda}_{\bar{2}}:=\sup_{\Phi\in B_{2\delta}^6(0)}\{\lambda_2(\Phi),\lambda_3(\Phi)\},\quad \underline{\lambda}_{\bar{2}}:=\inf_{\Phi\in B_{2\delta}^6(0)}\{\lambda_2(\Phi),\lambda_3(\Phi)\}
\end{equation*}
\begin{equation*}
\bar{\lambda}_{\bar{5}}:=\sup_{\Phi\in B_{2\delta}^6(0)}\{\lambda_4(\Phi),\lambda_5(\Phi)\},\quad \underline{\lambda}_{\bar{5}}:=\inf_{\Phi\in B_{2\delta}^6(0)}\{\lambda_4(\Phi),\lambda_5(\Phi)\}
\end{equation*}
and
\begin{equation*}
  \sigma:=\min_{\alpha<\beta\atop \alpha,\beta\in\{1,\bar{2},\bar{5},6\}}(\underline{\lambda}_\alpha-\bar{\lambda}_\beta).
\end{equation*}
According to \eqref{y3}, with $\delta$ sufficiently small we have that
$\sigma$ has a uniform positive lower bound.
By the defining equation of characteristics \eqref{flow}, for $\alpha\in\{1,\bar{2},\bar{5},6\}$, $z\in[\eta,2\eta]$, we have
\begin{equation*}
  z+\underline{\lambda}_\alpha t\leq X_\alpha(z,t)\leq z+\bar{\lambda}_\alpha t.
\end{equation*}
And for all $\alpha<\beta$, with $\alpha,\beta\in\{1,\bar{2},\bar{5},6\}$, it holds
\begin{equation*}
  \begin{split}
   X_\alpha(\eta,t)-X_\beta(2\eta,t)\geq (\eta+\underline{\lambda}_\alpha t)-(2\eta+\bar{\lambda}_\beta t)=-\eta+(\underline{\lambda}_\alpha-\bar{\lambda}_\beta)t\geq-\eta+\sigma t.
  \end{split}
\end{equation*}
Note that the above difference is strictly positive when
\begin{equation}\label{t0}
  t>t_0^{(\eta)}:=\frac\eta\sigma.
\end{equation}
This implies that the four characteristic strips in \eqref{4strip2} are well separated when $t>t_0^{(\eta)}$.

{\bf $\bullet$ Estimates for $t\in[0,t_0^{(\eta)}]$}

In the non-separated region before $t=t_0^{(\eta)}$, all the characteristic strips $\mathcal{R}_i$ are overlapped. But even then, the inverse foliation density of characteristics still obey positive lower bounds in this time region.

Let $\Gamma$ be the maximum of all the coefficients $\{c_{im}^i\}$ and $\{\gamma_{km}^i\}$ in \eqref{y14}-\eqref{y13}. According to the calculations in Section \ref{yf}, we have
\begin{equation}\label{order1}
  \Gamma=O(1).
\end{equation}
In the following part, we will bound $W(t)$, $V(t)$, $S(t)$, $J(t)$ and $\bar{U}(t)$ defined in \eqref{46}-\eqref{410}.
We estimate $W(t)$ first. For $t\in[0,t_0^{(\eta)}]$, by \eqref{y12}, we have
\begin{equation*}
  \frac{\partial}{\partial s_i}|w^i|\leq \Gamma W^2.
\end{equation*}
Comparing with solutions to
\begin{equation*}
  \left\{\begin{array}{ll}
  \frac{d}{dt}Y=\Gamma Y^2,\\
  Y(0)=W_0^{(\eta)},
  \end{array}
  \right.
\end{equation*}
we have
\begin{equation}\label{wt0}
  |w^i|\leq Y(t)=\frac{W_0^{(\eta)}}{1-\Gamma W_0^{(\eta)} t},\quad \text{for}\quad t<\min\Big\{\frac1{\Gamma W_0^{(\eta)}},T\Big\}.
\end{equation}
With \eqref{t0} and \eqref{order1}, we obtain
\begin{equation}\label{gwt}
  \Gamma W_0^{(\eta)} t_0=O(\eta W_0^{(\eta)}).
\end{equation}
Applying \eqref{gwt} to \eqref{wt0}, it holds
\begin{equation*}
  |w^i(x,t)|\leq (1+\varepsilon)W_0^{(\eta)},\quad \forall\ x\in\mathbb{R},\ t\in[0,t_0^{(\eta)}],
\end{equation*}
for some $\varepsilon>0$ small. This could be achieved by requiring $\theta$ in \eqref{data} sufficiently small.
Back to the definition in \eqref{48}, this implies that
\begin{equation}\label{y15}
  |W(t)|\leq (1+\varepsilon)W_0^{(\eta)},\quad \forall\  t\in[0,t_0^{(\eta)}].
\end{equation}

We proceed to bound $V(t)$. Any $(x',t')\notin\mathcal{R}_i$ can be characterised by the corresponding characteristic coordinates $(z_i',s_i')$ satisfying $z'_i\notin[\eta,2\eta]$. Since our constructed initial data are supported in $[\eta,2\eta]$, for $z'_i\notin[\eta,2\eta]$ we have $w^i_{(\eta)}(z_i,0)=0$. Integrating  \eqref{y12} along the characteristic $C_i$, we obtain
\begin{equation*}
  V(t)=O(\int_0^{t_0^{(\eta)}}w^iw^jds_i)=O(\eta [W(t)]^2)=O(\eta [W_0^{(\eta)}]^2).
\end{equation*}

We then estimate $S(t)$. From \eqref{y14}, we have
\begin{equation*}
  \frac{\partial\rho_i}{\partial s_i}=O(\rho_iW).
\end{equation*}
Integrating the above equation along $\mathcal{C}_i$, we get
\begin{equation}\label{rho i}
  \rho_i(z_i,t)=\rho_i(z_i,0)\exp\big(O(tW(t))\big).
\end{equation}
Note that by definitions of \eqref{flow} and \eqref{dense} we have
\begin{equation}\label{319}
  \rho_i(z_i,0)=1,
\end{equation}
then via (\ref{rho i}) it holds $\rho_i(z_i,t)>0$.
Moreover, by \eqref{y15}, we obtain
\begin{equation}\label{320}
  \rho_i(z_i,t)=\exp\big(O(\eta W_0^{(\eta)})\big),\quad\forall\ t\in[0,t_0^{(\eta)}].
\end{equation}
For $\eta$ being small, we can choose sufficiently small $\theta$ such that
\begin{equation}\label{321}
  1-\varepsilon\leq\exp\big(O(\eta W_0^{(\eta)})\big)\leq1+\varepsilon.
\end{equation}
Inserting \eqref{321} into \eqref{320}, we get
\begin{equation*}
 1-\varepsilon\leq\rho_i(z_i,t)\leq1+\varepsilon,\quad \forall\ t\in[0,t_0^{(\eta)}].
\end{equation*}
So we have
\begin{equation}\label{y16}
  S(t)=O(1),\quad \forall\ t\in[0,t_0^{(\eta)}].
\end{equation}

For $J(t)$, from \eqref{y13}, we have
\begin{equation*}
\frac{\partial v^i}{\partial s_i}=O(S(t)[W(t)]^2).
\end{equation*}
Using \eqref{y15} and \eqref{y16}, we get
\begin{equation}
  J(t)=O(W_0^{(\eta)}+t[W(t)]^2)=O(W_0^{(\eta)}+\eta [W_0^{(\eta)}]^2)=O(W_0^{(\eta)}),\quad \forall\ t\in[0,t_0^{(\eta)}].
\end{equation}

Next, we give an estimate for $\bar{U}$. By \eqref{uxx}, we have
\begin{equation}\label{phi}
  \Phi(x,t)=\int_{X_6(\eta,t)}^x\frac{\partial \Phi(x',t)}{\partial x} dx'=\int_{X_6(\eta,t)}^x\sum_{k}w^kr_k(x',t)dx'.
\end{equation}
since
\begin{equation*}
  r_k(\Phi)=O(1),
\end{equation*}
equality \eqref{phi} implies
\begin{equation}\label{y19}
  |\Phi(x,t)|=O\Big(\sum_{k}\int_{X_6(\eta,t)}^{X_1(2\eta,t)|}|w^k(x',t)|dx'\Big).
\end{equation}
By \eqref{y15} and definition \eqref{48}, we get
\begin{equation*}
   \bar{U}(t)=O\big(W(t)(\eta+(\bar{\lambda}_1-\underline{\lambda}_6)t)\big)=O(\eta W_0^{(\eta)}),\quad \forall\ t\in[0,t_0^{(\eta)}].
\end{equation*}

In summary, we have proved
\begin{align}
  W(t)&=O(W_0^{(\eta)}), &\forall\  t\in[0,t_0^{(\eta)}],\label{531}\\
  V(t)&=O(\eta [W_0^{(\eta)}]^2), &\forall\  t\in[0,t_0^{(\eta)}],\\
   S(t)&=O(1), &\forall\ t\in[0,t_0^{(\eta)}],\\
   J(t)&=O(W_0^{(\eta)}), &\forall\ t\in[0,t_0^{(\eta)}],\\
   \bar{U}(t)&=O(\eta W_0^{(\eta)}),&\forall\ t\in[0,t_0^{(\eta)}].
\end{align}

{\bf $\bullet$ Estimates for  $t\in[t_0^{(\eta)},T]$}

Though six characteristic strips only separate partially because of the non-strict hyperbolicity, the aforementioned four characteristic strips \eqref{4strip2} are well separated after $t>t_0^{(\eta)}$. In the following part we estimate geometric quantities in different strips, and we will show that a shock forms in $\mathcal{R}_1$ along $\mathcal{C}_1$. Subtle structures of \eqref{infmw}-\eqref{infmv} yield cancellations of some potentially dangerous terms.

Based on the bounds for $t\in[0,t_0^{(\eta)}]$, we carry on estimates in the region $t\in[t_0^{(\eta)},T]$.
We first estimate $S(t)$, i.e., the supremum of inverse foliation densities. For $\alpha=1$ or $6$, if $(x,t)\in \mathcal{R}_\alpha$, we have
\begin{equation}\label{335}
  \frac{\partial \rho_\alpha}{\partial s_\alpha}=O(J_\alpha+ V S_\alpha),\quad \alpha=1\ \text{or}\ 6.
\end{equation}
Thus, by integrating \eqref{335} along the characteristic $\mathcal{C}_\alpha$, we conclude
\begin{equation*}
  S_\alpha(t)=O(1+tJ_\alpha+ tV S_\alpha),\quad \alpha=1\ \text{or}\ 6.
\end{equation*}
For $\rho_2$, if $(x,t)\in \mathcal{R}_2$, the characteristic $\mathcal{C}_2$ crossing $(x,t)$ may also intersect $\mathcal{R}_3$. By a crucial structure displayed in \eqref{y14}
\begin{equation*}
  c_{23}^2=(\lambda_2-\lambda_3)O(1),
\end{equation*}
we have
\begin{equation}\label{338}
  \frac{\partial \rho_2}{\partial s_2}=O( V S_2+(\lambda_2-\lambda_3)w^3\rho_2).
\end{equation}
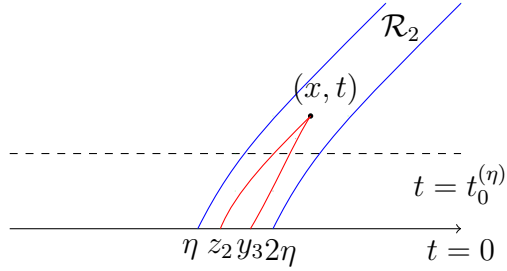
\begin{center}
\begin{tikzpicture}
\draw[->](0,0)--(6,0)node[left,below]{$t=0$};
\draw[dashed](0,1)--(6,1)node[left,below]{$t=t_0^{(\eta)}$};
\node [below]at(3.6,0){$2\eta$};
\node [below]at(2.4,0){$\eta$};
\filldraw [green] (3.5,0) circle [radius=0.01pt]
(4,1) circle [radius=0.01pt]
(4.5,1.5) circle [radius=0.01pt]
(6,3) circle [radius=0.01pt];
\draw [color=blue](3.5,0)..controls (4,1) and (4.5,1.5)..(6,3);

\filldraw [green] (2.5,0) circle [radius=0.01pt]
(3,1) circle [radius=0.01pt]
(3.5,1.5) circle [radius=0.01pt]
(5,3) circle [radius=0.01pt];
\draw [color=blue](2.5,0)..controls (3,1) and (3.5,1.5)..(5,3);
\node [below] at(5.2,3){$\mathcal{R}_2$};

\filldraw [green] (2.8,0) circle [radius=0.01pt]
(3,0.5) circle [radius=0.01pt]
(3.5,1) circle [radius=0.01pt]
(4,1.5) circle [radius=0.01pt];
\draw [color=red](2.8,0)..controls (3,0.5) and (3.5,1)..(4,1.5);
\node [below]at(2.8,0){$z_2$};

\filldraw [black] (3.2,0) circle [radius=0.01pt]
(3.5,0.5) circle [radius=0.01pt]
(3.8,1.2) circle [radius=0.01pt]
(4,1.5) circle [radius=0.8pt];
\draw [color=red](3.2,0)..controls (3.5,0.5) and (3.8,1.2)..(4,1.5);
\node [below]at(3.2,0){$y_3$};
\node [above]at(4.2,1.5){$(x,t)$};
\end{tikzpicture}
\end{center}
Applying \eqref{y17} and using the fact that $dz_2=0$, we have
\begin{equation}\label{y18}
\begin{split}
  &\int_0^t(\lambda_2-\lambda_3)w^3\big(X_2(z_2,t'),t'\big)\rho_2\big(X_2(z_2,t'),t'\big)dt'\\
  =&O\Big(\int_{[\eta,2\eta]}(\lambda_2-\lambda_3)\frac{\rho_3}{\lambda_2-\lambda_3}w^3\rho_2dy_3\Big)=O(\eta S_2J_3).
  \end{split}
\end{equation}
Hence, by integrating \eqref{338} along the characteristic $\mathcal{C}_2$, we get
\begin{equation}\label{540}
  S_2=O(1+tVS_2+\eta S_2J_3).
\end{equation}
For $\rho_3$, since we have
\begin{equation*}
  c_{32}^3=0,
\end{equation*}
it follows
\begin{equation}\label{542}
  S_3=O(1+tJ_3+tVS_3).
\end{equation}
For $\rho_4$ and $\rho_5$, if $(x,t)\in \mathcal{R}_{\bar{5}}$, we proceed as for $\rho_3$ and $\rho_2$, respectively. Note
\begin{equation*}
  c_{45}^4=0,\quad \text{and}\quad
  c_{54}^5=(\lambda_4-\lambda_5)O(1).
\end{equation*}
As estimates obtained in \eqref{542}, we have
\begin{equation*}
  S_4=O(1+t J_4+tVS_4).
\end{equation*}
And proceed as in \eqref{338}-\eqref{540}, we also obtain
\begin{equation}
  S_5=O(1+tVS_5+\eta S_5J_4).
\end{equation}
So we have
\begin{equation}
    S=O(1+tVS+tJ+\eta SJ).
\end{equation}

We then bound $J(t)$, the supremum of $\{v^i\}_{i=1,\cdots,6}$. For $\alpha=1,6$, and $(x,t)\in \mathcal{R}_\alpha$, with the cancellations of the deleted terms in \eqref{infmv}, we have
\begin{equation*}
  \frac{\partial v^\alpha}{\partial s_\alpha}=O(VJ_\alpha+ V^2 S_\alpha).
\end{equation*}
Thus,
\begin{equation}\label{548}
  J_\alpha(t)=O(W_0^{(\eta)}+tVJ_\alpha+ tV^2 S_\alpha),\quad \text{for}\ \alpha=1,6.
\end{equation}
For $v^3$, $(x,t)\in \mathcal{R}_3$, we have
\begin{equation*}
  \frac{\partial v^3}{\partial s_3}=O\Big(VJ_3+S_3 V^2+VS_3w^2\Big).
\end{equation*}
Since
\begin{equation*}
  \int_0^tVS_3w^2dt'=\int_0^tVS_3\frac{v^2}{\rho_2}dt'=O\Big(\frac{tV S_3J}{\underline{S}}\Big),
\end{equation*}
we obtain
\begin{equation}\label{551}
  J_3(t)=O\Big(W_0^{(\eta)}+tVJ_3+ tV^2 S_3+\frac{tV S_3J}{\underline{S}}\Big).
\end{equation}
For $v^4$, in the same fashion we get
\begin{equation}\label{552}
  J_4(t)=O\Big(W_0^{(\eta)}+tVJ_4+ tV^2 S_4+\frac{tV S_4J}{\underline{S}}\Big).
\end{equation}
For $v^2$, $(x,t)\in \mathcal{R}_2$, we have
\begin{equation*}
  \frac{\partial v^2}{\partial s_2}=O((\lambda_2-\lambda_3)w^3J_2+V^2 S_2).
\end{equation*}
By a similar estimate to \eqref{y18}, it holds
\begin{equation}\label{553}
\begin{split}
  &\int_0^t(\lambda_2-\lambda_3)w^3\big(X_2(z_2,t'),t'\big)J_2\big(X_2(z_2,t'),t'\big)dt'\\
  =&O\Big(\int_{[\eta,2\eta]}(\lambda_2-\lambda_3)\frac{\rho_3}{\lambda_2-\lambda_3}w^3J_2dy_3\Big)\\
  =&O\Big(J_2\int_{[\eta,2\eta]}v^3dy_3\Big)=O(\eta J_2J_3).
  \end{split}
\end{equation}
Integrating along characteristic $\mathcal{C}_2$,  we arrive at
\begin{equation}\label{555}
  J_2=O(W_0^{(\eta)}+\eta J_2J_3+tV^2 S_2).
\end{equation}
Analogously, for $v^5$ we obtain
\begin{equation}\label{556}
  J_5=O(W_0^{(\eta)}+\eta J_4J_5+tV^2 S_5).
\end{equation}
In conclusion, from \eqref{548}, \eqref{551}-\eqref{552}, and \eqref{555}-\eqref{556}, we have
\begin{align}
    J=O\Big(W_0^{(\eta)}+t VJ+tV^2 S+\eta J^2+\frac{tV SJ}{\underline{S}}\Big)).
\end{align}

We next bound $V(t)$, i.e., the supremum of $w^i$ outside the corresponding characteristic strip $\mathcal{R}_i$.
We first estimate $w^1(x,t)$ with $(x,t)\notin\mathcal{R}_1$ and $w^6(x,t)$ with $(x,t)\notin\mathcal{R}_6$. Let $i=\{1,6\}$. From \eqref{infmw}, we have
\begin{equation}\label{558}
  \frac{\partial w^i}{\partial{s_i}}=O(V^2)+O\Big(\sum_{k\neq i}w^k\Big)V+O\Big(\sum_{m\neq i,k\neq i\atop m\neq k}w^mw^k\Big).
\end{equation}
Note that $\mathcal{C}_i$ starts from $z_i\notin[\eta,2\eta]$ and ends at $(x,t)\notin\mathcal{R}_i$. When $t'\geq t_0^{(\eta)}$, for any point $\big(X_i(z_i,t'),t'\big)$ on $\mathcal{C}_i$, it holds either $\big(X_i(z_i,t'),t'\big)\in\big(\mathbb{R}\times[t_0^{(\eta)},t]\big)\setminus\bigcup_{k}\mathcal{R}_k$ or $\big(X_i(z_i,t'),t'\big)\in\mathcal{R}_k$ for some $k\neq i$.
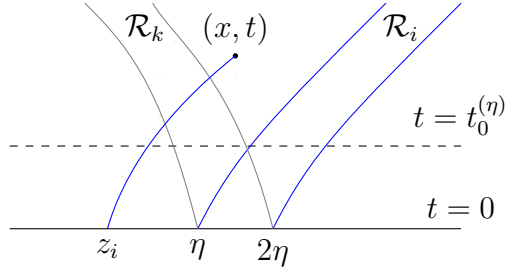
\begin{center}
\begin{tikzpicture}
\draw(0,0)--(6,0)node[left,above]{$t=0$};
\draw[dashed](0,1.1)--(6,1.1)node[left,above]{$t=t_0^{(\eta)}$};
\node [below]at(3.5,0){$2\eta$};
\node [below]at(2.5,0){$\eta$};

\filldraw [green] (3.5,0) circle [radius=0.01pt]
(4,1) circle [radius=0.01pt]
(4.5,1.5) circle [radius=0.01pt]
(6,3) circle [radius=0.01pt];
\draw [color=blue](3.5,0)..controls (4,1) and (4.5,1.5)..(6,3);

\filldraw [green] (2.5,0) circle [radius=0.01pt]
(3,1) circle [radius=0.01pt]
(3.5,1.5) circle [radius=0.01pt]
(5,3) circle [radius=0.01pt];
\draw [color=blue](2.5,0)..controls (3,1) and (3.5,1.5)..(5,3);
\node [below] at(5.2,3){$\mathcal{R}_{i}$};

\filldraw [gray] (3.5,0) circle [radius=0.01pt]
(3,2) circle [radius=0.01pt]
(2.2,2.5) circle [radius=0.01pt]
(1.9,3) circle [radius=0.01pt];
\draw [color=gray](3.5,0)..controls (3,2) and (2.2,2.5)..(1.9,3);
\node [below] at(1.8,3){$\mathcal{R}_{k}$};

\filldraw [gray] (2.5,0) circle [radius=0.01pt]
(2,2) circle [radius=0.01pt]
(1.5,2.5) circle [radius=0.01pt]
(1,3) circle [radius=0.01pt];
\draw [color=gray] (2.5,0)..controls (2,2) and (1.5,2.5)..(1,3);

\filldraw [black] (1.3,0) circle [radius=0.01pt]
(1.5,0.7) circle [radius=0.01pt]
(2,1.5) circle [radius=0.01pt]
(3,2.3) circle [radius=0.8pt];
\draw [color=blue](1.3,0)..controls (1.5,0.7) and (2,1.5)..(3,2.3);
\node[above]at (3,2.3){$(x,t)$};
\node[below]at (1.3,0){$z_i$};

\end{tikzpicture}
\end{center}
Let $I_k^i=\{t'\in[t_0^{(\eta)},t]:(x,t')\in\mathcal{C}_i\bigcap\mathcal{R}_k\}$  for $k\neq i$.
Integrating \eqref{558} along $\mathcal{C}_i$ and using $w^i_{(\eta)}(z_i,0)=0$, we have
\begin{equation}\label{559}
\begin{split}
    w^i(x,t)=&O\Big(tV^2+V\sum_{k\neq i}\int_0^tw^k\big(X_i(z_i,t'),t'\big)dt'\Big)\\
    =&O\Big(tV^2+V\sum_{k\neq i}\int_0^{t_0^{(\eta)}}w^k\big(X_i(z_i,t'),t'\big)dt'\Big)\\
    &+O\Big(V\sum_{k\neq i}\int_{t_0^{(\eta)}}^t w^k\big(X_i(z_i,t'),t'\big)dt'\Big)\\
    =&O\Big(tV^2+\eta[W_0^{(\eta)}]^2+V\sum_{k\neq i}\underbrace{\int_{I_k^i}w^k\big(X_i(z_i,t'),t'\big)dt'}_{M}\Big).
\end{split}
\end{equation}
Here  we use the fact $V(t)\leq W(t)=O(W_0^{(\eta)})$ for $t\leq t_0^{(\eta)}$. We proceed to bound $M$.
When $\big(X_i(z_i,t'),t'\big)\in I_k^i$ for some $k\neq i$, the picture is as below
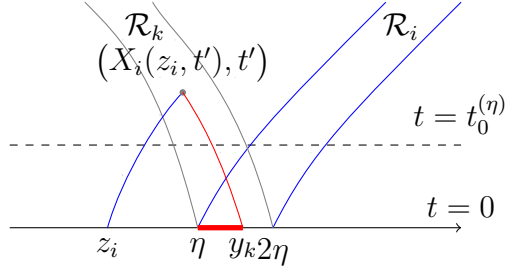
\begin{center}
\begin{tikzpicture}
\draw(0,0)--(2.5,0);
\draw[->](3.1,0)--(6,0)node[left,above]{$t=0$};
\draw[dashed](0,1.1)--(6,1.1)node[left,above]{$t=t_0^{(\eta)}$};
\node [below]at(3.5,0){$2\eta$};
\node [below]at(2.5,0){$\eta$};

\filldraw [green] (3.5,0) circle [radius=0.01pt]
(4,1) circle [radius=0.01pt]
(4.5,1.5) circle [radius=0.01pt]
(6,3) circle [radius=0.01pt];
\draw [color=blue](3.5,0)..controls (4,1) and (4.5,1.5)..(6,3);

\filldraw [green] (2.5,0) circle [radius=0.01pt]
(3,1) circle [radius=0.01pt]
(3.5,1.5) circle [radius=0.01pt]
(5,3) circle [radius=0.01pt];
\draw [color=blue](2.5,0)..controls (3,1) and (3.5,1.5)..(5,3);
\node [below] at(5.2,3){$\mathcal{R}_{i}$};

\filldraw [gray] (3.5,0) circle [radius=0.01pt]
(3,2) circle [radius=0.01pt]
(2.2,2.5) circle [radius=0.01pt]
(1.9,3) circle [radius=0.01pt];
\draw [color=gray](3.5,0)..controls (3,2) and (2.2,2.5)..(1.9,3);
\node [below] at(1.8,3){$\mathcal{R}_{k}$};

\filldraw [gray] (2.5,0) circle [radius=0.01pt]
(2,2) circle [radius=0.01pt]
(1.5,2.5) circle [radius=0.01pt]
(1,3) circle [radius=0.01pt];
\draw [color=gray] (2.5,0)..controls (2,2) and (1.5,2.5)..(1,3);

\filldraw [black] (1.3,0) circle [radius=0.01pt]
(1.5,0.7) circle [radius=0.01pt]
(2,1.5) circle [radius=0.01pt]
(2.3,1.8) circle [radius=0.8pt];
\draw [color=blue](1.3,0)..controls (1.5,0.7) and (2,1.5)..(2.3,1.8);
\node[above]at (2.3,1.8){$\big(X_i(z_i,t'),t'\big)$};
\node[below]at (1.3,0){$z_i$};

\filldraw [gray] (3.1,0) circle [radius=0.01pt]
(2.8,1) circle [radius=0.01pt]
(2.5,1.5) circle [radius=0.01pt]
(2.3,1.8) circle [radius=1pt];
\draw [color=red](3.1,0)..controls (2.8,1) and (2.5,1.5)..(2.3,1.8);
\node[below]at (3.1,0){$y_k$};

\draw[line width=2pt,color=red](2.5,0)--(3.1,0);
\end{tikzpicture}
\end{center}
we then employ the bi-characteristic coordinates and get
\begin{equation}\label{560}
  \begin{split}
    &\int_{I_k^i}w^k\big(X_i(z_i,t'),t'\big)dt'\\
    =&O\Big(\int_{y_k\in[\eta,2\eta]}\big|\frac{\rho_k\big(y_k,t'(y_i,y_k)\big)}{\lambda_i-\lambda_k}w^k\big(y_k,t'(y_i,y_k)\big)\big|dy_k\Big)\\
    =&O(\eta J_k).
  \end{split}
\end{equation}
Together with \eqref{559} and \eqref{560}, for $i=1,6$, we hence obtain
\begin{equation}\label{355}
  w^i(x,t)=O(tV^2+\eta [W_0^{(\eta)}]^2+\eta VJ),\quad\text{for}\ \forall\ (x,t)\notin\mathcal{R}_i.
\end{equation}
We next estimate $w^3(x,t)$ with $(x,t)\notin \mathcal{R}_2\cup\mathcal{R}_3$ and $w^4(x,t)$ with $(x,t)\notin \mathcal{R}_4\cup\mathcal{R}_5$. From \eqref{infmw}, we have that $w^4w^5$ and $w^2w^3$ vanish in $\partial_{s_3}w^3$ and $\partial_{s_4}w^4$. In the same fashion as for \eqref{559}, \eqref{560} and \eqref{355}, we hence obtain
\begin{equation}\label{563}
  w^i(x,t)=O(tV^2+\eta [W_0^{(\eta)}]^2+\eta VJ),\quad\text{for}\quad (x,t)\notin\mathcal{R}_i,\ i=3,4.
\end{equation}
We then deal with $w^2(x,t)$ for $(x,t)\notin \mathcal{R}_2\cup\mathcal{R}_3$. We discuss the following case first: $\mathcal{R}_2$ and $\mathcal{R}_3$ do not separate before time $t$. Note that $\lambda_2\geq\lambda_3$ by \eqref{211}-\eqref{212} and $a>b$. For $\mathcal{C}_2$ being a characteristic curve ending at $(x,t)$ and starting from $z_2\notin[\eta,2\eta]$, for any point $\big(X_2(z_2,t'),t'\big)$ on $\mathcal{C}_2$, the fact $\lambda_2\geq \lambda_3$ implies $\big(X_2(z_2,t'),t'\big)\notin\mathcal{R}_3$. Thus when $t'\in[t_0^{(\eta)},t]$, it holds that either $(X_2(z_2,t'),t'\big)$ stays in $\big(\mathbb{R}\times[t_0^{(\eta)},t]\big)\setminus\bigcup_{k}\mathcal{R}_k$ or it lies in one of the characteristic strips $\{\mathcal{R}_1,\mathcal{R}_6,\mathcal{R}_4\cup\mathcal{R}_5\}$. Integrating $\partial_{s_2}w^2$ along $\mathcal{C}_2$, we hence get
\begin{equation}\label{5.63}
\begin{split}
    w^2(x,t)=&O\Big(tV^2+V\sum_{k=1,6}\int_0^tw^k\big(X_2(z_2,t'),t'\big)dt'\Big)\\
    &+O\Big(\int_0^t\big|(\lambda_2-\lambda_3)w^4\big(X_2(z_2,t'),t'\big)w^5\big(X_2(z_2,t'),t'\big)\big|dt'\Big)\\
    =&O\Big(tV^2+\eta[W_0^{(\eta)}]^2+V\sum_{k=1,6}\int_{I_k^2}w^k\big(X_2(z_2,t'),t'\big)dt'\Big)\\
    &+O\Big(\int_{I_4^2\cup I_5^2}\big|(\lambda_2-\lambda_3)w^4\big(X_2(z_2,t'),t'\big)w^5\big(X_2(z_2,t'),t'\big)\big|dt'\Big).
\end{split}
\end{equation}
Similarly to \eqref{560}, it holds that
\begin{equation}\label{564}
  \begin{split}
    &\sum_{k=1,6}\int_{I_k^2}w^k\big(X_2(z_2,t'),t'\big)dt'\\
    =&O\Big(\sum_{k=1,6}\int_{y_k\in[\eta,2\eta]}\big|\frac{\rho_k\big(y_k,t'(y_2,y_k)\big)}{\lambda_2-\lambda_k}w^k\big(y_k,t'(y_2,y_k)\big)\big|dy_k\Big)\\
    =&O(\eta J).
  \end{split}
\end{equation}
Using bi-characteristic coordinates $(y_2,y_5)$ and together with \eqref{347}, we have
\begin{equation}\label{565}
\begin{split}
  &\int_{I_4^2\cup I_5^2}\big|(\lambda_2-\lambda_3)w^4\big(X_2(z_2,t'),t'\big)w^5\big(X_2(z_2,t'),t'\big)\big|dt'\\
  =&O\Big(\varepsilon\int_{y_5\in[\eta,2\eta]}\Big|w^4\big(y_5,t'(y_2,y_5)\big)w^5\big(y_5,t'(y_2,y_5)\big)\frac{\rho_5\big(y_5,t'(y_2,y_5)\big)}{\lambda_2-\lambda_5}dy_5\Big|\Big)\\
  =&O\Big(\varepsilon\int_{y_5\in[\eta,2\eta]}\Big|\frac{v^4 v^5}{\rho_4}dy_5\Big|\Big)=O\big(\eta\frac{\varepsilon}{\underline{S}}J^2\big).
\end{split}
\end{equation}
If $\mathcal{R}_2$ and $\mathcal{R}_3$ separate before $t$, for $(x,t)\notin \mathcal{R}_2\cup\mathcal{R}_3$, the characteristic curve $\mathcal{C}_2$ ending at $(x,t)$  may overlap $\mathcal{R}_3$ but does not intersect $\mathcal{R}_2$. In this case, we calculate the integration along $\mathcal{C}_2$ inside $\mathcal{R}_3$
\begin{equation}\label{5.42}
  \begin{split}
    &\int_{I_3^2}\big|(\lambda_2-\lambda_3)w^2\big(X_2(z_2,t'),t'\big)w^3\big(X_2(z_2,t'),t'\big)\big|dt'\\
    \leq&\int_{y_3\in[\eta,2\eta]}\Big|(\lambda_2-\lambda_3)w^2w^3\big(y_3,t'(y_2,y_3)\big)\frac{\rho_3\big(y_3,t'(y_2,y_3)\big)}{\lambda_2-\lambda_3}\Big|dy_3\\
    \leq&\int_{y_3\in[\eta,2\eta]}\Big|\frac{v^2}{\rho_2}v^3\big(y_3,t'(y_2,y_3)\big)\Big|dy_3\\
   =& O\big(\eta\frac{\varepsilon}{\underline{S}}J^2\big).
  \end{split}
\end{equation}

In sum, by \eqref{5.63}, \eqref{564}, \eqref{565} and \eqref{5.42}, we consequently obtain
\begin{equation}\label{357}
  w^2(x,t)=O\Big(\eta [W_0^{(\eta)}]^2+tV^2+\eta VJ+\eta\frac{\varepsilon}{\underline{S}}J^2\Big),\quad \forall\ (x,t)\notin \mathcal{R}_2\cup\mathcal{R}_3.
\end{equation}
Analogously, it also holds
\begin{equation}\label{358}
  w^5(x,t)=O\Big(\eta [W_0^{(\eta)}]^2+tV^2+\eta VJ+\eta\frac{\varepsilon}{\underline{S}}J^2\Big),\quad \forall\ (x,t)\notin \mathcal{R}_4\cup\mathcal{R}_5.
\end{equation}
Hence, from \eqref{355}-\eqref{563}, \eqref{357} and \eqref{358}, we get
\begin{equation}
  V=O\Big(\eta [W_0^{(\eta)}]^2+tV^2+\eta VJ+\eta\frac{\varepsilon}{\underline{S}}J^2\Big).
\end{equation}

Finally, we bound $\Phi$. If $(x,t)$ does not belong to any characteristic strip, we go back to \eqref{y19} and then obtain
\begin{equation}\label{u1}
 \bar{U}(t)=O\big((\eta+(\bar{\lambda}_1-\underline{\lambda}_6)t)V\big).
\end{equation}
If $(x,t)\in\mathcal{R}_k$ for some $k$, with characteristic coordinates we have
\begin{equation}\label{u2}
\begin{split}
  |\Phi(x,t)|=&\Big|\int_{X_k(\eta,t)}^x\frac{\partial \Phi(x',t)}{\partial x} dx'\Big|=\Big|\int_{X_k(\eta,t)}^x\sum_{k}w^kr_k(x',t)dx'\Big|\\
  \leq&\int_{X_k(\eta,t)}^{X_k(2\eta,t)}|w^k(x',t)|dx'=O\Big(\int_\eta^{2\eta}|w^k(x',t)|\rho_kdz_k\Big)\\
  =&O(\eta J).
  \end{split}
\end{equation}
Together with \eqref{u1} and \eqref{u2}, we get
\begin{equation}
   \bar{U}(t)=O(\eta J+\eta V+\eta tV),\quad \forall\ t\in[t_0^{(\eta)},T].
\end{equation}

In summary, for all $ t\in[0,T]$, we have
\begin{align}
    S=&O(1+tVS+tJ+\eta SJ),\label{y25}\\
    J=&O\Big(W_0^{(\eta)}+t VJ+tV^2 S+\eta J^2+\frac{tV SJ}{\underline{S}}\Big),\label{364}\\
    V=&O\Big(\eta [W_0^{(\eta)}]^2+tV^2+\eta VJ+\eta\frac{\varepsilon}{\underline{S}}J^2\Big), \label{366}\\
    \bar{U}=&O(\eta J+\eta V+\eta tV). \label{y31}
\end{align}

\noindent{\bf Remark:} The terms $\{\eta SJ,\eta J^2,\eta VJ,\frac{tV SJ}{\underline{S}},\eta\frac{\varepsilon}{\underline{S}}J^2\}$ in \eqref{y25}-\eqref{y31} are coming from system \eqref{y2} being non-strictly hyperbolic.
\section{Bootstrap argument}\label{bootargu}
We now design a bootstrap argument to bound $S$, $J$, $V$ and to improve assumption $|\Phi|\leq2\delta$.

When $t=0$, by \eqref{319} and definitions \eqref{47}-\eqref{48} we have
\begin{equation*}
 S(0)=1,\quad J(0)=W_0^{(\eta)},\quad V(0)=0.
\end{equation*}
Based on estimates \eqref{y25}-\eqref{366}, for a fixed constant $0<\kappa<\frac1{100}$ our goal is to prove
\begin{align}
  S(t)=&O(1),\label{63}\\
   J(t)=&O(W_0^{(\eta)}),\label{64}\\
   V(t)=&O(\eta  [W_0^{(\eta)}]^2+\eta\frac\varepsilon{\kappa}\theta^{-\frac13} [W_0^{(\eta)}]^2),\label{66}
\end{align}
for $t\in[0,T_\eta^*)$, where
\begin{equation}\label{368}
  T_\eta^*\leq\frac{C}{W_0^{(\eta)}}.
\end{equation}
And $C$ is a uniform constant. Once these estimates are achieved, we can improve bound for $\Phi$ and obtain
\begin{equation}\label{68}
\begin{split}
  &|\Phi|\leq\bar{U}(t)\\
  =&O(\eta W_0^{(\eta)}+\eta^2 W_0^{(\eta)}+\eta^2 [W_0^{(\eta)}]^2+\eta^2\frac\varepsilon\kappa \theta^{-\frac13}W_0^{(\eta)}+\eta^2\frac\varepsilon\kappa \theta^{-\frac13}[W_0^{(\eta)}]^2)\\
  =&O(\eta W_0^{(\eta)}+\eta^2\frac\varepsilon\kappa \theta^{-\frac13}W_0^{(\eta)}),\ t\in[0,T_\eta*].
  \end{split}
\end{equation}
Choosing $\theta$ being sufficiently small, we hence prove $\Phi\in B_\delta^6(0)$.

To obtain \eqref{63}-\eqref{66}, with $\theta$ small, we use bootstrap assumptions:
\begin{align}
  tV\leq& \theta^{\frac12},\label{y33}\\
   J\leq&\theta^{-\frac13} W_0^{(\eta)}\label{assum}\\
  \underline{S}:=&\min_{i=2,\cdots,6}\inf_{(z'_i,s'_i)\atop z'_i\in[\eta,2\eta],\ 0\leq s'_i\leq t} \rho_i(z'_i,s'_i)\geq \frac{\kappa}{2},\label{sbar}
\end{align}
where $0<\kappa<\frac1{100}$ is a uniform constant independent of $\theta$. In the argument below, we first improve bootstrap assumptions in \eqref{y33} and \eqref{assum}. An improvement of \eqref{sbar} for $\underline{S}$ will be given in Section \ref{lows}.

Applying \eqref{y33} and \eqref{assum} to \eqref{y25}, we have
\begin{align}\label{y34}
 S=O(1+tJ+\theta^{\frac12}S+\eta\theta^{-\frac13} W_0^{(\eta)} S)\quad
\Rightarrow\quad S=O(1+tJ).
\end{align}
For $J(t)$, we go back to \eqref{364}. With \eqref{y33}-\eqref{assum} we obtain
\begin{align}
   &J=O(W_0^{(\eta)}+\theta^{\frac12}J+\theta^{\frac12}V S+\eta\theta^{-\frac13} W_0^{(\eta)} J+\frac{\theta^{\frac12}}{\kappa}SJ)\notag\\
\Rightarrow\quad &J=O(W_0^{(\eta)}+\theta^{\frac12}V S+\theta^{\frac12}SJ).\label{y35}
\end{align}
Together with \eqref{y34}, \eqref{assum} and \eqref{y33}, equality \eqref{y35} implies
\begin{align}
 &J=O(W_0^{(\eta)}+\theta^{\frac12}V +\theta^{\frac12}tVJ+\theta^{\frac12}J +\theta^{\frac12}tJ^2)\notag\\
 &\ =O(W_0^{(\eta)}+\theta^{\frac12}V +\theta^{\frac14}J+\theta^{\frac12}J +\theta^{\frac16}J)\notag\\
 \Rightarrow \quad &J=O(W_0^{(\eta)}+\theta^{\frac12}V).\label{612}
\end{align}
For $V(t)$, employing \eqref{y33}-\eqref{assum} to \eqref{366}, we get
\begin{align*}
  V=O\Big(\eta [W_0^{(\eta)}]^2+\theta^{\frac12}V+\eta\theta^{-\frac13} W_0^{(\eta)} V+\eta\frac\varepsilon{\kappa}\theta^{-\frac13} W_0^{(\eta)}J\Big).
\end{align*}
Using \eqref{612}, we hence have
\begin{equation*}
 V=O\Big(\eta [W_0^{(\eta)}]^2+\eta\frac\varepsilon{\kappa}\theta^{-\frac13}[ W_0^{(\eta)}]^2+\eta\frac\varepsilon{\kappa}\theta^{\frac16} W_0^{(\eta)}V\Big).
\end{equation*}
This yields the desired bound of \eqref{66}
\begin{equation}\label{615}
 V=O\Big(\eta [W_0^{(\eta)}]^2+\eta\frac\varepsilon{\kappa}\theta^{-\frac13}[ W_0^{(\eta)}]^2\Big).
\end{equation}
Furthermore, by \eqref{368}, we have
\begin{equation}\label{tv}
 tV=O\Big(\eta W_0^{(\eta)}+\eta\frac\varepsilon{\kappa}\theta^{-\frac13} W_0^{(\eta)}\Big)=O\Big(\theta\eta(\ln\eta)^\alpha+\frac\varepsilon{\kappa}\theta^{\frac23}\eta(\ln\eta)^\alpha\Big)<O(\theta^{\frac12}),
\end{equation}
with $0<\alpha<\frac12$. Estimate \eqref{tv} improves bootstrap assumption \eqref{y33}. Back to \eqref{612} and \eqref{y34} with \eqref{615}, we obtain the desired bounds of $J(t)$ and $S(t)$:
\begin{equation}\label{J}
J=O(W_0^{(\eta)}+\theta^{\frac12}\eta [W_0^{(\eta)}]^2+\eta\frac\varepsilon{\kappa}\theta^{\frac16}[ W_0^{(\eta)}]^2)=O(W_0^{(\eta)}),
\end{equation}
\begin{equation}
 S=O(1+tW_0^{(\eta)})=O(1).
\end{equation}
And \eqref{J} improves bootstrap assumption in \eqref{assum}.

\section{Shock formation}
In this section, we show that in $\mathcal{R}_1$ the inverse foliation density $\rho_1$ goes to zero as time goes to a certain $T_\eta^*$.

With the following transport equation for $\rho_1$
\begin{equation*}
  \frac{\partial\rho_1}{\partial s_1}=c_{11}^1(\Phi)v^1+O\Big(\sum_{ k\neq 1}w^k\Big)\rho_1,
\end{equation*}
and the fact $c_{11}^1(\Phi)<0$, we have
\begin{equation}\label{381}
  -|c_{11}^1||v^1|-\Big|O\Big(\sum_{ k\neq 1}w^k\Big)\Big|\rho_1\leq\frac{\partial\rho_1}{\partial s_1}\leq-|c_{11}^1||v^1|+\Big|O\Big(\sum_{ k\neq 1}w^k\Big)\Big|\rho_1.
\end{equation}
By \eqref{68}, $|\Phi|=O(\eta W_0^{(\eta)}+\eta^2\frac\varepsilon\kappa \theta^{-\frac13}W_0^{(\eta)})\leq\delta$, choosing $\theta$ sufficiently small, it holds
\begin{equation}\label{y37}
 (1-\varepsilon)|c_{11}^1(0)|\leq |c_{11}^1(\Phi)|\leq (1+\varepsilon)|c_{11}^1(0)|.
\end{equation}
With bi-characteristic coordinates, we have
\begin{equation*}
  \int_0^{t}\sum_{k\neq 1}w^k(X_1(z_1,t'),t')dt'=O(\eta W_0^{(\eta)}+\eta J)=O(\eta W_0^{(\eta)}).
\end{equation*}
For $0<\eta\ll 1$, this implies
\begin{equation}\label{385}
  1-\varepsilon\leq \exp{\Big(\int_0^{t}O\big(\sum_{k\neq 1}w^k(X_1(z_1,t'),t')\big)dt'\Big)}\leq1+\varepsilon,
\end{equation}
and
\begin{equation}\label{386}
  1-\varepsilon\leq \exp{\Big(-\int_0^{t}O\big(\sum_{k\neq 1}w^k(X_1(z_1,t'),t')\big)dt'\Big)}\leq1+\varepsilon.
\end{equation}

Employing Gr\"{o}nwall inequality to \eqref{381} and combining with \eqref{y37}, \eqref{385}-\eqref{386},
we get
\begin{equation}\label{y38}
\begin{split}
  (1-\varepsilon)\Big(1-(1+\varepsilon)^2|c_{11}^1(0)|&\int_0^t|v^1(z_1,t')|dt'\Big)\leq\rho_1(z_1,t)\\
  &\leq (1+\varepsilon)\Big(1-(1-\varepsilon)^2|c_{11}^1(0)|\int_0^t|v^1(z_1,t')|dt'\Big)\\
  &\leq 1+\varepsilon.
\end{split}
\end{equation}
For $v^1$, we integrate
\begin{equation}\label{y39}
\frac{\partial v^1}{\partial s_1}=O\Big(\sum_{m\neq 1}w^m\Big)v^1+O\Big(\sum_{m\neq1,k\neq1\atop\{m,k\}\neq\{2,3\},\{m,k\}\neq\{4,5\}}w^mw^k\Big)\rho_1
\end{equation}
along $\mathcal{C}_1$ and obtain
\begin{equation}\label{79}
  v^1(z_1,t)\leq w^1_{(\eta)}(z_1,0)+O(tVJ+tV^2S)=w^1_{(\eta)}(z_1,0)+O(\eta [W_0^{(\eta)}]^2+\eta\frac\varepsilon\kappa\theta^{-\frac13}[W_0^{(\eta)}]^2).
\end{equation}
For sufficiently small $\theta$ and $\eta$, the above inequality yields
\begin{equation*}
  v^1(z_0,t)\leq (1+\varepsilon) W_0^{(\eta)}.
\end{equation*}
By using the first inequality of \eqref{y38}, we arrive at
\begin{equation}\label{711}
  \rho_1(z_0,t)\geq(1-\varepsilon)\Big(1-(1+\varepsilon)^3|c_{11}^1(0)|tW_0^{(\eta)}\Big),
\end{equation}
which shows that $\rho_1(z_0,t)>0$ when
\begin{equation*}
  t<\frac1{(1+\varepsilon)^3|c_{11}^1(0)|W_0^{(\eta)}}.
\end{equation*}
Meanwhile,
applying $\rho_1\leq 1+\varepsilon$ to \eqref{y39}, we have
\begin{equation*}
\frac{\partial v^1}{\partial s_1}\geq-\Big|O\Big(\sum_{m\neq 1}w^m\Big)\Big|v^1-(1+\varepsilon)\Big|O\Big(\sum_{m\neq1,k\neq1\atop\{m,k\}\neq\{2,3\},\{m,k\}\neq\{4,5\}}w^mw^k\Big)\Big|,
\end{equation*}
provided $v^1>0$.
Since
\begin{equation*}
  O\Big(\sum_{m\neq 1}\int_0^tw^mdt'\Big)=O(\eta W_0^{(\eta)}+\eta J)=O(\eta W_0^{(\eta)})
\end{equation*}
and
\begin{equation*}
  O\Big(\sum_{m\neq1,k\neq1\atop\{m,k\}\neq\{2,3\},\{m,k\}\neq\{4,5\}}\int_0^tw^mw^k\rho_1dt'\Big)=O(St V^2)=O(\eta^2 [W_0^{(\eta)}]^3+\eta^2\frac{\varepsilon^2}{\kappa^2}\theta^{-\frac23}[W_0^{(\eta)}]^3),
\end{equation*}
via Gr\"{o}nwall inequality, we have
\begin{equation} \label{vlow1}
 v^1(z_1,t)\geq(1-\varepsilon)\Big[w^1_{(\eta)}(z_1,0)-(1+\varepsilon)^2\big(\eta^2 [W_0^{(\eta)}]^3+\eta^2\frac{\varepsilon^2}{\kappa^2}\theta^{-\frac23}[W_0^{(\eta)}]^3\big)\Big].
\end{equation}
Choosing $\theta$ sufficiently small, we have
\begin{equation} \label{vlow2}
  (1+\varepsilon)^2O\big(\eta^2 [W_0^{(\eta)}]^3+\eta^2\frac{\varepsilon^2}{\kappa^2}\theta^{-\frac23}[W_0^{(\eta)}]^3\big)\leq \varepsilon W_0^{(\eta)},
\end{equation}
and setting $z_1=z_0$, it holds
\begin{equation*}
 v^1(z_0,t)\geq(1-\varepsilon)[W_0^{(\eta)}-\varepsilon W_0^{(\eta)}]=(1-\varepsilon)^2W_0^{(\eta)}.
\end{equation*}
By the second inequality of \eqref{y38}, we get
\begin{equation*}
  \rho_1(z_0,t)\leq (1+\varepsilon)\Big(1-(1-\varepsilon)^4|c_{11}^1(0)|tW_0^{(\eta)}\Big).
\end{equation*}
Together with \eqref{711}, we conclude that there exists $T_\eta^*$ (shock formation time) such that
\begin{equation*}
  \lim_{t\rightarrow T_\eta^*} \rho_1(z_0,t)=0.
\end{equation*}
And $T_\eta^*$ obeys
\begin{equation} \label{Tshock}
 \frac1{(1+\varepsilon)^3|c_{11}^1(0)|W_0^{(\eta)}}\leq T_\eta^*\leq\frac{1}{(1-\varepsilon)^4|c_{11}^1(0)|W_0^{(\eta)}},
\end{equation}
which is consistent with the requirement \eqref{368}
\begin{equation*}
  T_\eta^*\leq\frac{C}{W_0^{(\eta)}}.
\end{equation*}

\section{Lower bound for $\underline{S}$}\label{lows}
We next improve bootstrap assumption \eqref{sbar} via obtaining a uniformly positive lower bound for  $\{\rho_i\}_{i=2,3,4,5,6}$.

We start from $\rho_6$. Similarly to \eqref{y38}, we have
\begin{equation*}
\begin{split}
 \rho_6(z_6,t)\geq(1-\varepsilon)\Big(1-(1+\varepsilon)^2|c_{66}^6(0)|\int_0^t|v^6(z_6,t')|dt'\Big).
\end{split}
\end{equation*}
Taking $i=6$ in \eqref{y13} and integrating it along $\mathcal{R}_6$, proceeding as in \eqref{79} we get
\begin{equation*}
  v^6(z,t)\leq w^6_{(\eta)}(z,0)+O(tVJ+tV^2S)=w^6_{(\eta)}(z,0)+O(\eta [W_0^{(\eta)}]^2+\eta\frac{\varepsilon}{\kappa}\theta^{-\frac13}[W_0^{(\eta)}]^2).
\end{equation*}
With initial data condition \eqref{273} and the fact $0<\theta,\eta\ll 1$, the above inequality implies
\begin{equation*}
  v^6(z,t)\leq\frac{(1-\varepsilon)^4}{2(1+\varepsilon)^2}W_0^{(\eta)}.
\end{equation*}
 Noting that $c_{66}^6(0)=-c_{11}^1(0)$, we obtain a lower bound for $\rho_6$:
\begin{equation*}
\begin{split}
 \rho_6(z_6,t)&\geq(1-\varepsilon)\Big(1-(1+\varepsilon)^2|c_{11}^1(0)|\frac{(1-\varepsilon)^4}{2(1+\varepsilon)^2}W_0^{(\eta)}t \Big)\\
 &>(1-\varepsilon)\Big(1-|c_{11}^1(0)|\frac{(1-\varepsilon)^4 W_0^{(\eta)}}{2}\frac{1}{(1-\varepsilon)^4|c_{11}^1(0)|W_0^{(\eta)}} \Big)\\
 &=\frac{(1-\varepsilon)}{2}>0,\quad\text{for any}\ t<T_\eta^*.
\end{split}
\end{equation*}

For $\rho_3$ and $\rho_4$, we use the structure of equations in \eqref{infmrho}, i.e., $\rho_3 w^2$ and $\rho_4 w^5$ vanish in $\partial_{s_3}\rho_3$ and $\partial_{s_4}\rho_4$, respectively. Using similar arguments as above, with initial data condition \eqref{initial} we have
\begin{equation*}
\begin{split}
 \min_{i=3,4}\rho_i(z_i,t)&\geq(1-\varepsilon)\Big(1-(1+\varepsilon)^2\varepsilon\frac{(1-\varepsilon)^4|c_{11}^1(0)|}{(1+\varepsilon)^2}W_0^{(\eta)}t \Big)\\
 &>(1-\varepsilon)\Big(1-\varepsilon(1-\varepsilon)^4 |c_{11}^1(0)|W_0^{(\eta)}\frac{1}{(1-\varepsilon)^4|c_{11}^1(0)|W_0^{(\eta)}} \Big)\\
 &=(1-\varepsilon)^2>0,\quad \text{for}\ t\leq T_\eta^*.
\end{split}
\end{equation*}

For $\rho_2$ and $\rho_5$, since  $c_{22}^2=c_{55}^5=0$, a direct calculation yields
\begin{equation*}
  \rho_2(z,t)\geq 1-\varepsilon>0,\quad \text{and}\quad \rho_5(z,t)\geq 1-\varepsilon>0.
\end{equation*}

In conclusion, we obtain the following lower bound for $\underline{S}$
\begin{equation}\label{lbound}
    \underline{S}(s):=\min_{i=2,\cdots,6}\inf_{(z'_i,s'_i)\atop z'_i\in[\eta,2\eta],\ 0\leq s'_i\leq t} \rho_i(z'_i,s'_i)\geq \frac{(1-\varepsilon)^2}{2}.
\end{equation}
With $\varepsilon\in(0,\frac1{100}]$ and $\kappa\in(0,\frac1{100}]$, the inequality \eqref{lbound} improves \eqref{sbar}.\\

We then move to prove that $\{\rho_i\}_{i=1,\cdots,6}$ obey a positive lower bound at points $(x,t)$ outside $\mathcal{R}_i$. For $t<T_{\eta}^*$ we have
$$\min_{i=1,\cdots,6}\inf_{(z'_i,s'_i)\atop  z_i'\not\in [\eta, 2\eta],\ 0\leq s'_i\leq t} \rho_i(z'_i,s'_i)\geq \frac{(1-\varepsilon)^2}{2} \mbox{ for some } \varepsilon\in(0,\frac1{100}].$$

\noindent Proof: For $(x,t)$ outside of $\mathcal{R}_i$, lying at the intersection of $\mathcal{C}_i$ and $\{\mathcal{C}_k\}_{k\neq i}$, we estimate $\rho_i$ along $\mathcal{C}_i$.
\begin{center}
\begin{tikzpicture}
\draw[->](0,0)--(9,0)node[right,below]{$t=0$};
\draw[dashed](0,1)--(9,1)node[right,below]{$t=t_0^{(\eta)}$};
\node [below]at(3.6,0){$2\eta$};
\node [below]at(2.4,0){$\eta$};
\filldraw [black] (3.5,0) circle [radius=0.01pt]
(5,0.8) circle [radius=0.01pt]
(6,1) circle [radius=0.01pt]
(8,1.8) circle [radius=0.01pt];
\draw (3.5,0)..controls (5,0.8)and(6,1)..(8,1.8);

\filldraw [black] (2.5,0) circle [radius=0.01pt]
(4,1) circle [radius=0.01pt]
(5.5,1.5) circle [radius=0.01pt]
(6,1.68) circle [radius=0.01pt];
\draw (2.5,0)..controls (4,1) and (5.5,1.5)..(6,1.68);
\node [above] at(7.5,1.6){$\mathcal{R}_i$};

\filldraw [black] (2,0) circle [radius=0.01pt]
(2.8,0.5) circle [radius=0.01pt]
(3.5,1) circle [radius=0.01pt]
(4.5,1.5) circle [radius=0.8pt];
\draw [color=blue](2,0)..controls (2.8,0.5) and (3.5,1)..(4.5,1.5);
\node [above]at(4.5,1.5){$(x,t)$};
\node [above,color=blue]at(4,1){$\mathcal{C}_i$};
\node[below] at (2,0){$z_i$};

\filldraw [black] (6,0) circle [radius=0.01pt]
(5.5,0.5) circle [radius=0.01pt]
(5,1) circle [radius=0.01pt]
(4.5,1.5) circle [radius=0.8pt];
\draw [color=green](6,0)..controls (5.5,0.5) and (5,1)..(4.5,1.5);
\node [above]at(4.5,1.5){$(x,t)$};
\node [above,color=green]at(5,1){$\mathcal{C}_k$};
\draw[dashed,color=red](0,1.6)--(9,1.6)node[right,above]{$t=T_\eta^*$};
\end{tikzpicture}
\end{center}
From \eqref{y14}, we have
\begin{equation}\label{91}
 \frac{\partial\rho_i}{\partial s_i}\geq-\Big|O\Big(\sum_k w^k\Big)\Big|\rho_i.
\end{equation}
For $z_i\notin [\eta,2\eta]$, $t\in[t_0^{(\eta)},T_\eta^*)$, utilizing \eqref{y15} and \eqref{560} we obtain
\begin{equation}
\begin{split}
&\int_0^t\Big|O\Big(\sum_k w^k(z_i,s)\Big)\Big|ds\\
=&\int_0^{t_0^{(\eta)}}\Big|O\Big(\sum_k w^k(z_i,s)\Big)\Big|ds+\int_{t_0^{(\eta)}}^t\Big|O\Big(\sum_{k} w^k(z_i,s)\Big)\Big|ds\\
=&O\Big(\eta W_0^{(\eta)}+tV+\eta J\Big).
\end{split}
\end{equation}
By the bound \eqref{tv} for $tV$ and bound \eqref{64} for $J$, we have
\begin{equation}
\begin{split}
\int_0^t\Big|O\Big(\sum_k w^k(z_i,s)\Big)\Big|ds=O\Big(\eta W_0^{(\eta)}+\eta\varepsilon\theta^{-\frac13} W_0^{(\eta)}\Big).
\end{split}
\end{equation}
Employing Gr\"{o}nwall inequality for \eqref{91}, by choosing $\theta$ small, we conclude
\begin{equation}
  \rho_i(z_i,t)\geq\exp\Big(-\int_0^t\Big|O\Big(\sum_k w^k(z_i,s)\Big)\Big|ds\Big)\geq 1-\varepsilon,
\end{equation}
for $z_i\notin [\eta,2\eta],\ t\in[t_0^{(\eta)},T_\eta^*)$. This estimate illustrates that there is no shock formed outside $\mathcal{R}_1, \cdots, \mathcal{R}_6$. And by the lower bound of $\underline{S}(t)$, we also conclude that no shock emerges in  $\mathcal{R}_2, \cdots, \mathcal{R}_6$. For the whole spacetime region before $t=T^*_{\eta}$, the \underline{only} shock happens in $\mathcal{R}_1$, i.e., $\rho_1(t)\rightarrow 0$ as $t\rightarrow T^*_{\eta}$.

For $i\in\{1,\cdots,6\}$, recall that $w^i$ is uniformly bounded outside $\mathcal{R}_i$. Utilizing the fact that $v^i=\rho_i w^i$ is uniformly bounded inside $\mathcal{R}_i$, together with the lower bound of $\{\rho_i\}_{i\in\{2,\cdots,6\}}$, we conclude that $\{w^i\}_{i\in\{2,\cdots,6\}}$ are also uniformly bounded inside $\{\mathcal{R}_i\}_{i\in\{2,\cdots,6\}}$. For the whole spacetime region before $t=T^*_{\eta}$, the \underline{only} singularity happens in $\mathcal{R}_1$, i.e., $w^1(t)\rightarrow +\infty$ as $t\rightarrow T^*_{\eta}$. \\

We summarize the above conclusions into
\begin{prop}\label{wave dynamics}
 For $t\leq T_\eta^*$, we have
 \begin{itemize}
   \item 
   $\{\rho_i\}_{i=2,\cdots,6}$ are bounded away from zero in the whole $(x,t)$-plane. $\rho_1$ obeys  a positive lower bound outside $\mathcal{R}_1$ and within $\mathcal{R}_1$ it holds
   $\rho_1\rightarrow 0$ as $t\rightarrow T^*_{\eta}$.  And the first singularity (shock) forms.

   \item $\{w^i\}_{i=2,\cdots,6}$ are uniformly bounded in the whole $(x,t)$-plane. $w^1$ obeys a  uniform bound outside $\mathcal{R}_1$. Within $\mathcal{R}_1$ it holds
   $w^1\rightarrow +\infty$ as $t\rightarrow T^*_{\eta}$ and $w^1$ being finite for $t<T^*_{\eta}$.

\end{itemize}
\end{prop}

\noindent These further imply that the solutions to the Cauchy problem of system \eqref{y1} are smooth before time $T_\eta^*$.

\section{Estimate for $\partial_{z_1}\rho_1$}\label{supbound}
In this section, we fix $\eta$ and employ characteristic coordinates and bi-characteristic coordinates introduced in \eqref{cc}, \eqref{37} and \eqref{38} for $i\neq j$:
\begin{equation*}
  (x,t)=\big(X_i(z_i,s_i),s_i\big)=\big(X_i(y_i,t'(y_i,y_j)),t'(y_i,y_j)\big)=\big(X_j(y_j,t'(y_i,y_j)),t'(y_i,y_j)\big).
\end{equation*}
For any smooth function $f(x,t)=f\big(X_i(y_i,t'(y_i,y_j)),t'(y_i,y_j)\big)$,  with \eqref{35}, \eqref{38}, we calculate
\begin{equation}\label{9.2}
\begin{split}
  \partial_{y_i}f=&(\partial_{z_i}X_i+\partial_{s_i}X_i\partial_{y_i}t')\partial_{x}f+\partial_{y_i}t'\partial_{t}f\\
  =&\Big(\rho_i+\lambda_i\frac{\rho_i}{\lambda_j-\lambda_i}\Big)\partial_{x}f+\frac{\rho_i}{\lambda_j-\lambda_i}\partial_{t}f\\
  =&\rho_i\partial_{x}f+\frac{\rho_i}{\lambda_j-\lambda_i}(\partial_{t}+\lambda_i\partial_{x})f\\
  =&\partial_{z_i}f+\frac{\rho_i}{\lambda_j-\lambda_i}\partial_{s_i}f.
  \end{split}
\end{equation}
Bi-characteristic coordinates also means $f(x,t)=f\big(X_j(y_j,t'(y_i,y_j)),t'(y_i,y_j)\big)$, thus we deduce
\begin{equation}\label{9.3}
\begin{split}
  \partial_{y_i}f=&\partial_{s_j}X_j\partial_{y_i}t'\partial_{x}f+\partial_{y_i}t'\partial_{t}f\\
  =&\lambda_j\frac{\rho_i}{\lambda_j-\lambda_i}\partial_{x}f+\frac{\rho_i}{\lambda_j-\lambda_i}\partial_{t}f\\
  =&\frac{\rho_i}{\lambda_j-\lambda_i}(\partial_{t}+\lambda_j\partial_{x})f\\
  =&\frac{\rho_i}{\lambda_j-\lambda_i}\partial_{s_j}f.
  \end{split}
\end{equation}
From \eqref{9.2} and \eqref{9.3}, we hence have that
transformations between these coordinates satisfy
\begin{equation}\label{92}
  \partial_{y_i}=\frac{\rho_i}{\lambda_j-\lambda_i}\partial_{s_j}=\partial_{z_i}+\frac{\rho_i}{\lambda_j-\lambda_i}\partial_{s_i}.
\end{equation}
We fix $y_1$ along $\mathcal{C}_1$ and choose $y_6$ acting as a parameter. From \eqref{92}, we have
\begin{equation*}
  \partial_{y_1}=\partial_{z_1}+\frac{\rho_1}{\lambda_6-\lambda_1}\partial_{s_1}=\partial_{z_1}-\frac{\rho_1}{2\lambda_1}\partial_{s_1},
\end{equation*}
where we use $\lambda_6=-\lambda_1$. Hence,
\begin{equation}\label{rhod}
 \partial_{z_1}\rho_1=\partial_{y_1}\rho_1+\frac{\rho_1}{2\lambda_1}\partial_{s_1}\rho_1.
\end{equation}

To bound $\partial_{z_1}\rho_1$, we start with controlling $\partial_{y_1}\rho_1$ by studying its evolution equation. Let
\begin{equation*}
\tau_{1} ^{(6)}:=\partial_{y_1}\rho_1,\quad\pi_{1}^{(6)}:=\partial_{y_1}v^1.
\end{equation*}
With the help of $[\partial_1,\partial_6]=0$, we have
\begin{equation}\label{tau1}
  \partial_{y_6}\tau_{1}^{(6)}= \partial_{y_6}\partial_{y_1}\rho_1=\partial_{y_1}\partial_{y_6}\rho_1.
\end{equation}
Since
\begin{equation}\label{96}
  \partial_{y_6}=\frac{\rho_6}{\lambda_1-\lambda_6}\partial_{s_1}=\frac{\rho_6}{2\lambda_1}\partial_{s_1},
\end{equation}
with \eqref{tau1} and \eqref{96}, by calculation we obtain
\begin{equation}\label{97}
\begin{split}
  \partial_{y_6}\tau_{1}^{(6)}=&\partial_{y_1}\Big(\frac{\rho_6}{2\lambda_1}\partial_{s_1}\rho_1\Big)\overset{\eqref{y14}}{=}\partial_{y_1}\Big(\frac{\rho_1\rho_6}{2\lambda_1}\sum_m c_{1m}^1w^m\Big)\\
  =&\frac{\rho_6}{2\lambda_1}\Big(c_{11}^1\partial_{y_1}v^1+\sum_{m\neq1}c_{1m}^1w^m \partial_{y_1}\rho_1\Big)\\
  &+\frac{\rho_1}{2\lambda_1}\Big(\sum_{m\neq6}c_{1m}^1w^m\partial_{y_1}\rho_6+c_{16}^1\partial_{y_1}v^6\Big)\\
  &+\frac{\rho_1\rho_6}{2\lambda_1}\Big(\sum_{m\neq1,m\neq6}c_{1m}^1\frac1{\rho_m}\partial_{y_1}v^m-\sum_{m\neq1,m\neq6}c_{1m}^1\frac{w^m}{\rho_m}\partial_{y_1}\rho_m\Big)\\
  &-\frac{\partial_{y_1}\lambda_1}{2\lambda_1^2}\rho_1\rho_6\sum_m c_{1m}^1w^m+\frac{\rho_1\rho_6}{2\lambda_1}\sum_m \partial_{y_1}c_{1m}^1w^m.
\end{split}
\end{equation}
With bi-characteristic coordinates $(y_1,y_m)$ ($m\neq1$), we have
\begin{equation}\label{9.11}
  \begin{split}
    \partial_{y_1}\lambda_1=&\nabla_\Phi\lambda_1\cdot\partial_{y_1}\Phi=\nabla_\Phi\lambda_1\cdot[\partial_{s_m}X_m\partial_{y_1}t'\partial_{x}\Phi+\partial_{y_1}t'\partial_{t}\Phi]\\
    =&\nabla_\Phi\lambda_1\cdot[\lambda_m\frac{\rho_1}{\lambda_m-\lambda_1}\sum_kw^kr_k+\frac{\rho_1}{\lambda_m-\lambda_1}(-A(\Phi)\sum_kw^kr_k)]\\
    =&O(v^1+\rho_1\sum_{k\neq1}w^k)
  \end{split}
\end{equation}
and
\begin{equation}\label{9.12}
  \begin{split}
    \partial_{y_1}c_{1m}^1=&\nabla_\Phi c_{1m}^1\cdot\partial_{y_1}\Phi=\nabla_\Phi c_{1m}^1\cdot[\partial_{s_m}X_m\partial_{y_1}t'\partial_{x}\Phi+\partial_{y_1}t'\partial_{t}\Phi]\\
    =&\nabla_\Phi c_{1m}^1\cdot[\lambda_m\frac{\rho_1}{\lambda_m-\lambda_1}\sum_kw^kr_k+\frac{\rho_1}{\lambda_m-\lambda_1}(-A(\Phi)\sum_kw^kr_k)]\\
    =&O(v^1+\rho_1\sum_{k\neq1}w^k).
  \end{split}
\end{equation}
With \eqref{92}, it also holds
\begin{equation*}
  \partial_{y_1}=\frac{\rho_1}{\lambda_m-\lambda_1}\partial_{s_m},\quad\text{for}\quad m\neq1.
\end{equation*}
By \eqref{y14} and \eqref{y13}, we obtain
\begin{equation}\label{9.14}
  \partial_{y_1}\rho_m=\frac{\rho_1}{\lambda_m-\lambda_1}\partial_{s_m}\rho_m=O\Big(\rho_mv^1+\rho_1\rho_m\sum_{k\neq 1}w^k\Big),
\end{equation}
and
\begin{equation}\label{9.15}
  \partial_{y_1}v^m=\frac{\rho_1}{\lambda_m-\lambda_1}\partial_{s_m}v^m=O\Big(\rho_mv^1\sum_{k\neq1}w^k+\rho_1\rho_m\sum_{j\neq1,k\neq1\atop j\neq k}w^jw^k\Big).
\end{equation}
Together with estimates in Section \ref{3d}, \eqref{9.11}-\eqref{9.12}, \eqref{9.14}-\eqref{9.15} imply that the right hand side of \eqref{97} are bounded as below
\begin{equation}\label{linear}
\begin{split}
  \partial_{s_1}\tau_{1}^{(6)}=B_{11}^\eta\tau_1^{(6)}+B_{12}^\eta\pi_1^{(6)}+B_{13}^\eta,
\end{split}
\end{equation}
where $B_{11}^\eta,B_{12}^\eta,B_{13}^\eta$ are uniformly bounded constants depending on $\eta$.

Similarly, we have the evolution equation for $\pi_1^{(6)}$ with bounded coefficients:
\begin{equation}\label{917}
  \begin{split}
  \partial_{y_6}\pi_{1}^{(6)}=&\frac{\rho_6 }{2\lambda_1}\Big(\sum_{p\neq 1,q\neq 1\atop p\neq q}\gamma_{pq}^1w^pw^q \tau_{1}^{(6)}+\sum_{p\neq 1}\gamma_{1p}^1w^p\pi_{1}^{(6)}\Big)\\
  &-\frac{\partial_{y_1}(\lambda_1)}{4\lambda_1^2}\Big(\sum_{p\neq 1}\gamma_{1p}^1 w^p v^1\rho_6+\sum_{p\neq 1,q\neq 1\atop p\neq q}\gamma_{pq}^1w^pw^q\rho_6\rho_1\Big)\\
    &+\frac{\rho_6\rho_1}{2\lambda_1}\Big(\sum_{p\neq 1}\partial_{y_1}\gamma_{1p}^1 w^p w^m+\sum_{p\neq 1,q\neq 1\atop p\neq q}\partial_{y_1}\gamma_{pq}^1w^pw^q\Big)\\
    &+\frac{\rho_1}{2\lambda_1 }\Big(\sum_{p\neq 1}\gamma_{1p}^1 w^p w^1+\sum_{p\neq 1,q\neq 1\atop p\neq q}\gamma_{pq}^1w^pw^q\Big)\partial_{y_1}\rho_6\\
    &+\frac{\rho_6\rho_1}{2\lambda_1 }\Big(\sum_{p\neq 1,p\neq6}\gamma_{1p}^1\frac{w^1}{\rho_p}+\sum_{p\neq 1,q\neq 1\atop p\neq q}\gamma_{pq}^1\frac{w^q}{\rho_p}\Big)(\partial_{y_1}v^p-w^p\partial_{y_1}\rho_p)\\
    &++\frac{\rho_6\rho_1}{2\lambda_1 }\sum_{p\neq 1,q\neq 1\atop p\neq q}\gamma_{pq}^1\frac{w^p}{\rho_q}(\partial_{y_1}v^q-w^q\partial_{y_1}\rho_q)
    +\frac{\rho_1}{2\lambda_1 }\gamma_{16}^1 w^1\partial_{y_1}v^6\\
    =&B_{21}^\eta\tau_1^{(6)}+B_{22}^\eta\pi_1^{(6)}+B_{23}^\eta,
  \end{split}
\end{equation}
where $B_{21}^\eta,B_{22}^\eta,B_{23}^\eta$ are uniformly bounded constants depending on $\eta$.

We next check that the initial data of $\tau_1^{(6)}$ and $\pi_m^{(6)}$ are bounded.
Since
\begin{equation*}
  \rho_1(z_1,0)=1, \quad v^1(z_1,0)=w^1(z_1,0).
\end{equation*}
using \eqref{92}, for fixed $\eta$ we have
\begin{equation*}
  \begin{split}
    \tau_1^{(6)}(z_1,0):=&\partial_{y_1}\rho_1(z_1,0)\\
    =&\partial_{z_1}\rho_1(z_1,0)-\frac{\rho_1(z_1,0)}{2\lambda_1}\partial_{s_1}\rho_1(z_1,0)\\
    =&-\frac{1}{2\lambda_1}\sum_{k}c_{1k}^1w^k(z_1,0)\rho_1(z_1,0)\\
    =&-\frac{1}{2\lambda_1}\sum_{k}c_{1k}^1w^k(z_1,0)\\
    =&O(W_0^{(\eta)})<+\infty,
  \end{split}
\end{equation*}
and
\begin{equation*}
  \begin{split}
    \pi_1^{(6)}(z_1,0):=&\partial_{y_1}v^1(z_1,0)\\
    =&\partial_{z_1}v^1(z_1,0)-\frac{\rho_1(z_1,0)}{2\lambda_1}\partial_{s_1}v^1(z_1,0)\\
    =&\partial_{z_1}w^1(z_1,0)-\frac{1}{2\lambda_1}\sum_{q\neq 1,q\neq p}\gamma_{pq}^1w^p(z_1,0)w^q(z_1,0)\rho_1(z_1,0)\\
    =&O(\partial_{z_1}w^1(z_1,0)+[W_0^{(\eta)}]^2)<+\infty.
  \end{split}
\end{equation*}
Applying Gr\"owall inequality to \eqref{linear} and \eqref{917}, for $t\leq T_\eta^*$ we have that $\tau_{1}^{(6)}:=\partial_{y_1}\rho_1$ is bounded by a constant depending on $\eta$.

Back to \eqref{rhod} and \eqref{y14}, we have
\begin{equation*}
  \partial_{z_1}\rho_1=\partial_{y_1}\rho_1+O(v^1+\sum_{m\neq 1}w^m\rho_1).
\end{equation*}
With the bounds for $J(t)$, $S(t)$ and $V(t)$ in Section \ref{3d}, consequently, we conclude $\partial_{z_1}\rho_1$ is bounded by $C_\eta$, a constant depending on $\eta$.

\section{Ill-posedness mechanism}\label{ill}
From our construction of modified Lindblad-type initial data \eqref{W0}-\eqref{data} and the bounds \eqref{Tshock} obtained for the shock formation time $T^*_\eta$, we immediately conclude that $T^*_\eta \to 0$ as $\eta \to 0$. This gives the ill-posedness in Theorem \ref{3D}.

In the following part of this section, we further show that the $H^2$ norm of the solutions to elastic waves \eqref{y1} is infinite at shock formation time $T^*_\eta$ when restricted to a spatial region $\Omega_{T_\eta^*}$. In the picture below, $\Omega_{T_\eta^*}:=\{(x,T_\eta^*):x=X_1(z,T_\eta^*)\ \text{and}\ \eta\leq z\leq2\eta\}$.
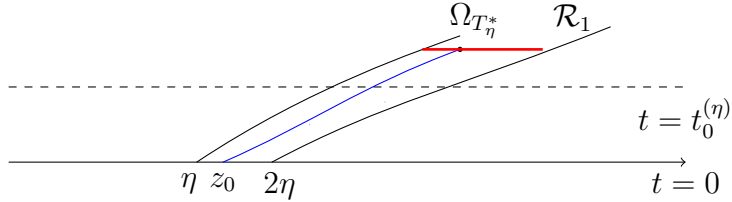
\begin{center}
\begin{tikzpicture}
\draw[->](0,0)--(9,0)node[right,below]{$t=0$};
\draw[dashed](0,1)--(9,1)node[right,below]{$t=t_0^{(\eta)}$};
\node [below]at(3.6,0){$2\eta$};
\node [below]at(2.4,0){$\eta$};
\filldraw [black] (3.5,0) circle [radius=0.01pt]
(5,0.8) circle [radius=0.01pt]
(6,1) circle [radius=0.01pt]
(8,1.8) circle [radius=0.01pt];
\draw (3.5,0)..controls (5,0.8)and(6,1)..(8,1.8);

\filldraw [black] (2.5,0) circle [radius=0.01pt]
(4,1) circle [radius=0.01pt]
(5.5,1.5) circle [radius=0.01pt]
(6,1.68) circle [radius=0.01pt];
\draw (2.5,0)..controls (4,1) and (5.5,1.5)..(6,1.68);
\node [above] at(7.5,1.6){$\mathcal{R}_1$};

\filldraw [black] (2.85,0) circle [radius=0.01pt]
(4,0.5) circle [radius=0.01pt]
(4.65,1) circle [radius=0.01pt]
(6,1.5) circle [radius=0.8pt];
\draw [color=blue](2.85,0)..controls (4,0.5) and (4.65,1)..(6,1.5);
\node [below]at(2.85,0){$z_0$};

\node [above]at(6.2,1.5){$\Omega_{T_\eta^*}$};

\draw[line width=1pt, color=red](5.5,1.5)--(7.1,1.5);
\end{tikzpicture}
\end{center}

Considering acoustic metric $g_{c_1}$ according to the wave equation for $U^1$ in \eqref{y1}, we have
\begin{align*}
  \begin{split}
    {(g_{c_1}^{-1})}^{\alpha\beta}=\left(
              \begin{array}{cccccc}
                1 & 0 & 0 &0 \\
                0 & -(c_1^2+2\sigma_0 \partial_x u^1) & \bar{\sigma}\partial_x u^2 &\bar{\sigma}\partial_x u^3  \\
                0 & \bar{\sigma}\partial_x u^2 & -(c_2^2+2\sigma_1 \partial_x u^1) &0   \\
                0 & \bar{\sigma}\partial_x u^3 &0 &-(c_2^2+2\sigma_1 \partial_x u^1) \\
              \end{array}
            \right),
  \end{split}
\end{align*}
where $\bar{\sigma}=2\sigma_1+\sigma_2+\sigma_3-\sigma_4$. For $\Phi\in B_\delta^6(0)$, $g_{c_1}$ is a small perturbation of $g_{c_1}^{(Flat)}$
$$
g_{c_1}^{(Flat)}=dt^2-{c_1}^{-2}d x^2-{c_2}^{-2}(dY^2)^2-{c_2}^{-2}(dY^3)^2.
$$
Using {\color{black}$g_{c_1}^{(Flat)}$} as the induced metric on $\Omega_{T_\eta^*}$, we have that:
In 3D $\Omega_{T_\eta^*}$ to be an ellipsoidal ball, centered at
$$P=\big(\frac{X_1(2\eta, T_\eta^*)+X_1(\eta, T_\eta^*)}{2},0,0\big)=(c_1T_\eta^*+\frac32\eta,0,0)$$
with major axis $X_1(2\eta, T_\eta^*)-X_1(\eta,T_\eta^*)=\eta$ and two minor axes $\frac{c_2}{c_1}[X_1(2\eta, T_\eta^*)-X_1(\eta,T_\eta^*)]=\frac{c_2\eta}{c_1}$.

To calculate $\int_{\Omega_{T_\eta^*}}|w^1(Y,T_\eta^*)|^2dxdY^2dY^3$, we first compute:
\begin{prop}\label{pro}
Let $\Omega_{T_\eta^*}^r=\{(Y^2,Y^3):(x,Y^2,Y^3)\in\Omega_{T_\eta^*}\}$. For $(x,Y^2,Y^3)\in\Omega_{T_\eta^*}$ with $x$ along $\mathcal{C}_1$ starting from $z$, we have
\begin{equation} \label{volume}
\int_{\Omega_{T_\eta^*}^r} dY^2 dY^3 \sim\big(z-\eta+O(\varepsilon)\eta\big)\big(2\eta-z+O(\varepsilon)\eta\big),
\end{equation}
and
\begin{equation}\label{V}
  |\Omega_{T_\eta^*}|=\int_{\Omega_{T_\eta^*}} dxdY^2 dY^3\sim \eta^3.
\end{equation}
\end{prop}
{\it Proof.}
In fact, $\Omega_{T_\eta^*}$ is given by
\begin{equation}\label{ball}
\frac{{c_2^2(x-c_1T_\eta^*-\frac{3\eta}{2})}^2}{c_1^2}+(Y^2)^2+(Y^3)^2 \leq \frac{c_2^2}{c_1^2} {\left(\frac{\eta}{2} \right)}^2.
\end{equation}
It follows,
\begin{equation}\label{102}
\begin{split}
  (Y^2)^2+(Y^3)^2\leq&\frac{c_2^2}{c_1^2} \Big[{\left(\frac{\eta}{2}\right)}^2-{\left(x-c_1T_\eta^*-\frac{3\eta}{2}\right)}^2\Big]
\end{split}\end{equation}
Since $\lambda_1(0)=c_1$, it holds that
\begin{equation*}
  \lambda_1(\Phi)=c_1(1+O(\varepsilon)),\quad \text{for}\quad \Phi\in B_\delta^6(0)
\end{equation*}
with sufficiently small $\delta$. Hence, along $\mathcal{C}_1$ starting from $z$, we have
\begin{equation}\label{104}
  x=z+(1+O(\varepsilon))c_1T_\eta^*.
\end{equation}
Employing \eqref{104} to \eqref{102}, we obtain
\begin{equation*}
\begin{split}
  \int_{\Omega_{T_\eta^*}^r} dY^2 dY^3=&\frac{c_2^2}{c_1^2} \Big[{\left(\frac{\eta}{2}\right)}^2-{\left(x-c_1T_\eta^*-\frac{3\eta}{2}\right)}^2\Big]\\
  \sim&\frac{c_2^2}{c_1^2} \Big[{\left(\frac{\eta}{2}\right)}^2-{\left(z-\frac{3\eta}{2}+O(\varepsilon)T_\eta^*\right)}^2\Big]\\
  \sim&\big(z-\eta+O(\varepsilon)\eta\big)\big(2\eta-z+O(\varepsilon)\eta\big).
\end{split}\end{equation*}
So we complete proof of \eqref{volume}. Estimate \eqref{V} is obvious by the expression of $\Omega_{T_\eta^*}$ in \eqref{ball}.\hfill$\Box$

With bi-characteristic coordinates and Proposition \eqref{pro}, we have
\begin{equation*}
  \int_{\Omega_{T_\eta^*}}|w^1|^2dxdY^2dY^3  \geq C\int_{\eta}^{2\eta}\Big|\frac{v^1}{\rho_1}(z,T_\eta^*)\Big|^2\rho_1(z,T_\eta^*) \big(z-\eta+O(\varepsilon)\eta\big)\big(2\eta-z+O(\varepsilon)\eta\big)dz.
\end{equation*}
By \eqref{vlow1}, \eqref{vlow2}, we have
 \begin{equation} \label{vlow}
 v^1(z,t)\geq(1-\varepsilon)[w^1_{(\eta)}(z,0)-\varepsilon W_0^{(\eta)}].
\end{equation}
Restrict the spatial integration region to a subinterval $(z_0,z_0^*]\subseteq[\eta,2\eta]$, where $w^1_{(\eta)}(z,0)>\frac{1}{2}W_0^{(\eta)}$ for $z\in(z_0,z_0^*]$. Then by \eqref{vlow} and fact that $\rho_1(z_0,T_\eta^*)=0$, we have
\begin{equation*}
\begin{split}
    &\|w^1(\cdot,T_\eta^*)\|_{L^2(\Omega_{T_\eta^*})}^2 \\
    &\geq C\int_{z_0}^{z_0^*}\frac{{(v^1)}^2}{\rho_1(z,T_\eta^*)} \big(z-\eta+O(\varepsilon)\eta\big)\big(2\eta-z+O(\varepsilon)\eta\big) dz\\
  &\geq C{(1-\varepsilon)}^2{(1-2\varepsilon)}^2[W_0^{(\eta)}]^2\int_{z_0}^{z_0^*}\frac{1}{\rho_1(z,T_\eta^*)} \big(z-\eta+O(\varepsilon)\eta\big)\big(2\eta-z+O(\varepsilon)\eta\big)dz\\
  &\geq C(z_0-\eta)(2\eta-z_0^*)[W_0^{(\eta)}]^2\int_{z_0}^{z_0^*}\frac{1}{\rho_1(z,T_\eta^*)-\rho_1(z_0,T_\eta^*)}dz\\
  &\geq C(z_0-\eta)(2\eta-z_0^*)[W_0^{(\eta)}]^2\int_{z_0}^{z_0^*}\frac{1}{(\sup_{z\in(z_0,z_0^*]}|\partial_z\rho_1|)(z-z_0)}dz.
\end{split}\end{equation*}
With the crucial boundedness in Section \ref{supbound} for $\partial_z\rho_1$, we conclude that
\begin{equation*}
 \|w^1(\cdot,T^*_\eta)\|_{L^2(\Omega_{T^*_\eta})}\geq C_\eta\int_{z_0}^{z_0^*}\frac1{z-z_0}dz=+\infty.
\end{equation*}

For the component $U^1$ of elastic waves, we derive from \eqref{uxx}, \eqref{y4} and \eqref{right}
\begin{equation*}
  \begin{split}
    \|\partial_x^2U^1\|_{L^2(\Omega_{T_\eta^*})}=&\|\sum_{k=1}^6 w^kr_{k1}\|_{L^2(\Omega_{T_\eta^*})}\geq\|w^1r_{11}\|_{L^2(\Omega_{T_\eta^*})}-\sum_{k=2}^6\|\frac{v^k}{\rho_k}r_{k1}\|_{L^2(\Omega_{T_\eta^*})}\\
    \geq&C \Big[\|w^1\|_{L^2(\Omega_{T_\eta^*})}-\|\frac{v^3}{\rho_3}\|_{L^2(\Omega_{T_\eta^*})}-\|\frac{v^4}{\rho_4}\|_{L^2(\Omega_{T_\eta^*})}-\|\frac{v^6}{\rho_6}\|_{L^2(\Omega_{T_\eta^*})}\Big]\\
    \geq&C\Big[\|w^1\|_{L^2(\Omega_{T_\eta^*})}-3\varepsilon\eta^\frac32\frac{W_0^{(\eta)}}{(1-\varepsilon)^2}\Big].
  \end{split}
\end{equation*}
Finally, since $\|w^1\|_{L^2(\Omega_{T_\eta^*})}=+\infty$, we obtain
\begin{equation*}
  \|\partial_x^2U^1\|_{L^2(\Omega_{T_\eta^*})}=+\infty.
\end{equation*}

In summary, with initial data requirements:
 \begin{equation*}
 \left\{\begin{aligned}
  &W_0^{(\eta)}=\max_{i=1,\cdots,6}\sup_z|w^i_{(\eta)}(z_i,0)|=w^1_{(\eta)}(z_0,0),\\
   &\max_{i=3,4}\sup_z|w^i_{(\eta)}(z_i,0)|\leq \min\Big\{\frac{(1-\varepsilon)^4|c_{11}^1(0)|}{(1+\varepsilon)^3}W_0^{(\eta)},W_0^{(\eta)}\Big\},\\
   &\sup_z|w^6_{(\eta)}(z_6,0)|\leq \frac{(1-\varepsilon)^4}{2(1+\varepsilon)^3}W_0^{(\eta)},
  \end{aligned}\right.
\end{equation*}
and  a general condition
$$c_{11}^1(\Phi)<0,\quad\forall \ \Phi\in B_{2\delta}^6(0),$$
we prove
\begin{thm}\label{detail}
The Cauchy problems of the 3D elastic wave equations \eqref{y1} are ill-posed in $H^3(\mathbb{R}^3)$ in the following sense: We construct a family of compactly supported  smooth initial data $(U_0^{(\eta)}, U_1^{(\eta)})$ with $${\|U_0^{(\eta)}\|}_{{\dot{H}}^3(\mathbb{R}^3)}+{\|U_1^{(\eta)}\|}_{{\dot{H}}^2(\mathbb{R}^3)} \to 0,\quad \text{as} \quad \eta\to0.$$ Let $T_\eta^*$ be the largest time such that \eqref{y1} $($with a general condition \eqref{sigma}$)$ has a solution $U_\eta\in C^\infty(\mathbb{R}^3\times[0,T_\eta^*) )$. As $\eta \to 0$, we have $T_\eta^* \to 0$.

Moreover, in a spatial region $\Omega_{T_\eta^*}$ the $H^2$ norm of the solution to elastic waves \eqref{y1}  blows up at shock formation time $T^*_{\eta}$:
\begin{equation*}
  \|U_\eta(\cdot,T_\eta^*)\|_{H^2(\Omega_{T_\eta^*})}=+\infty.
\end{equation*}
\end{thm}

\section{2D Case}\label{sec2d}
The aim of this section is to prove low regularity ill-posedness for Cauchy problems of two dimensional elastic waves.

For a 2D elastic waves model with plane symmetry, equations in \eqref{y1}
\begin{equation}\label{111}
    \partial_t^2 U-c_2^2\Delta U-(c_1^2-c_2^2)\nabla(\nabla\cdot U)=N(\nabla U,\nabla^2U)
\end{equation}
can be reduced to
\begin{align}\label{2d1}
    \left\{\begin{array}{lll}
    \partial_t^2u^1-c_1^2\partial_x^2u^1=\sigma_0\partial_x(\partial_xu^1)^2+\sigma_1\partial_x(\partial_xu^2)^2,\\
    \partial_t^2u^2-c_2^2\partial_x^2u^2=2\sigma_1\partial_x(\partial_xu^1\partial_xu^2).
    \end{array}\right.
\end{align}
Let
\begin{align*}
  \phi_1:=\partial_xu^1,\quad \phi_2:=\partial_xu^2,\quad \phi_3:=\partial_tu^1,\quad \phi_4:=\partial_tu^2.
\end{align*}
And set $\Phi:=(\phi_1,\phi_2,\phi_3,\phi_4)^T$. Then system \eqref{2d1} is equivalent to
\begin{align}\label{2d2}
  \partial_t\Phi+A(\Phi)\partial_x\Phi=0,
\end{align}
where
\begin{align*}
  \begin{split}
    A(\Phi)=\left(
              \begin{array}{cccccc}
                0 & 0  &-1 & 0  \\
                0 & 0 &0  & -1  \\
                -(c_1^2+2\sigma_0\phi_1) & -2\sigma_1\phi_2 &0  & 0  \\
                -2\sigma_1\phi_2 & -(c_2^2+2\sigma_1\phi_1)&0  & 0  \\
              \end{array}
            \right).
  \end{split}
\end{align*}
\begin{lem}\label{2dlem1}
 There exists a small $\delta>0$. For $\Phi\in B_{2\delta}^6(0)$, we have that the reduced system \eqref{2d2} of the two dimensional elastic wave equations under planar symmetry is uniformly strictly hyperbolic. Moreover, the $1^{\text{st}}$ and $4^{\text{th}}$ characteristics are genuinely non-linear:
 \begin{equation*}
   \nabla_\Phi\lambda_i \cdot r_i\neq 0,\quad i=1,4,\quad\forall\ \Phi\in B_{2\delta}^4.
 \end{equation*}
\end{lem}
{\it Proof.} Let $a,b,c$ be given by \eqref{y2.6}. We compute
\begin{align*}
  \begin{split}
      \det{(\lambda I-A)}&=\lambda^4-(a+b)\lambda^2+(ab-c^2)\\
      &=(\lambda^2-a)(\lambda^2-b)-c^2,
  \end{split}
\end{align*}
and obtain the eigenvalues of $(A(\Phi))_{4\times4}$:
\begin{align}\label{118}
  \begin{array}{lll}
    \lambda_1=&\sqrt{\frac12(a+b)+\frac12\sqrt{(a-b)^2+4c^2}},\\
    \lambda_2=&\sqrt{\frac12(a+b)-\frac12\sqrt{(a-b)^2+4c^2}},\\
    \lambda_3=&-\sqrt{\frac12(a+b)-\frac12\sqrt{(a-b)^2+4c^2}},\\
    \lambda_4=&-\sqrt{\frac12(a+b)+\frac12\sqrt{(a-b)^2+4c^2}}.
  \end{array}
\end{align}
These four eigenvalues are completely distinct,
\begin{align*}
    \lambda_4(\Phi)<\lambda_3(\Phi)<\lambda_2(\Phi)<\lambda_1(\Phi).
\end{align*}
Then we calculate their corresponding right eigenvectors  and get:
\begin{align*}
  r_1=\left(
                \begin{array}{c}
                  \frac{\lambda_1^2-b}{2\sigma_1} \\
                  \phi_2 \\
                  -\frac{\lambda_1(\lambda_1^2-b)}{2\sigma_1} \\
                  -\lambda_1\phi_2 \\
                \end{array}
              \right),\ r_2=\left(
                                     \begin{array}{c}
                                       \frac{\lambda_2^2-b}{2\sigma_1} \\
                                       \phi_2 \\
                                       -\frac{\lambda_2(\lambda_2^2-b)}{2\sigma_1} \\
                                       -\lambda_2\phi_2\\
                                     \end{array}
                                   \right),\ r_3=\left(
                                                   \begin{array}{c}
                                                     \frac{\lambda_2^2-b}{2\sigma_1} \\
                                                     \phi_2 \\
                                                     \frac{\lambda_2(\lambda_2^2-b)}{2\sigma_1} \\
                                                     \lambda_2\phi_2 \\
                                                   \end{array}
                                                 \right),\ r_4=\left(
                                                                 \begin{array}{c}
                                                                    \frac{\lambda_1^2-b}{2\sigma_1} \\
                                                                   \phi_2 \\
                                                                   \frac{\lambda_1(\lambda_1^2-b)}{2\sigma_1}\\
                                                                   \lambda_1\phi_2  \\
                                                                 \end{array}
                                                               \right).
\end{align*}
Since system \eqref{2d2} is strictly hyperbolic for $|\Phi|<2\delta$ with small $\delta$, Christodoulou-Perez's formulas and discussion of shock formation in \cite{christodoulou} are applicable. Here we follow our argument in Section \ref{yf}- Section \ref{ill} to study the ill-posedness theory of system \eqref{2d2}.

By calculation, we obtain
\begin{align*}
\begin{array}{lll}
  l^1=&\frac1{K}\Big(\frac{\lambda_1^2-b}{2\sigma_1},\phi_2,-\frac{\lambda_1^2-b}{2\sigma_1\lambda_{1}},-\frac{\phi_2}{\lambda_{1}}\Big),\
  &l^2=\frac1{N}\Big(\frac{\lambda_2^2-b}{2\sigma_1},\phi_2,-\frac{\lambda_2^2-b}{2\sigma_1\lambda_{2}},-\frac{\phi_2}{\lambda_{2}}\Big),\\ l^3=&\frac1{N}\Big(\frac{\lambda_2^2-b}{2\sigma_1},\phi_2,\frac{\lambda_2^2-b}{2\sigma_1\lambda_{2}},\frac{\phi_2}{\lambda_{2}}\Big)\ &l^4=\frac1{K}\Big(\frac{\lambda_1^2-b}{2\sigma_1},\phi_2,\frac{\lambda_1^2-b}{2\sigma_1\lambda_{1}},\frac{\phi_2}{\lambda_{1}}\Big),
\end{array}
\end{align*}
where
 \begin{align*}
   K=&\frac{(a-b)^2+4c^2}{4\sigma_1^2}+\frac{(a-b)\sqrt{(a-b)^2+4c^2}}{4\sigma_1^2},\\
   N=&\frac{(a-b)^2+4c^2}{4\sigma_1^2}-\frac{(a-b)\sqrt{(a-b)^2+4c^2}}{4\sigma_1^2}.
 \end{align*}
Moreover, since
\begin{equation*}
  (\lambda_1^2-b)(\lambda_2^2-b)+c^2=0,
\end{equation*}
we verify that
\begin{equation*}
  l^i(\Phi)\cdot r_j(\Phi)=\delta_j^i.
\end{equation*}

As discussed in Section \ref{yf}, we employ wave decomposition and explore the coefficients of equations \eqref{y14}-\eqref{y13}.
Denote
\begin{equation*}
  \Delta=(a-b)^2+4c^2.
\end{equation*}
Then, by calculation it follows
\begin{align*}
  &\nabla_\Phi\lambda_1=-\nabla_\Phi\lambda_4=\Big(\frac{(\sigma_0+\sigma_1)\sqrt{\Delta}+(a-b)(\sigma_0-\sigma_1)}{2\lambda_1\sqrt{\Delta}},\frac{4\sigma_1^2\phi_2}{\lambda_1\sqrt{\Delta}},0,0\Big),\\
  &\nabla_\Phi\lambda_2=-\nabla_\Phi\lambda_3=\Big(\frac{(\sigma_0+\sigma_1)\sqrt{\Delta}-(a-b)(\sigma_0-\sigma_1)}{2\lambda_2\sqrt{\Delta}},-\frac{4\sigma_1^2\phi_2}{\lambda_2\sqrt{\Delta}},0,0\Big).
\end{align*}
And from \eqref{y8} and \eqref{gamma}, we have
\begin{equation*}
  c_{11}^1(\Phi)=-c_{44}^4(\Phi)=\nabla_\Phi\lambda_1 \cdot r_1=\frac{2\sigma_0(a-b)(\lambda_1^2-b)+(2\sigma_0+6\sigma_1)c^2}{4\sigma_1\lambda_1\sqrt{\Delta}},
\end{equation*}
\begin{equation*}
  c_{22}^2(\Phi)=-c_{33}^3(\Phi)=\nabla_\Phi\lambda_2 \cdot r_2=-\frac{2\sigma_0(a-b)(\lambda_2^2-b)-(2\sigma_0+6\sigma_1)c^2}{4\sigma_1\lambda_2\sqrt{\Delta}}.
\end{equation*}
By \eqref{118} and definitions of $a,b,c$ in \eqref{y2.6}, for $\Phi\in B_{2\delta}^6(0)$, we have
\begin{equation*}
  a-b\sim c_1^2-c_2^2,\quad \lambda_1^2-b\sim c_1^2-c_2^2,\quad c\sim 0.
\end{equation*}
This implies
\begin{equation*}
  c_{11}^1(\Phi)=-c_{44}^4(\Phi)\neq 0.
\end{equation*}
This concludes the proof of Lemma \ref{2dlem1}.\hfill$\Box$

By the same definitions for $\rho_i$, $w^i$, $v^i$ as in \eqref{dense} and \eqref{y4} for $i=1,2,3,4$, we have decomposition of waves corresponding to this strictly hyperbolic system \eqref{2d2}:
\begin{align*}
  \partial_{s_i}w^i=&-c_{ii}^i(w^i)^2+\Big(\sum_{m\neq i}(-c_{im}^i+\gamma_{im}^i)w^m\Big)w^i+\sum_{m\neq i,k\neq i\atop m\neq k}\gamma_{km}^iw^kw^m,\\
  \partial_{s_i}\rho_i=&c_{ii}^iv^i+\Big(\sum_{m\neq i}c_{im}^iw^m\Big)\rho_i,\\
  \partial_{s_i}v^i=&\Big(\sum_{m\neq i}\gamma_{im}^iw^m\Big)v^i+\sum_{m\neq i,k\neq i\atop m\neq k}\gamma_{km}^iw^kw^m\rho_i.
\end{align*}

We construct the initial data satisfying
\begin{equation*}
  W_0^{(\eta)}:=\max_i\sup_{z_i}|w^i_{(\eta)}(z_i,0)|=w^1_{(\eta)}(z_0,0).
\end{equation*}
We further choose
\begin{equation*}
   w^1_{(\eta)}(z,0)=\theta\int_\mathbb{R}\zeta_{\frac\eta{10}}(y)| \ln (x-y)|^\alpha \chi(x-y)dy,\quad\text{for}\quad0<\alpha<1,
\end{equation*}
and require
\begin{align*}
  \max_{i=2,3,4}\sup_{z_i}|w^i_{(\eta)}(z_i,0)|\leq \frac{(1-\varepsilon)^4}{2(1+\varepsilon)^3}W_0^{(\eta)}.
  \end{align*}
Let \begin{align*}
  S(t):=&\max_{i}\sup_{(z'_i,s'_i)\atop z'_i\in[\eta,2\eta],\ 0\leq s'_i\leq t}\rho_i(z'_i,s'_i),& J(t):=&\max_{i}\sup_{(z'_i,s'_i)\atop z'_i\in[\eta,2\eta]\ 0\leq s'_i\leq t}|v^i(z'_i,s'_i)|,\\
  V(t):=&\max_i\sup_{(x',t')\notin\mathcal{R}_i,\atop 0\leq t'\leq t}|w^i(x',t')|,&\bar{U}(t):=&\sup_{(x',t')\atop 0\leq t'\leq t}|\Phi(x',t')|.
\end{align*}

For $\Phi\in C^2(\mathbb{R}\times[0,T],B_{2\delta}^4(0))$ being a solution to \eqref{111} for some $T>0$, with a general assumption $c_{11}^1(0)<0$, by analogous arguments as in Section \ref{3d}, we obtain
estimates:
\begin{align*}
  S&=O(1+tJ+tVS), &J&=O(W_0^{(\eta)}+t VJ+tV^2 S),\\
 V&=O\Big(\eta [W_0^{(\eta)}]^2+tV^2+\eta VJ\Big),&\bar{U}&=O(\eta J+\eta V+\eta tV),
\end{align*}
where $t<T$.
For $tW_0^{(\eta)}\leq C$, by a bootstrap argument we hence have the following bounds:
\begin{align*}
  S&=O(1), &J&=O(W_0^{(\eta)}),\\\
 V&=O\Big(\eta [W_0^{(\eta)}]^2\Big),&\bar{U}&=O(\eta W_0^{(\eta)}).
\end{align*}
These imply there exists $T_\eta^*$ satisfying
\begin{equation}
 \frac1{(1+\varepsilon)^3|c_{11}^1(0)|W_0^{(\eta)}}\leq T_\eta^*\leq\frac{1}{(1-\varepsilon)^4|c_{11}^1(0)|W_0^{(\eta)}}.
\end{equation}
And a shock forms at time $T_\eta^*$,
\begin{equation*}
  \lim_{t\rightarrow T_\eta^*} \rho_1(z_0,t)=0.
\end{equation*}
Moreover, deriving estimates as in Section \ref{supbound}, we can show $\partial_{z_1}\rho_1$ is bounded.

We construct $\Omega_{T_\eta^*}\subseteq\{(x,Y^2,t):t=T_\eta^*\}$ to be the region
\begin{equation*}
\frac{{c_2^2(x-c_1T_\eta^*-\frac{3\eta}{2})}^2}{c_1^2}+(Y^2)^2 \leq \frac{c_2^2}{c_1^2} {\left(\frac{\eta}{2} \right)}^2,
\end{equation*}
hence we obtain
\begin{equation*}
\int_{\Omega_{T_\eta^*}^r} dY^2 \sim(1+O(\varepsilon))\sqrt{\big(z-\eta+O(\varepsilon)\eta\big)\big(2\eta-z+O(\varepsilon)\eta\big)},
\end{equation*}
where $\Omega_{T_\eta^*}^r=\{Y^2:(x,Y^2)\in\Omega_{T_\eta^*}\}$.
This yields
\begin{equation*}
\begin{split}
  &\|w^1(\cdot,T_\eta^*)\|_{L^2(\Omega_{T_\eta^*})}^2\\
  \geq &C\int_{\eta}^{2\eta}\Big|\frac{v^1}{\rho_1}(z,T_\eta^*)\Big|^2\rho_1(z,T_\eta^*) \sqrt{\big(z-\eta+O(\varepsilon)\eta\big)\big(2\eta-z+O(\varepsilon)\eta\big)}dz\\
\geq &C{(1-\varepsilon)}^2{(1-2\varepsilon)}^2[W_0^{(\eta)}]^2\int_{z_0}^{z_0^*}\frac{1}{\rho_1(z,T_\eta^*)} \sqrt{\big(z-\eta+O(\varepsilon)\eta\big)\big(2\eta-z+O(\varepsilon)\eta\big)}dz\\
\geq& C\sqrt{(z_0-\eta)(2\eta-z_0^*)}[W_0^{(\eta)}]^2\int_{z_0}^{z_0^*}\frac{1}{\rho_1(z,T_\eta^*)-\rho_1(z_0,T_\eta^*)}dz\\
\geq& C\sqrt{(z_0-\eta)(2\eta-z_0^*)}[W_0^{(\eta)}]^2\int_{z_0}^{z_0^*}\frac{1}{(\sup_{z\in(z_0,z_0^*]}|\partial_z\rho_1|)(z-z_0)}dz.
  \end{split}
\end{equation*}
With the boundedness of $\partial_{z_1}\rho_1$, we arrive at
\begin{equation*}
  \|U_\eta(\cdot,T_\eta^*)\|_{H^2(\Omega_{T_\eta^*})}\gtrsim\|w^1(\cdot,T_\eta^*)\|_{L^2(\Omega_{T_\eta^*})}=+\infty.
\end{equation*}

In summary, we prove the following result.
\begin{thm}
The Cauchy problems of 2D elastic wave equations are ill-posed in $H^\frac52(\mathbb{R}^2)$ in the following sense: We construct a family of compactly supported smooth initial data $(U_0^{(\eta)}, U_1^{(\eta)})$ with $${\|U_0^{(\eta)}\|}_{{\dot{H}}^\frac52(\mathbb{R}^2)}+{\|U_1^{(\eta)}\|}_{{\dot{H}}^\frac32(\mathbb{R}^2)} \to 0,\quad \text{as} \quad \eta\to0.$$ Let $T_\eta^*$ be the largest time such that \eqref{y1} has a solution $U_\eta\in C^\infty(\mathbb{R}^2\times[0,T_\eta^*))$. As $\eta \to 0$, we have $T_\eta^* \to 0$.

Moreover, in a spatial region $\Omega_{T_\eta^*}$ the $H^2$ norm of the solution to elastic waves \eqref{y1} blows up at shock formation time $T^*_{\eta}$:
\begin{equation*}
  \|U_\eta(\cdot,T_\eta^*)\|_{H^2(\Omega_{T_\eta^*})}=+\infty.
\end{equation*}
\end{thm}

\appendix
\section{Appendix: Norms for initial data}
In this appendix, we construct examples of initial data in 3D and 2D satisfying requirements in Section \ref{datanorm}  and in Section \ref{sec2d}.

Assume $\eta>0$ and $\theta>0$. Let $\chi(x)$ be the characteristic function
\begin{equation}
  \chi(x)=\left\{\begin{array}{ll}
  1,\quad x\in [\frac65\eta,\frac{9}5\eta],\\
  0,\quad x\notin[\frac65\eta,\frac95\eta].
  \end{array}\right.
\end{equation}
Set $\zeta_{\frac\eta{10}}(x)$ to be a test function in $C_c^\infty(\mathbb{R})$ satisfying:
\begin{equation*}
  \text{supp}\ \zeta_{\frac\eta{10}}(x)\subseteq\{x:|x|\leq\frac\eta{10}\},\,\,\, 0\leq\zeta_{\frac\eta{10}}(x)\leq\frac{20}{\eta},  \mbox{  and  }  \int_{\mathbb{R}}\zeta_{\frac\eta{10}}(x)dx=1.
\end{equation*}
We start with the 3D case.
Denote $B_{\frac\eta2}^3:=\{(x,Y^2,Y^3):(x-\frac{3\eta}{2})^2+(Y^2)^2+(Y^3)^2\leq(\frac\eta2)^2\}$. We have:
\begin{lem}\label{alem1}
Let $0<\alpha<\frac12$. We can construt a function $\hat{w}\in H^1(\mathbb{R}^3)$ satisfying
\begin{equation}\label{app1}
   \hat{w}(x,Y^2,Y^3)=\theta\int_\mathbb{R}\zeta_{\frac\eta{10}}(y)| \ln (x-y)|^\alpha \chi(x-y)dy, \ \text{in}\ B_{\frac\eta2}^3,
\end{equation}
and
\begin{equation}\label{a44}
  \|\hat{w}\|_{\dot{H}^1(\mathbb{R}^3)}\lesssim\theta\sqrt{\eta}(1+|\ln \frac{6\eta}{5}|^{\alpha}+|\ln \frac{9\eta}{5}|^{\alpha}\big).
\end{equation}
\end{lem}
{\it Proof.}
By construction, we have $\text{supp}\ \hat{w}\subset\{(x,Y^2,Y^3):\eta\leq x\leq 2\eta\}$. For $\Omega_0^r=\{(Y^2,Y^3):(x,Y^2,Y^3)\in B_{\frac\eta2}^3\}$, we have
\begin{equation}\label{a4}
\begin{split}
  \int_{\Omega_0^r}dY^2dY^3=&(\frac\eta2)^2-(x-\frac{3\eta}{2})^2=(x-\eta)(2\eta-x)\\
  =&\big[(x-\frac\eta{10})-\frac{9\eta}{10}\big]\big[\frac{19\eta}{10}-(x-\frac{\eta}{10})\big]\\
  \leq&\big[(x-\frac\eta{10})-\frac{9\eta}{10}\big]\frac{19\eta}{10}\\
  \leq&\frac{19}{10}\eta (x-\frac\eta{10}).
\end{split}
\end{equation}
It follows that
\begin{equation*}
\begin{split}
  &\int_{B_{\frac\eta2}^3}|\partial_x\hat{w}|^2dxdY^2dY^3\\
  \leq&\underbrace{\int_{B_{\frac\eta2}^3}\Big|\int_{|y|\leq\frac{\eta}{10}} \zeta_{\frac\eta{10}}(y)\frac{\alpha\theta}{(x-y)|\ln (x-y)|^{1-\alpha}}\chi(x-y)dy\Big|^2dxdY^2dY^3}_{L_1}\\
  &+\underbrace{\int_{B_{\frac\eta2}^3}\Big|\int_{|y|\leq\frac{\eta}{10}} \zeta_{\frac\eta{10}}(y)\theta|\ln (x-y)|^{\alpha}\delta(x-y-\frac{6\eta}{5})dy\Big|^2dxdY^2dY^3}_{L_2}\\
    &+\underbrace{\int_{B_{\frac\eta2}^3}\Big|\int_{|y|\leq\frac{\eta}{10}} \zeta_{\frac\eta{10}}(y)\theta|\ln (x-y)|^{\alpha}\delta(x-y-\frac{9\eta}{5})dy\Big|^2dxdY^2dY^3}_{L_3}
\end{split}\end{equation*}
For $L_1$, by \eqref{a4}, we have
\begin{equation}\label{a5}
\begin{split}
 L_1\lesssim&\int_\eta^{2\eta}\frac{\alpha^2\theta^2\cdot(x-\eta)(2\eta-x)}{(x-\frac\eta{10})^2|\ln (x-\frac\eta{10})|^{2-2\alpha}}dx\\
  \lesssim&\int_\eta^{2\eta}\frac{\alpha^2\theta^2\eta}{(x-\frac\eta{10})|\ln (x-\frac\eta{10})|^{2-2\alpha}}dx\\
  \lesssim&\theta^2\eta.
\end{split}\end{equation}
For $L_2+L_3$, we proceed and get
\begin{equation}\label{a55}
\begin{split}
 L_2+L_3\lesssim&\int_{B_{\frac\eta2}^3}\Big|\theta\sup_z|\zeta(z)|\int_{\mathbb{R}}  |\ln(x-y)|^{\alpha}\delta(x-y-\frac{6\eta}{5})dy\Big|^2dxdY^2dY^3\\
 &+\int_{B_{\frac\eta2}^3}\Big|\theta\sup_z|\zeta(z)|\int_{\mathbb{R}} |\ln(x-y)|^{\alpha}\delta(x-y-\frac{9\eta}{5})dy\Big|^2dxdY^2dY^3\\
  \lesssim&\theta^2\frac{1}{\eta^2}\cdot\big(|\ln \frac{6\eta}{5}|^{2\alpha}+|\ln \frac{9\eta}{5}|^{2\alpha}\big)\cdot\eta^3\lesssim\theta^2\eta\cdot\big(|\ln \frac{6\eta}{5}|^{2\alpha}+|\ln \frac{9\eta}{5}|^{2\alpha}\big).
\end{split}\end{equation}
Together with \eqref{a5} and \eqref{a55}, we conclude
\begin{equation*}
\begin{split}
  \int_{B_{\frac\eta2}^3}|\partial_x\hat{w}|^2dxdY^2dY^3\lesssim\theta^2\eta\big(1+|\ln \frac{6\eta}{5}|^{2\alpha}+|\ln \frac{9\eta}{5}|^{2\alpha}\big).
\end{split}\end{equation*}

We hence have $\|\hat{w}\|_{\dot{H}^1(B_{\eta/2}^3)}\lesssim\theta\sqrt{\eta}\big(1+|\ln \frac{6\eta}{5}|^{\alpha}+|\ln \frac{9\eta}{5}|^{\alpha}\big)$. We then extend $\hat{w}$ to a compactly supported $H^1$ function in the whole region satisfying  \eqref{a44}. This finishes the proof.\hfill$\Box$

For 2D case, we will use the following identity: for $0<s<1$, it holds that
\begin{equation}\label{a6}
  \|f\|_{\dot{H}^s(\mathbb{R}^d)}^2=C_s\iint_{\mathbb{R}^d\times\mathbb{R}^d}\frac{[f(x+y)-f(x)]^2}{|y|^{2s+d}}dxdy.
\end{equation}
Denote $B_{\frac\eta2}^2:=\{(x,Y^2):(x-\frac{3\eta}{2})^2+(Y^2)^2\leq(\frac\eta2)^2\}$. We have:
\begin{lem}
Let $0<\alpha<1$. We can construct a function $\hat{w}\in H^{\frac12}(\mathbb{R}^2)$ satisfying
\begin{equation}\label{app2}
   \hat{w}(x,Y^2)=\theta\int_\mathbb{R}\zeta_{\frac\eta{10}}(y)| \ln (x-y)|^\alpha \chi(x-y)dy, \ \text{in}\ B_{\frac\eta2}^2,
\end{equation}
and
\begin{equation}\label{com}
  \|\hat{w}\|_{\dot{H}^{\frac12}(\mathbb{R}^2)}\lesssim\theta\eta^{\frac14}\big(1+|\ln \frac{6\eta}{5}|^{\alpha}+|\ln \frac{9\eta}{5}|^{\alpha}\big).
\end{equation}
\end{lem}
{\it Proof.}
By construction, we have $\text{supp}\ \hat{w}\subset\{(x,Y^2):\eta\leq x\leq 2\eta\}$. For $\Omega_1^r=\{Y^2:(x,Y^2)\in B_{\frac\eta2}^2\}$, it holds that
\begin{equation*}
  \int_{\Omega_1^r}dY^2=\sqrt{(x-\eta)(2\eta-x)}\leq\frac12\eta.
\end{equation*}
Let $w(x)=\theta (\zeta_{\frac\eta{10}}\ast(| \ln |\cdot||^\alpha \chi))(x)$. We have $\text{supp}\ w\subseteq[\eta,2\eta]$ and
\begin{equation*}
  \begin{split}
    \int_{B_{\frac\eta2}^2}\big||\nabla_x|^{\frac12}\hat{w}\big|^2dxdY^2\leq&\frac12\eta\int_\eta^{2\eta}\big||\nabla_x|^{\frac12}w\big|^2dx.
  \end{split}
\end{equation*}
Applying \eqref{a6} for $d=1$, we get
\begin{equation}\label{a11}
  \begin{split}
    \int_{B_{\frac\eta2}^2}\big||\nabla_x|^{\frac12}\hat{w}\big|^2dxdY^2\lesssim\eta\iint_{\mathbb{R}\times\mathbb{R}}\frac{|w(x+y)-w(x)|^2}{|y|^2}dxdy.
  \end{split}
\end{equation}
Now we split the domain of integration on the right hand side of \eqref{a11} into six parts:
\begin{align*}
  &I_1=\{(x,y):x\in[\eta,2\eta],y\in[-\sqrt{\eta},\sqrt{\eta}]\},\\
  &I_2=\{(x,y):x\in[\eta,2\eta],y\in(\sqrt{\eta},+\infty)\},\\
  &I_3=\{(x,y):x\in[\eta,2\eta],y\in(-\infty,-\sqrt{\eta})\},\\
  &I_4=\{(x,y):x\in(2\eta,+\infty),x+y\in[\eta,2\eta]\}\\
   &I_5=\{(x,y):x\in(-\infty,\eta),x+y\in[\eta,2\eta]\},\\
   &I_6=\mathbb{R}^2\setminus\cup_{k\in\{1,\cdots,5\}}I_k.
\end{align*}
In this way, $w(x+y)-w(x)$ is supported in $I_j, j=1,\cdots,5$. P

In $I_1$, we have
\begin{equation}\label{a7}
\begin{split}
\eta\iint_{I_1}\frac{|w(x+y)-w(x)|^2}{|y|^2}dxdy\leq&\eta\int_{-\sqrt{\eta}}^{\sqrt{\eta}}\int_\eta^{2\eta}\frac{|\int_0^y\partial_xw(x+z)dz|^2}{|y|^2}dxdy\\
    \leq&\eta\int_{-\sqrt{\eta}}^{\sqrt{\eta}}\int_\eta^{2\eta}\sup_{z\in[\eta,2\eta]}|\partial_zw(z)|^2dxdy.
\end{split}
\end{equation}
Since
\begin{equation*}
  \begin{split}
   \sup_{z\in[\eta,2\eta]}| \partial_zw(z)|\leq&\sup_{z\in[\eta,2\eta]}\int_{|y|\leq\frac{\eta}{10}} \zeta_{\frac\eta{10}}(y)\frac{\alpha\theta}{(z-y)|\ln (z-y)|^{1-\alpha}}dy\\
   &+\sup_{z\in[\eta,2\eta]}\int_{|y|\leq\frac{\eta}{10}} \zeta_{\frac\eta{10}}(y)\alpha\theta|\ln (z-y)|^{\alpha}\delta(z-y-\frac{6\eta}{5})dy\\
   &+\sup_{z\in[\eta,2\eta]}\int_{|y|\leq\frac{\eta}{10}} \zeta_{\frac\eta{10}}(y)\alpha\theta|\ln (z-y)|^{\alpha}\delta(z-y-\frac{9\eta}{5})dy\\
   \lesssim&\sup_{z\in[\eta,2\eta]}\frac{\theta}{(z-\frac\eta{10})|\ln (z-\frac\eta{10})|^{1-\alpha}}+\frac{1}{\eta}\theta\cdot\big(|\ln \frac{6\eta}{5}|^{\alpha}+|\ln \frac{9\eta}{5}|^{\alpha}\big)\\
   \lesssim&\frac{\theta}{\eta}\big(1+|\ln \frac{6\eta}{5}|^{\alpha}+|\ln \frac{9\eta}{5}|^{\alpha}\big),
  \end{split}
\end{equation*}
from \eqref{a7}, we get
\begin{equation}\label{aa}
  \begin{split}
   \eta\iint_{I_1}\frac{|w(x+y)-w(x)|^2}{|y|^2}dxdy\lesssim \theta^2\sqrt{\eta}(1+|\ln \frac{6\eta}{5}|^{2\alpha}+|\ln \frac{9\eta}{5}|^{2\alpha}\big).
  \end{split}
\end{equation}
In $I_2$, we have $w(x+y)=0$. Thus, it holds that
\begin{equation}\label{a8}
\begin{split}
&\eta\iint_{I_2}\frac{|w(x+y)-w(x)|^2}{|y|^2}dxdy
=\eta\int_{\sqrt{\eta}}^{+\infty}\int_{\eta}^{2\eta}\frac{|w(x)|^2}{|y|^2}dxdy\\
    \lesssim&\eta\big|\sup_{x\in[\eta,2\eta]}w(x)\big|^2\int_{\sqrt{\eta}}^{+\infty}\int_{\eta}^{2\eta}\frac{1}{|y|^2}dxdy\lesssim\eta\theta^2|\ln\eta|^{2\alpha}\frac{\eta}{\sqrt{\eta}}\\
    \lesssim&\eta^{\frac32}\theta^2|\ln\eta|^{2\alpha}.
\end{split}
\end{equation}
By the same calculation as in \eqref{a8}, we obtain
\begin{equation}\label{a10}
\begin{split}
\eta\iint_{I_3}\frac{|w(x+y)-w(x)|^2}{|y|^2}dxdy\lesssim&\eta^{\frac32}\theta^2|\ln\eta|^{2\alpha}.
\end{split}
\end{equation}
In $I_4$, it holds that $w(x)=0$. Note $\text{supp}\ w\subseteq[\eta,2\eta]$. We have
\begin{equation}\label{a9}
\begin{split}
&\eta\iint_{I_4}\frac{|w(x+y)-w(x)|^2}{|y|^2}dxdy\\
=&\eta\int_{2\eta}^{+\infty}\int_{-\infty}^{-\eta}\frac{|w(x+y)|^2}{|y|^2}dydx+\eta\int_{2\eta}^{+\infty}\int_{-\eta}^{0}\frac{|w(x+y)-w(x)|^2}{|y|^2}dydx\\
\leq&\eta\int_{\eta}^{2\eta}\int_{-\infty}^{-\eta}\frac{|w(z)|^2}{|y|^2}dydz+\eta\int_{2\eta}^{+\infty}\int_{-\eta}^{0}\frac{|\int_0^y\partial_xw(x+z)dz|^2}{|y|^2}dydx\\
\lesssim&\eta\sup_{z\in[\eta,2\eta]}|w(z)|^2\int_{\eta}^{2\eta}\int_{-\infty}^{-\eta}\frac{1}{|y|^2}dydz+\eta\sup_{z\in[\eta,2\eta]}|\partial_zw(z)|^2\eta\int_{-\eta}^{0}dy\\
\lesssim&\eta\theta^2|\ln\eta|^{2\alpha}+\eta^3\frac{\theta^2}{\eta^2}\big(1+|\ln \frac{6\eta}{5}|^{2\alpha}+|\ln \frac{9\eta}{5}|^{2\alpha}\big)\\
\lesssim&\eta\theta^2\big(1+|\ln \frac{6\eta}{5}|^{2\alpha}+|\ln \frac{9\eta}{5}|^{2\alpha}\big).
\end{split}
\end{equation}
For $I_5$, it holds that $y\in[0,+\infty)$. We divide $[0,+\infty)$ into $[\eta,+\infty)$ and $[0,\eta)$. With a similar calculation as in \eqref{a9}, we also get
\begin{equation}\label{a13}
\begin{split}
\eta\iint_{I_5}\frac{|w(x+y)-w(x)|^2}{|y|^2}dxdy
\lesssim\eta\theta^2\big(1+|\ln \frac{6\eta}{5}|^{2\alpha}+|\ln \frac{9\eta}{5}|^{2\alpha}\big).
\end{split}
\end{equation}
Finally, in $I_6$, since $w(x+y)=w(x)=0$, we have
\begin{equation}\label{a14}
\begin{split}
\eta\iint_{I_6}\frac{|w(x+y)-w(x)|^2}{|y|^2}dxdy=0.
\end{split}
\end{equation}
Consequently, together with \eqref{a11}, \eqref{aa}, \eqref{a8}, \eqref{a10}, \eqref{a9}, \eqref{a13} and \eqref{a14}, we obtain
\begin{equation*}
  \begin{split}
    \int_{B_{\frac\eta2}^2}\big||\nabla_x|^{\frac12}\hat{w}\big|^2dxdY^2\lesssim \theta^2\sqrt{\eta}\big(1+|\ln \frac{6\eta}{5}|^{2\alpha}+|\ln \frac{9\eta}{5}|^{2\alpha}\big).
  \end{split}
\end{equation*}
We hence have $\|\hat{w}\|_{\dot{H}^{\frac12}(B_{\eta/2}^2)}\lesssim\theta\eta^{\frac14}\big(1+|\ln \frac{6\eta}{5}|^{\alpha}+|\ln \frac{9\eta}{5}|^{\alpha}\big)$. We then extend $\hat{w}$ to a compactly supported $H^\frac12$ function in the whole region satisfying  \eqref{com}. This completes the proof.\hfill$\Box$


\begin{thebibliography}{99}
\bibitem{Agemi}
Rentaro Agemi, {\it Global existence of nonlinear elastic waves}, Invent. Math. {\bf 142}(2), 225-250(2000).
\bibitem{Alinhac99} Serge Alinhac, \textit{Blowup of small data solutions for a quasilinear wave equation in two space dimensions}, Ann. of Math. {\bf 149}(1), 97-127(1999).
 \bibitem{Alinhac99II} Serge Alinhac, \textit{Blowup of small data solutions for a class of quasilinear wave equations in two space dimensions, \uppercase\expandafter{\romannumeral2}}, Acta. Math. {\bf 182}(1), 1-23(1999).
\bibitem{Alinhac01} Serge Alinhac, \textit{The null condition for quasilinear wave equations in two space dimensions \uppercase\expandafter{\romannumeral1}}, Invent. Math. {\bf 145}(3), 597-618(2001).
\bibitem{bahouri-chemin1} Hajer Bahouri, Jean-Yves Chemin, \textit{\'{E}quations d'ondes quasilin\'{e}aires et effet dispersif}, Internat. Math. Res. Notices. no. 21, 1141-1178(1999).
\bibitem{bahouri-chemin2} Hajer Bahouri, Jean-Yves Chemin, \textit{\'{E}quations d'ondes quasilin\'{e}aires et estimations de Strichartz}, Amer. J. Math. {\bf 121}(6), 1337-1377(1999).
\bibitem{christodoulou10} Demetrios Christodoulou, {\it The Formation of Shocks in 3-dimensional Fluids, EMS Monographs in Mathematics}, Z\"{u}rich: European Mathematical Society Publishing House, viii, 992 pp(2007).
\bibitem{christodoulou} Demetrios Christodoulou, Daniel Raoul Perez, {\it On the formation of shocks of electromagnetic plane waves in non-linear crystals}, J. Math. Phys. 57, 081506(2016).
\bibitem{Disconzi}  Marcelo M Disconzi, Chenyun Luo, Guizy Mazzone, Jared Speck, \textit{Rough sound waves in 3D compressible Euler flow with vorticity}, arXiv:1909.02550.

\bibitem{lindblad17} Boris Ettinger, Hans Lindblad, \textit{A Sharp Counterexample to Local Existence of Low Regularity Solutions to Einstein's Equations in Wave  Coordinates}, Ann. of Math. {\bf 185}(1), 311-330(2017).
\bibitem{Ross}  Ross Granowski, \textit{Asymptotically Stable Ill-posedness of Geometric Quasilinear Wave Equations}, Thesis, 132 pp(2018).


\bibitem{zha}
    Kunio Hidano, Dongbing Zha, {\it Space-time $L^2$ estimates, regularity and almost global existence for elastic waves}, Forum Math. {\bf 30}(5), 1291-1307(2018).
\bibitem{kato}
 Thomas J R Hughes, Tosio Kato, Jerrold E Marsden, {\it Well-posed quasi-linear second-order hyperbolic systems with applications to nonlinear elastodynamics and general relativity}, Arch. Rational Mech. Anal. {\bf 63}(3), 273-294(1976).

\bibitem{john74} Fritz John, \textit{Formation of singularities in one-dimensional nonlinear wave propagation}, Comm. Pure. Appl. Math. {\bf 27}, 377-405(1974).
\bibitem{john84} Fritz John, \textit{Formation of singularities in elastic waves}, Lec. Notes in Phys. {\bf 195}, Springer, 194-210(1984).
\bibitem{john88} Fritz John, \textit{Almost global existence of elastic waves of finite amplitude arising from small initial disturbances}, Comm. Pure. Appl. Math. {\bf 41}, 615-666(1988).
\bibitem{klainerman-igor}  Sergiu Klainerman, Igor Rodnianski, \textit{Improved local well-posedness for quasilinear wave equations in dimension three}, Duke. Math. J. {\bf 117}(1), 1-124(2003).
 \bibitem{klainerman-igor05}  Sergiu Klainerman, Igor Rodnianski, \textit{Rough solutions of the Einstein-vacuum equations}, Ann. of Math. {\bf 161}(3), 1195-1243(2005).
 \bibitem{klainerman-igor-szeftel}  Sergiu Klainerman, Igor Rodnianski, Jeremie Szeftel, \textit{The bounded $L^2$ curvature conjecture}, Invent. Math. {\bf 202}(1), 91-216(2015).
\bibitem{klainerman}  Sergiu Klainerman, Thomas C Sideris, \textit{On almost global existence for nonrelativistic wave equations in 3D}, Comm. Pure. Appl. Math. {\bf 49}, 307-321(1996).
\bibitem{lindblad93} Hans Lindblad, \textit{A Sharp Counterexample to Local Existence of Low-Regularity Solutions to Nonlinear Wave Equations}, Duke. Math. J {\bf 72}(2), 43-110(1993).
\bibitem{Lind96}
Hans Lindblad, {\it Conterexamples to local existence for semi-linear wave equations}, Amer. J. Math. {\bf 118}(1), 1-16(1996).
\bibitem{Lind98}
Hans Lindblad, {\it Conterexamples to local existence for quasilinear wave equations}, Math. Research Letters 5, 605-622(1998).

\bibitem{Speck-luk}
Jonathan Luk, Jared Speck, \textit{Shock Formation in solutions to the 2D compressible Euler equations in the presence of non-zero vorticity},  Invent. Math. {\bf 214}(1), 1-169(2018).

\bibitem{miao}Shuang Miao, Pin Yu, \textit{On the fromation of shocks for quasilinear wave equations}, Invent. Math. {\bf 207}(2), 697-831(2017).

\bibitem{Sideris96}
Thomas C Sideris, {\it The null condition and global existence of nonlinear elastic waves}, Invent. Math. {\bf 123}, 323-342(1996).
\bibitem{Sideris00}
Thomas C Sideris, {\it Nonresonance and global existence of prestressed nonlinear elastic waves}, Ann. of  Math. {\bf 151}(2), 849-874(2000).
\bibitem{speckbk}
Jared Speck, {\it Shock formation in small-data solutions to 3D quasilinear wave equations}, Math. Surveys and Monographs 214, Amer. Math. Soc., Providence, RI, 2016.
\bibitem{Speck18}
Jared Speck, \textit{Shock Formation for 2D quasilinear wave systems featuring multiple speeds: Blowup for the fastest wave, with non-trivial interactions up to the singularity},  Ann. PDE, {\bf 4}(1), 131 pp(2018).
\bibitem{speck20}
Jared Speck, {\it Stable ODE-type blowup for some quasilinear wave equations with derivative-quadratic nonlinearities}, Anal. PDE, {\bf 13}(1), 93šC146(2020).
\bibitem{Speck16}
Jared Speck, Gustav Holzegel, Jonathan Luk, Willie Wong, \textit{Stable Shock Formation for Nearly Simple Outgoing Plane Symmetry Waves},  Ann. PDE, {\bf 2}(2), 198 pp(2016).


\bibitem{S1}
Jeremie Szeftel, {\it Parametrix for wave equations on a rough background \uppercase\expandafter{\romannumeral1}: regularity of the phase
at initial time}, arXiv:1204.1768(2012).
\bibitem{S2}
Jeremie Szeftel, {\it Parametrix for wave equations on a rough background \uppercase\expandafter{\romannumeral2}: construction of the
parametrix and control at initial time}, arXiv:1204.1769(2012).
\bibitem{S3}
Jeremie Szeftel, {\it Parametrix for wave equations on a rough background \uppercase\expandafter{\romannumeral3}: space-time regularity
of the phase}, arXiv:1204.1770(2012).
\bibitem{S4}
Jeremie Szeftel, {\it Parametrix for wave equations on a rough background \uppercase\expandafter{\romannumeral4}: control of the error
term}, arXiv:1204.1771(2012).


\bibitem{shadi}
Shadi A Tahvildar-Zadeh, {\it Relativistic and nonrelativistic elastodynamics with small shear strains},
Annales de Inst. I'I.H.P. Physique th\'eorique. {\bf 69}(3), 275-C307 (1998).

\bibitem{tataru1}
Daniel Tataru, \textit{Strichartz estimates for operators with nonsmooth coefficients and the nonlinear wave equation},  Amer. J. Math. {\bf 122}(2), 349-376(2000).
\bibitem{tataru2}
Daniel Tataru, \textit{Strichartz estimates for second order hyperbolic operators with nonsmooth coefficients, II},  Amer. J. Math. {\bf 123}(3), 385-423(2001).
\bibitem{tataru3}
Daniel Tataru, \textit{Strichartz estimates for second order hyperbolic operators with nonsmooth coefficients, III},  J. Amer. Math. Soc. {\bf 15}(2), 419-442(2002).
\bibitem{tataru}
Daniel Tataru, Hart F Smith, \textit{Sharp local well-posedness results for the nonlinear wave equation},  Ann. of Math. {\bf 162}(1), 291-366(2005).
\bibitem{Wang}
Qian Wang, \textit{A geometric approach for sharp local well-posedness of quasilinear wave equations},  Ann. PDE, {\bf 3}(1), 108 pp(2016).
\bibitem{Wang19}
Qian Wang, {\it Rough solutions of the 3-D compressible Euler equations}, arXiv: 1911.05038.
\bibitem{Zhou}
Yi Zhou, {\it Low regularity solutions for linearly degenerate hyperbolic systems}, Nonlinear Anal. {\bf 26}(11), 1843-1857(1996).
\end{thebibliography}
\end{document}